\documentclass[12pt,leqno]{article}

\usepackage{amsthm}

\usepackage{amsmath}
\usepackage{amscd}
\usepackage{amsmath}
\usepackage{amssymb}
\usepackage{latexsym}
\makeatletter

\@addtoreset{equation}{section}
\makeatother

\newtheorem{thm}{Theorem}[section]
\newtheorem{pr}[thm]{Proposition}
\newtheorem{df}[thm]{Definition}
\newtheorem{lm}[thm]{Lemma}
\newtheorem{cor}[thm]{Corollary}

\newtheorem{rmk}[thm]{Remark}

\newtheorem{df-lm}[thm]{Definition-Lemma}

\newcommand{\A}{\mathbb{A}}
\newcommand{\OO}{\mathcal{O}}
\newcommand{\F}{\mathbb{F}}
\newcommand{\Proj}{\mathbb{P}}

\newcommand{\Ql}{\mathbb{Q}_\ell}

\newcommand{\mQ}{\mathbb{Q}}

\input xy \xyoption{all} \CompileMatrices

\begin{document}

\title{Symmetric bilinear forms and local epsilon factors of isolated 
singularities in positive characteristic}
\author{Daichi Takeuchi}
\date{}
\maketitle

\footnote[0]{\textit{Keywords}: Isolated singularity; Vanishing cycles; Local epsilon factors}
\footnote[0]{2020 \textit{Mathematics Subject Classification}. 
 Primary: 14B05, 14G17; Secondary: 11S40.}
\begin{abstract}
Let $f\colon X\to\mathbb{A}^1_k$ be a morphism from a smooth variety to an affine line with an isolated singular point. For such a 
singularity, 
we have two invariants. One is a 
non-degenerate symmetric bilinear form (de Rham), and the other 
is the vanishing cycles complex (\'etale). 

In this article, we give a formula which expresses the local epsilon 
factor of the vanishing cycles complex in terms of the bilinear form. In particular, 
the sign of the local epsilon factor is determined by the discriminant 
of the bilinear form. This formula can be thought as 
a refinement of the Milnor formula, which compares 
the total dimension 
of the vanishing cycles and 
the rank of the bilinear form. 

In characteristic $2$, we find a generalization of the Arf invariant, which can be regarded as an invariant for non-degenerate 
quadratic singularities, to general isolated singularities. 

\end{abstract}
\section{Introduction}

Let $k$ be a perfect field of characteristic $p>0$ and $X$ be a smooth 
$k$-scheme. Let $f\colon X\to C$ be a morphism of $k$-schemes to a smooth 
$k$-curve with an isolated singular point. 
To such a singularity, 
one can attach an arithmetic invariant, namely {\it the vanishing cycles complex}.

Fix a prime number $\ell$ different from $p$. Let $x\in X$ be an isolated singular point with respect to $f$. 
Then, the vanishing 
cycles complex $R\Phi_f(\Ql)_x$ supported at $x$ is 
a bounded complex of finite dimensional $\Ql$-representations of the 
absolute Galois group of the local field of $C$ at $f(x)\in C$. It is expected and observed in many cases that 
the Galois action on $R\Phi_f(\Ql)_x$ respects some complexity of the singularity; as the simplest 
example, it is acyclic if $x$ is a smooth point. Aside from such special cases, 
it is usually difficult to describe $R\Phi_f(\mathbb{Q}_\ell)_x$ explicitly as Galois representations. On the other hand, 
since singular points are of 
geometric nature, it is natural to ask whether there are relations between $R\Phi_f(\Ql)_x$ and 
some invariants of singularities which arise geometrically.

In this direction, Deligne considers an analogue of the Milnor formula in 
the complex geometry to an arithmetic setting in SGA $7$ \cite{Mil}. There he considers the total dimension 
${\rm dimtot}R\Phi_f(\Ql)_x$ and proves in positive characteristic cases that it coincides with the Milnor number
\begin{equation}\label{intromil}
(-1)^{n+1}{\rm dimtot}R\Phi_f(\Ql)_x=\mu(f,x). 
\end{equation}
Here $n$ is the dimension of $X$ and $\mu(f,x)$ is the Milnor number, which is defined purely algebro-geometrically from the isolated singular point. 

In this article, we give a formula which 
expresses {\it the local epsilon factors}, instead of the total dimensions, 
of the vanishing cycles complexes in terms of coherent sheaves (more precisely, symmetric bilinear forms 
as described below) arising from singularities. 
The local epsilon factor is defined by Langlands and Deligne to be a candidate of Galois-theoretic counterparts of the constant terms in the 
functional equations of automorphic local L-functions, in view of the local Langlands correspondence. 
As is well-known, and is what follows from our main results, the total dimensions of the vanishing cycles complexes has 
as much information as the (any) complex absolute values of 
the local epsilon factors. In this sense, our formula for local epsilon factors can be regarded as a refinement of the Milnor formula (\ref{intromil}). 

To motivate us to give such a refinement, let us consider the Hasse-Weil zeta function $Z(X/\F_q;t)$ of an $n$-dimensional 
projective smooth variety $X$ over a finite field $\mathbb{F}_q$. It is defined to be the infinite product 
\begin{equation*}
Z(X/\mathbb{F}_q;t)=\prod_{x\in|X|}\frac{1}{1-t^{{\rm deg}(x/\F_q)}}
\end{equation*}
indexed by the closed points of $X$. Similarly as the usual zeta function and L-functions, it admits a functional equation of the following type 
\begin{equation*}
Z(X/\F_q;t)=\varepsilon(X,\mathbb{Q}_\ell)t^{-\chi(X)}Z(X/\F_q;\frac{1}{q^nt}). 
\end{equation*}
Here $\chi(X)$ is the Euler characteristic, say, of the $\ell$-adic cohomology groups $H^i(X_{\overline{\F_q}},\Ql)$ of $X$ and $\varepsilon(X,\mathbb{Q}_\ell)=\prod_i
\det(-{\rm Frob}_q,H^i(X_{\overline{\F_q}},\Ql))^{(-1)^{n+1}}$
is the {\it global epsilon factor}, which is roughly the alternating product of the determinants of $H^i(X_{\overline{\F_q}},\Ql)$. It is known that these correcting terms of the functional equation admit {\it Euler product expansions}, 
as those of the usual zeta function have. 

To explain this precisely, take a map $f\colon X\to C$ to a projective smooth curve. According to 
the well-known philosophy that curves are geometric analogues of the rings of integers of number fields, a pair $(X,f)$ of a variety $X$ with a map to a curve in geometry should correspond to a variety over the ring of integers in arithmetic, and taking such a map should give decompositions of $\chi(X),\varepsilon(X,\mathbb{Q}_\ell)$ into local contributions. 
Indeed, these decompositions are realized by the Grothendieck--Ogg--Shafarevich formula (for $\chi(X)$) and Laumon's product formula (for $\varepsilon(X,\mathbb{Q}_\ell)$) applied to the pushforward $Rf_\ast\mathbb{Q}_\ell$ to the curve. 
From this point of view, the total dimensions of the 
vanishing cycles are nothing but the local terms appearing in the decompositions of $\chi(X)$, and 
on the other hand the local epsilon factors are the one appearing in those of the global epsilon factor $\varepsilon(X,\mathbb{Q}_\ell)$. Our results in this paper are thus to determine such local factors of the global epsilon factors in terms of 
geometry of singularities, with the assumption that the singularities are isolated.

To formulate the refinement for local epsilon factors, we will 
use linear-algebraic enhancements of the Milnor numbers, {\it non-degenerate symmetric bilinear forms}. 
Here is our main result in the odd characteristic cases. 

\begin{thm}(Theorem \ref{mainMil}.1, $p$ is odd) \label{thm1}
Let $k$ be a finite field of odd characteristic $p$. Let $X$ be a smooth $k$-scheme of dimension $n$ and let 
$f\colon X\to \A^1_k$ be a $k$-morphism with an isolated singular point $x\in X$. The singular point $x$ is assumed $k$-rational for simplicity in the introduction. 
Let 
$\varepsilon_0(\A^1_{k,(f(x))},R\Phi_f(\Ql)_x,dt)$ be the local epsilon factor and $(\varphi_{f},B_{f,dt})$ be the non-degenerate 
symmetric bilinear form over $k$ associated with $(f,x)$ (Definition \ref{phiB}). Then we have 
\begin{equation*}
(-1)^{{\rm dimtot}R\Phi_f(\Ql)_x}\varepsilon_0(\A^1_{k,(f(x))},R\Phi_f(\Ql)_x,dt)=\Bigl(\frac{(-2)^{n\mu(f,x)}\cdot{\rm disc}B_{f,dt}}{k}\Bigr) 
\cdot\tau_\psi^{(-1)^{n+1}n\mu(f,x)} 
\end{equation*}
where $(\frac{}{k})$ denotes the Legendre symbol and $\tau_\psi$ is the quadratic Gauss sum 
$\tau_{\psi}=-\sum_{a\in k}\psi(a^2)$. Here $\psi\colon k\to\overline{\Ql}^\times$ is a non-trivial additive character 
which is used in defining local epsilon factors (cf. \cite[(3.1.5.4)]{Lau}). 
\end{thm}
In this article, we treat local epsilon factors multiplied with the sign as in the left-hand side of the theorem, rather than the local epsilon factors themselves. 
This is because the natural generalization of local epsilon factor to a general perfect field $k$ is obtained with this sign, as is observed in 
 \cite{geomep} and \cite{Y3}. It is a character of the absolute Galois group of $k$, and if $k$ is finite, the value at the geometric Frobenius is the local epsilon factor with the sign. For this reason, we will use the symbol $\varepsilon_{0,\bar{k}}$, where $\bar{k}$ stands for an algebraic closure of $k$ 
used to define the Galois group, for local epsilon factors 
as characters in \cite{geomep}, \cite{Y3} and keep in mind the relation 
\begin{equation*}
\varepsilon_{0,\bar{k}}(T,V,\omega)({\rm Frob}_k)=(-1)^{{\rm dimtot}V}\varepsilon_0(T,V,\omega)
\end{equation*}
with the classical ones, in order to avoid any confusion on signs, which are main subjects of this paper.  

Attaching bilinear forms to isolated singularities is considered for example in \cite{EKL}
 and in 
\cite[Section 12]{enumquad}, 
in the spirit of refining classical formulae 
of Euler characteristics, which are valued in the integers, to the Grothendieck--Witt ring, i.e. the ring of non-degenerate symmetric bilinear forms over $k$. 
In our situation, taking the ranks of the bilinear forms gives the classical Milnor formula for the total dimensions. From this 
point of view, it can be said that what is done in this paper is to give 
\'etale counterparts of the discriminants of them. 

The most surprising results at least to the author in this paper appear in characteristic $2$, which we now explain. The construction of the bilinear forms from 
isolated singularities always works over any base scheme, even in characteristic $2$. However, the discriminants 
live in the group $k^\times/(k^\times)^2$, which is trivial if $k$ is of characteristic $2$ (and perfect). Of course, in order to 
construct a character of order $2$ of the absolute Galois group of $k$ (which will be a sign in the finite field case, 
by evaluating it at the geometric Frobenius), 
we need to find an element in $k/\wp(k)$, rather than the multiplicative group, according to Artin--Schreier theory. 
Here $\wp\colon k\to k$ is the map sending $x\mapsto x^2-x$. 

To get elements in $k/\wp(k)$, we consider liftings to the $p$-adic base. Precisely speaking, we consider a formal scheme formally smooth and 
formally of finite type over the Witt ring $W(k)$ with a morphism to the formal  $p$-adic completion of $\A^1_{W(k)}$ whose reduction mod $p$ 
is isomorphic to the initial one.  
First, take such one lift and fix it. Then, as is in the scheme case, we obtain a non-degenerate 
symmetric bilinear form over $W(k)$, whose discriminant lives in $W(k)^\times/(W(k)^\times)^2$. 
As is similar as in the case of Arf invariant (cf. \cite[1.4]{FQE}), we can prove that the discriminant does {\it not} depend on the choice of lifts and it 
lives in $\pm(1+4W(k))$. Therefore, choosing the sign appropriately, which is done in $4.3$ in this paper, we obtain an element of the form 
$1+4[a]$, hence an element $a\in k$, which is determined up to modulo $\wp(k)$. Here $[-]$ denotes the Teichm\"{u}ler lift. The constant $a\in k/\wp(k)$ so obtained we call {\it the Arf invariant}, and we denote by ${\rm Arf}(f,x)$. For the detail, see Section $4$, especially Theorems 
\ref{z/2ext}, \ref{z/2extk}. 

 \sloppy 
The terminology of Arf invariant is justified 
by the following example: Let $Q\in k[T_0,\dots,T_n]$ be a non-degenerate quadratic form with $k$-coefficient. Then the map 
$ Q\colon \A^{n+1}_k\to\nolinebreak\A^1_k$ has an isolated singular point at the origin and {\it our} Arf invariant ${\rm Arf}(Q,0)$ coincides with the original one. 

Actually, for our purpose, 
it is enough to consider liftings to $W_3$, the ring of Witt vectors of length $3$, instead of full 
liftings to $W$. 
This is because the reduction map $W\to W_3$ induces 
an isomorphism of the classifying spaces (i.e., $H^1$) of $\mu_2$-torsors 
over $W$ and $W_3$ (Lemma \ref{mapz/2mu_2}), where the discriminants of 
bilinear forms live. For $\mu_2$-torsors whose induced 
line bundles under the inclusion 
$\mu_2\to\mathbb{G}_m$ are trivial, this is equivalent 
to saying that $(W^\times)^2$ contains $1+8W$. 
For this reason, we state and prove theorems in characteristic $2$ in terms of $W_3$-liftings in Section $4$, which allows notations to 
be relatively simple. 

Using this Arf invariant, we can state the main result in characteristic $2$ as follows. 
\begin{thm}(Theorem \ref{mainMil}.2, $p=2$)
Let $k$ be a finite field with $q$ elements which is a power of $2$. Let 
$X$ be a smooth $k$-scheme of dimension $n$ and let $f\colon X\to\A^1_k$ be a morphism with a $k$-rational 
isolated singular point $x\in X$. 
Then $\frac{n\mu(f,x)}{2}$ is an integer and the ratio
\begin{equation*}
(-1)^{{\rm dimtot}R\Phi_f(\Ql)_x}\varepsilon_0(\A^1_{k,(f(x))},R\Phi_f(\Ql)_x,dt)/q^{\frac{(-1)^{n+1}n\mu(f,x)}{2}}
\end{equation*}
is $\pm1$. This sign is determined by the Arf invariant $a={\rm Arf}(f,x)$. Namely, the sign is $1$ if and only if 
$a$ belongs to $\wp(k)$. 
\end{thm}

It is very interesting to seek a theory which deals with the vanishing cycles of arbitrary $\ell$-adic sheaves  
 and contains our results as the special case of constant sheaves. 
 For the total dimensions, this is achieved 
by the theory of characteristic cycles given by T. Saito \cite{Sai17}. For the local epsilon factors, 
 epsilon cycles are defined in \cite{charep} as refinements of characteristic cycles. They 
treat the local epsilon factors of the vanishing cycles of $\ell$-adic sheaves 
in a geometric way using cotangent bundles, but require that we ignore roots of unity, namely 
we treat the local epsilon factors in $\overline{\Ql}^\times\otimes\mQ$, a quotient of the multiplicative group 
$\overline{\Ql}^\times$. 
This paper tells us that, even for constant sheaves, we need such 
complicated enhancements as symmetric bilinear forms to treat the 
local epsilon factors without taking modulo roots of unity. 
The author hopes that this paper can shed some new 
light on the theory of epsilon cycles, so that we eventually handle the local epsilon factors themselves without taking modulo roots of unity.  
In \cite{geomephigh}, Guignard gives another method for computing global epsilon factors in higher dimension, in a different way from epsilon cycles. It will also be an interesting work to 
clarify a relation between his results and ours. 

We briefly explain the strategy of the proof in the cases of odd characteristic. 
The proof for the case of characteristic $2$ is quite similar, although 
it is more involved. 

We explain the proof of Theorem \ref{thm1}. In contrast to Deligne's global approach in \cite{Mil}, our approach is 
local. Instead of referring to global results, such as the product formula of epsilon factors 
by Laumon, we use the continuity of local epsilon factors \cite{contiep}. Consider a family of morphisms from smooth $k$-schemes to 
$\A^1_k$ parametrized by a smooth base scheme $S$. Namely we 
consider a commutative diagram
\begin{equation}\label{famiso}
\xymatrix{
Y\ar[rd]_\pi\ar[rr]^{\tilde{f}}&&\A^1_S\ar[ld]\\&S
}
\end{equation}
of $k$-schemes of finite type 
where $\pi$ is smooth and $S$ is a smooth $k$-scheme. 
If a $k$-rational point $s\in S(k)$ is specified, we call such a family 
{\it a deformation of} $\tilde{f}_s$, where $\tilde{f}_s$ is the fiber of $\tilde{f}$ over $s$. 
The continuity says that, if the singular locus of $\tilde{f}$ is finite over the base $S$, 
the local epsilon factors vary continuously over $S$. More precisely, 
they satisfy the reciprocity law and give a 
character of the fundamental group of $S$. See Theorem \ref{dta} for the precise statement, or \cite[Theorem 4.8]{contiep} in a more general form of it. A geometric counterpart of this continuity is the flatness 
of the bilinear forms over $S$, which is rather trivial from the definition. 

Using this continuity, the proof goes as follows. Let $f\colon X\to \A^1_k$ be a morphism with an isolated singularity $x$. 
We argue by induction on the Milnor number $\mu(f,x)$. 
We construct a 
 deformation of $f$ as in (\ref{famiso}) so that its generic fibers contain at least one ordinary quadratic point with Milnor number $1$ or $2$ (Lemma \ref{deform}). By the continuity of local epsilon factors, the flatness of the bilinear forms, and 
 the induction hypothesis on the Milnor number, we reduce the proof to the case where the singularity is ordinary quadratic with the Milnor number $1$ or $2$, 
 in which we manage to 
 compute explicitly the both sides of the equality in the theorem (Proposition \ref{compphi}, Lemma \ref{epex}). 
 
 Let us explain the construction of this paper. We collect some  preliminaries 
 on Witt rings in Section \ref{Witpre}. We prove that a finite \'etale algebra $A$ over a Witt ring is integrally closed in the fraction ring $A[\frac{1}{p}]$ (Lemma \ref{Wet}). We also prove that the property for 
 a finite \'etale covering over the generic fiber of a Witt ring to extend 
 to a finite \'etale covering of the Witt ring is Zariski-local on the special fiber, in a certain sense (Lemma \ref{aaa}). Using these results, we prove results on $\mu_2$-torsors and $\mathbb{Z}/2$-coverings 
 over a $2$-typical Witt ring in Lemma \ref{mapz/2mu_2} which 
 will be needed in Section $4$. 
 Section \ref{singbil} 
 develops our main tools. After recalling basic results on residue symbol, we construct bilinear forms from isolated singularities and 
 study some of their basic properties.  
 We also include calculations on quadratic singularities which are 
 the first pieces of key results for our main theorems. 
The latter part is devoted to develop a machinary by which 
we reduce theorems to the case of quadratic singularities.  
In Section $4$, we define a finite \'etale $\mathbb{Z}/2$-covering 
from an isolated singularities, especially in characteristic $2$. 
 Finally, in Section $5$, we state and prove the main theorem of 
 this paper. After proving the main theorem, 
as an example of our results, we record a computation on the discriminants of the bilinear forms associated to 
homogeneous functions in the last subsection 5.3. In particular, 
in characteristic $2$, the Arf invariants of homogeneous functions  relate with the determinants of the \'etale cohomologies of middle degree of hypersurfaces (cf. \cite{hypdet}).

We collect terminologies we use throughout the paper. 
\begin{enumerate}
\item For a scheme $X$, we write $k(x)$ for the residue field 
at a point $x\in X$. 
\item Let $S$ be a scheme. For a point $s\in S$, we write 
$S_{(s)}$ for the henselization at $s$. On the other hand, 
for a geometric point $\bar{s}$ of $S$, we write 
$S_{(\bar{s})}$ for the strict henselization in $\bar{s}$. 
\item Let $S$ be a scheme and $s\to S$ be a morphism from the 
spectrum of a field. Subscripts $(-)_s$ indicates the base change 
to $s$. For example, for a morphism 
$f\colon X\to Y$ of $S$-schemes, we write 
$f_s\colon X_s\to Y_s$ for the base change to $s$. 
When $s={\rm Spec}(k)$, we also write $(-)_k$ for $(-)_s$. 
\item For a henselian local ring $(R,\mathfrak{m}_R)$, we write $R\{t_1,\dots,t_n\}$ for the henselization of $R[t_1,\dots,t_n]$ at 
$(\mathfrak{m}_R,t_1,\dots, t_n)$. 
\end{enumerate}

\tableofcontents

\section{On finite \'etale algebras over Witt rings}\label{Witpre}

Let $p$ be a prime number, which we fix in the subsection \ref{Witcom}. Let $A$ be a perfect $\mathbb{F}_p$-algebra, i.e. the 
Frobenius map $A\xrightarrow{x\mapsto x^p}A$ is an isomorphism. 
Let $W(A)$ be the ($p$-typical) Witt ring with coefficients in $A$. 
 In this preliminary section, 
we give a criterion for a $W(A)$-algebra of a certain type to be finite \'etale (Lemma \ref{aaa}). Using this lemma, we give some results on 
$\mathbb{Z}/2$-coverings and 
$\mu_2$-torsors when $p=2$ (Lemma \ref{mapz/2mu_2}). 

\subsection{Finite modules on adic rings}
We start with recalling some basics on finite modules over adic rings. 
Let $R$ be a ring with a non-zero divisor $\varpi$ such that $R$ is $\varpi$-adically complete and separated. 
We write $R_n$ for $R/\varpi^{n+1}$. 
\begin{df}\label{dffmod}
We define ${\rm fMod}_\bullet(R)$ to be the category of projective systems $(M_n)_{n\geq0}$ indexed by the integers $n\geq0$ with the 
following properties: 
\begin{enumerate}
\item 
For each $n\geq0$, $M_n$ is an $R_n$-module and 
the transition map $M_{n+1}\to M_n$ is $R_{n+1}$-linear. $M_0$ is a finite $R_0$-module. 
\item 
The transition maps 
induce quasi-isomorphisms 
$M_{n+1}\otimes^L_{R_{n+1}}R_n\to M_n$ of complexes of $R_n$-modules. 
\end{enumerate}
\end{df}
For  a $\varpi$-torsion free $R$-module $M$, the system $(M/\varpi^{n+1}M)_n$ is an object of ${\rm fMod}_\bullet(R)$ provided that 
$M/\varpi M$ is finite.

\begin{rmk}\label{fres}
\begin{enumerate}
\item By the condition $2$, we especially have 
$M_n\otimes_{R_n}R_0\cong M_0$. Thus, $M_n$ is finite since 
$M_0$ is finite and the kernel of $R_n\to R_0$ is nilpotent.  

\item Let $n\geq i\geq0$ be integers. 
Composing the quasi-isomorphisms $M_{n+1}\otimes^L_{R_{n+1}}R_n\to M_n$ repeatedly, we know that 
the map $M_{n+1}\otimes^L_{R_{n+1}}R_i\to M_i$ 
is a quasi-isomorphism. 
Since $\varpi$ is not a zero-divisor in $R$, 
\begin{equation*}
\cdots\to R_{n+1}\xrightarrow{\varpi^{n+1-i}}R_{n+1}\xrightarrow{\varpi^{i+1}}R_{n+1}\xrightarrow{\varpi^{n+1-i}}R_{n+1}\xrightarrow{\varpi^{i+1}} R_{n+1}
\end{equation*}
gives an $R_{n+1}$-free resolution of $R_i$. Therefore, the condition that $M_{n+1}\otimes^L_{R_{n+1}}R_i\to M_i$ is a 
quasi-isomorphism for any $n\geq i\geq0$ is equivalent to the following two conditions: 
\begin{enumerate}
\item The map $M_{n+1}\to M_i$ induces an isomorphism $M_{n+1}/\varpi^{i+1}M_{n+1}\to M_i$. 
\item $M_{n+1}\xrightarrow{\varpi^{n+1-i}}M_{n+1}\xrightarrow{\varpi^{i+1}}M_{n+1}$ is exact. 
\end{enumerate}
\end{enumerate}
\end{rmk}
\begin{lm}\label{projlimf}
Let $(M_n)_n$ be an object in ${\rm fMod}_\bullet(R)$ (Definition \ref{dffmod}). Let $M:=\varprojlim_nM_n$ be the projective limit. 
\begin{enumerate}
\item The limit $M$ is a finite $\varpi$-torsion free $R$-module. 
The projection $M\to M_n$ induces an isomorphism $M/\varpi^{n+1}M\to M_n$. 
\item If $M_0$ is finitely presented as an $R_0$-module, $M$ is finitely presented as an $R$-module. 
\item If $M_0$ is projective as an $R_0$-module, $M$ is projective as an $R$-module. 
\end{enumerate}
\end{lm}
\proof{
$1$.  First we show that $M/\varpi^{n+1}M\cong M_n$. By Remark \ref{fres}, we have a short exact sequence 
\begin{equation*}
0\to M_s\xrightarrow{\varpi^{n+1}}M_{n+s+1}\to M_n\to0. 
\end{equation*}
Taking the limit with respect to $s$, we have 
\begin{equation*}
0\to M\xrightarrow{\varpi^{n+1}}M\to M_n\to0. 
\end{equation*}
Here the exactness on the right comes from that $(M_s)_s$ satisfies the Mittag-Leffler condition. Therefore we have 
$M/\varpi^{n+1}M\cong M_n$. 

Next we show that $M$ is finite. 
Since $M_0$ is finite, we can choose finitely many generators $x_{0,0},\dots,x_{0,m}\in M_0$. Inductively on $n$, we choose elements 
$x_{n,j}\in M_n$ such that the images of $x_{n+1,j}$ in $M_n$ is $x_{n,j}$. Since the kernel of 
$R_{n}\to R_0$ is nilpotent, $x_{n,0},\dots,x_{n,m}$ generate $M_n$. Then $x_j:=(x_{n,j})_n\in M$ generate $M$. 
Indeed, take an element $x=(x_n)\in M$. First choose $a_j\in R$ so that the images $\bar{a}_j\in R_0$ of $a_j$ satisfy the relation 
$x_0=\sum_j\bar{a}_jx_{0,j}$. Then $x':=x-\sum_ja_jx_{j}$ has $0$ in its $0$-th component. By the isomorphism $M/\varpi M\cong M_0$, 
$x'$ can be written in the form $\varpi y$ for some $y\in M$. 
 Applying the same procedure to $y$, we find $b_j\in R$ such that $y-\sum_jb_jx_j$ has $0$ in its $0$-th component. 
 Construct sequences $a_j,b_j,\dots$ for each $j$ inductively and we define $\alpha_j\in R$ to be the limits $a_j+\varpi b_j+\cdots$. 
 Then we have $x=\sum_j\alpha_jx_j$, which shows the finiteness.  
 
We prove that $M$ is $\varpi$-torsion free. 
Consider an exact sequence 
\begin{equation*}
0\to K_n\to M_n\xrightarrow{\varpi}M_n. 
\end{equation*}
By the left exactness of $\varprojlim_n$, the kernel of $\varpi\colon M\to M$ is isomorphic to $\varprojlim_nK_n$. 
Since $K_{n+1}$ maps to $0$ via $M_{n+1}\to M_n$ (Remark \ref{fres}), this limit vanishes, which implies the claim.

$2$. 
Let $R^m\to M$ be a surjection from a free module and let $N$ be its kernel. Taking modulo $\varpi^{n+1}$, we have a short exact sequence 
\begin{equation}\label{modpi2}
0\to N/\varpi^{n+1}N\to R_n^m\to M_n\to0. 
\end{equation}
Since $M_0$ is a finitely presented $R_0$-module, $N/\varpi N$ is finitely generated. 
Since $N$ is $\varpi$-torsion free, the maps  $N/\varpi^{n+2}N\otimes_{R_{n+1}}^LR_n\to N/\varpi^{n+1}N$ are 
quasi-isomorphisms. Thus, the system $(N/\varpi^{n+1}N)_n$ is an element of ${\rm fMod}_\bullet(R)$. 
Applying the five lemma to the limit of (\ref{modpi2}), $N$ is isomorphic to the limit $\varprojlim_nN/\varpi^{n+1}N$, which is finitely generated by 1. 

$3$. We know that $M_n$ is finitely presented as an $R_n$-module by $2$ and $R_n$-flat by \cite[22.3]{Mat}, hence projective. 

Take a surjection $R^m\to M$ and let $f_n\colon R_n^m\to M_n$ be its reduction mod $\varpi^{n+1}$. We construct a splitting 
$\varphi_n\colon M_n\to R_n^m$ for $f_n$ by induction on $n$. For $n=0$, we may take any splitting as $\varphi_0$. 
Suppose that we choose $\varphi_n$. By the projectivity of $M_{n+1}$, we find a map 
$\varphi'\colon M_{n+1}\to R_{n+1}^m$ which makes the diagram
\begin{equation*}
\xymatrix{
M_{n+1}\ar[d]^{\varphi'}\ar[r]&M_n\ar[d]^{\varphi_n}\\
 R_{n+1}^m\ar[r]&R_n^m
 }
 \end{equation*}
where the horizontal arrows are the transition maps commutative. Since the composition $f_{n+1}\circ\varphi'$ becomes the identity after taking modulo $\varpi^{n+1}$, it is an isomorphism. 
Then we can take $\varphi_{n+1}:=\varphi'\circ (f_{n+1}\circ\varphi')^{-1}$. 
\qed}

\begin{rmk}
Let ${\rm fMod}(R)$ be the category of finitely generated $R$-modules 
which are $\varpi$-torsion free. Lemma \ref{projlimf}.1 implies that 
the assignment $(M_n)_n\mapsto \varprojlim_nM_n$ defines a 
functor ${\rm fMod}_\bullet(R)\to{\rm fMod}(R)$. This 
is fully faithful and its essential image consists of 
the elements $M\in {\rm fMod}(R)$ which are $\varpi$-adically 
complete and separated. By Lemma \ref{projlimf}.1 and Nakayama's lemma, 
an element $M\in  {\rm fMod}(R)$ is always $\varpi$-adically complete, i.e. $M\to\varprojlim_nM/\varpi^{n+1}M$ is surjective, but 
not separated in general. 
However, the following holds. 
\end{rmk}
The following lemma is not used in the sequel. 
\begin{lm}
Assume that $M\in{\rm fMod}(R)$ is finitely presented. Then it is 
$\varpi$-adically complete and separated. 
\end{lm}
\proof{
Let us denote by $\hat{L}$ the $\varpi$-adic completion of an 
$R$-module $L$. 

The surjectivity of $M\to\hat{M}$ follows from 
Lemma \ref{projlimf}.1 and Nakayama's lemma, for which we merely 
use the property that $M$ is finitely generated. 
Take a short exact sequence
\begin{equation}\label{abc}
0\to N\to R^m\to M\to0
\end{equation}
where $N$ is a finitely generated $R$-module. 
As $M$ is $\varpi$-torsion free, 
the $\varpi^{n+1}$-mod reductions of (\ref{abc}) remain exact. 
As $(N/\varpi^{n+1}N)_n$ satisfies the Mittag-Leffler condition, the 
$\varpi$-adic completion of $(\ref{abc})$ is exact. 
We have a morphism of short exact sequences 
\begin{equation*}
\xymatrix{
0\ar[r]&N\ar[d]\ar[r]&R^m\ar[d]_{\cong}\ar[r]&M\ar[r]\ar[d]&0\\
0\ar[r]&\hat{N}\ar[r]&\hat{R^m}\ar[r]&\hat{M}\ar[r]&0. 
}
\end{equation*}
The injectivity of $M\to\hat{M}$ follows from the surjectivity of 
$N\to\hat{N}$. 
\qed}

\begin{pr}\label{etmod}
\begin{enumerate}
\item 
Let $(R'_n)_n$ be an object in ${\rm fMod}_\bullet(R)$. Assume that $R'_n$ are equipped with structures of $R_n$-algebras 
such that the transition maps $R'_{n+1}\to R'_n$ are morphisms of $R$-algebras. Further assume that $R'_0$ is a finite \'etale 
$R_0$-algebra. Then the limit 
$R':=\varprojlim_nR'_n$ is a finite \'etale $R$-algebra. 
\item The modulo-$\varpi$ reduction gives an equivalence of 
categories from finite \'etale algebras over $R$ to 
those over $R_0$. 
\end{enumerate}
\end{pr}
\proof{
 $1$. 
By Lemma \ref{projlimf}.3, we know that $R'$ is finite projective as an $R$-module. Thus it is enough to show that $\Omega^1_{R'/R}$ vanishes, which follows from the vanishing of $\Omega^1_{R'_0/R_0}$ 
and Nakayama's lemma. 

$2$. We construct a quasi-inverse functor. Let $R'_0$ be a finite \'etale 
$R_0$-algebra. For $n\geq0$, the closed immersion 
${\rm Spec}(R_0)\to{\rm Spec}(R_n)$ is defined by a nilpotent ideal. 
Hence $R'_0$ uniquely lifts to a finite \'etale $R_n$-algebra 
$R'_n$. Then the limit $R':=\varprojlim_nR'_n$ is a finite \'etale 
$R$-algebra by $1$, which is quasi-inverse to the modulo-$\varpi$ 
reduction. 
\qed}

We recall that the Picard groups of adic rings are isomorphic to those of 
the reductions. 
\begin{lm}\label{Gmisom}
Let $R$ be a ring which is $\varpi$-adically complete and separated. Then the reduction 
map $R\to R_0=R/\varpi$ induces an isomorphism of groups 
$H^1({\rm Spec}(R),\mathbb{G}_m)\to H^1({\rm Spec}(R_0),\mathbb{G}_m)$. 
\end{lm}
\proof{
The injectivity follows from Nakayama's lemma. We show the surjectivity. 
Take an invertible $R_0$-module $M_0$. We construct an invertible $R_n$-module 
$M_n$ for $n>0$ such that $M_n\otimes_{R_n}R_{n-1}\cong 
M_{n-1}$. Then the limit 
$M:=\varprojlim_nM_n$ is an invertible $R$-module with 
$M\otimes_{R}R_0\cong M_0$ by Lemmas \ref{projlimf}.1,3. 

The construction goes by induction on $n$. 
Suppose that we are given an invertible $R_n$-module $M_n$. 
To show that $M_n$ lifts to an invertible $R_{n+1}$-module, 
we show that the map 
$H^1(S_{n+1},\mathbb{G}_m)\to H^1(S_n,\mathbb{G}_m)$ is surjective, where we set $S_m={\rm Spec}(R_m)$ for integers $m\geq0$. We have a short exact sequence 
\begin{equation*}
0\to \varpi^{n+1}\OO_{S_{n+1}}\xrightarrow{x\mapsto1+x}\OO_{S_{n+1}}^\times\to\OO_{S_n}^\times\to0. 
\end{equation*}
Then the surjectivity follows from $H^2(S_{n+1},\varpi^{n+1}\OO_{S_{n+1}})=0$, as 
$S_{n+1}$ is affine. 
\qed}

\subsection{On Witt rings}\label{Witcom}
From now on, we suppose that $R=W(A)$ is the Witt ring of a perfect $\mathbb{F}_p$-algebra $A$.  
This is $p$-adically complete separated, and $p$ is not a zero-divisor in it. Therefore, we can apply the previous results to $(R,\varpi)=(W(A),p)$. 

First we show that finite \'etale algebras over the Witt rings are 
``normal''. 
\begin{lm}\label{Wet}
Let $R'$ be a finite \'etale $W(A)$-algebra. 
Then $R'$ is the integral closure of $W(A)$ in $R'[\frac{1}{p}]$. 
\end{lm}
\proof{
We show that $R'$ is integrally closed in $R'[\frac{1}{p}]$. 
Take $a\in R'[\frac{1}{p}]$ which is integral over $R'$. 
It can be written in the form $a=b/p^n$ for some $b\in R'$ and a non-negative integer $n\geq0$. We choose a presentation 
so that $n$ is minimal. Since $a$ is integral over $R'$, there exist 
$a_1,\dots,a_m\in R'$ such that 
\begin{equation*}
a^m+a_1a^{m-1}+\cdots+a_m=0. 
\end{equation*}
Multiplying $p^{nm}$, we get 
\begin{equation*}
b^m+a_1 p^nb^{m-1}+\cdots+a_mp^{nm}=0.
\end{equation*}
Therefore the class $\bar{b}$ of $b$ in $R'/pR'$ is nilpotent if 
$n>0$. But $R'/pR'$ is reduced since it is \'etale over $W(A)/p=A$, which is reduced. 
This implies that $\bar{b}$ is zero, i.e. $b$ is divisible by $p$ in $R'$, a contradiction to the minimality of $n$. 
Thus we have $n=0$, hence $a\in R'$. 
\qed}

Next we give a criterion for a normalization to be \'etale in Lemma 
\ref{aaa}. 
We start with the following lemmas. 
\begin{lm}\label{genred}
Let $\{\eta_i\}_{i\in I}$ be the set of generic points of $A$. Then the map 
$A\to\prod_ik(\eta_i)$ is injective. 
\end{lm}
\proof{
It follows since $A$ is reduced. 
\qed}

\begin{lm}\label{dual}
Let $R'$ be a finite $W(A)$-algebra which is projective as a $W(A)$-module. Assume that $R'[\frac{1}{p}]$ is finite \'etale over 
$W(A)[\frac{1}{p}]$. 
Let $R''$ be the integral closure of $R'$ in $R'[\frac{1}{p}]$. 
\begin{enumerate}
\item
There exists a finite projective $W(A)$-submodule 
$R'^\vee$ of $R'[\frac{1}{p}]$ which contains $R''$. 
\item $R''$ is $p$-adically complete and separated, i.e. the canonical map $R''\to\varprojlim_nR''/p^{n+1}R''$ 
is an isomorphism. 
\end{enumerate}
\end{lm}
\proof{
$1$. 
The trace on $R'$ induces a pairing 
\begin{equation*}
R'\times R'\to W(A),\ (x,y)\mapsto{\rm Tr}_{R'/W(A)}(xy). 
\end{equation*}
This induces a $W(A)$-linear map $f\colon R'\to R'^\vee$ where we set $R'^\vee={\rm Hom}_{W(A)}(R',W(A))$. 
As $R'[\frac{1}{p}]$ is finite \'etale over $W(A)[\frac{1}{p}]$, the map $f$ becomes an isomorphism after inverting $p$. 
Thus $R'^\vee$ can be seen as a $W(A)$-submodule of $R'[\frac{1}{p}]$. 

To show $R''\subset R'^\vee$, it suffices to show that, 
for $x\in R'[\frac{1}{p}]$ integral over $W(A)$, we have 
${\rm Tr}_{R'/W(A)}(x)\in W(A)$. By Lemma \ref{genred}, 
it is enough to show that the image of ${\rm Tr}_{R'/W(A)}(x)$ 
in $W(k(\eta))[\frac{1}{p}]$ for any generic point $\eta$ of 
$A$ is contained in $W(k(\eta))$. Hence 
we reduce it to the case where $A$ is a perfect field. 
In this case, it is well-known as $W(A)$ is a normal domain. 

$2$. We show that the $p$-adic topology on $R''$ induces the $p$-adic topology on $R'$ and that $R'$ is open in $R''$ in 
that topology. Indeed, for any integer $n\geq0$, there exists $N\geq0$ such that 
$p^NR'^\vee\subset p^nR'$ as 
$R'^\vee$ is finitely generated and 
$R'^\vee/ R'$ is $p$-torsion. 
Hence we have $p^NR''\subset p^nR'$ for such an $N$. 

We find that $R''$ contains an open submodule $R'$ which is complete and separated in the $p$-adic topology 
of $R''$. Thus $R''$ itself is $p$-adically complete and 
separated. 
\qed}

We denote the Teichm\"{u}ller lift $A\to W(A)$ by $a\mapsto[a]$. 
\begin{lm}\label{aaa}
Let $R'$ be a finite $W(A)$-algebra which is projective as a $W(A)$-module. Assume that $R'[\frac{1}{p}]$ is finite \'etale over $W(A)[\frac{1}{p}]$. 
Let $R''$ be the integral closure of $W(A)$ in $R'[\frac{1}{p}]$. 
\begin{enumerate}
\item Let $g\in A$ be an element. 
The $p$-adic completion of $R''\otimes_{W(A)}W(A)[\frac{1}{[g]}]$ is isomorphic to the $p$-adic completion 
of $R''\otimes_{W(A)}W(A_{g})$. 
We denote it by $R''\hat{\otimes}_{W(A)}W(A_{g})$. 
\item The $p$-adic completion $R''\hat{\otimes}_{W(A)}W(A_{g})$, defined in $1$, is the integral closure 
of $W(A_{g})$ in $R'[\frac{1}{p}]\otimes_{W(A)}W(A_{g})$. 
\item 
Let $g_1,\dots,g_n\in A$ be a sequence of elements which generates the unit ideal of $A$. 
If $R''\hat{\otimes}_{W(A)}W(A_{g_i})$ is a finite \'etale $W(A_{g_i})$-algebra for each $i$, $R''$ is a finite \'etale 
$W(A)$-algebra. 
\end{enumerate}
\end{lm}
\proof{
By Lemma \ref{dual}.1, there exists a finite projective $W(A)$-submodule $R'^\vee$ of $R'[\frac{1}{p}]$ which contains $R''$,~i.e.~$R'\subset R''\subset R'^\vee$. 

$1$. The assertion follows since the mod $p^{n+1}$ reductions of the two rings are isomorphic. 

$2$. The ring $R''\otimes_{W(A)}W(A)[\frac{1}{[g]}]$ is the integral closure 
of $W(A)[\frac{1}{[g]}]$ in $R'[\frac{1}{p}]\otimes_{W(A)}W(A)[\frac{1}{[g]}]$. 
We have the inclusions 
\begin{equation*}
R'\otimes_{W(A)}W(A)[\frac{1}{[g]}]\subset R''\otimes_{W(A)}W(A)[\frac{1}{[g]}]\subset R'^\vee
\otimes_{W(A)}W(A)[\frac{1}{[g]}]. 
\end{equation*}
This implies that 
the $p$-adic topology of $R''\otimes_{W(A)}W(A)[\frac{1}{[g]}]$ induces the same topology 
on $R'[\frac{1}{p}]\otimes_{W(A)}W(A)[\frac{1}{[g]}]$ as that of $R'\otimes_{W(A)}W(A)[\frac{1}{[g]}]$ does. 
Thus the $p$-adic completion $R''\hat{\otimes}_{W(A)}W(A_{g})$ is integrally closed in the completion of 
$R'[\frac{1}{p}]\otimes_{W(A)}W(A)[\frac{1}{[g]}]$, which is $R'[\frac{1}{p}]\otimes_{W(A)}W(A_{g})$. 

It remains to show that $R''\hat{\otimes}_{W(A)}W(A_{g})$ is integral over $W(A_{g})$. 
Take any element $x\in R''\hat{\otimes}_{W(A)}W(A_{g})$. Then there exists $y\in R''\otimes_{W(A)}W(A)[\frac{1}{[g]}]$ such 
that $x-y\in R'\otimes_{W(A)}W(A_{g})$ as $R'\otimes_{W(A)}W(A_{g})$ is an open subring. Thus $x$ is integral.  

$3$. Assume that $R''\hat{\otimes}_{W(A)}W(A_{g_i})$ is finite \'etale over $W(A_{g_i})$. 
By the Zariski descent, $R''/pR''$ is finite \'etale over $A$. Thus, by Proposition \ref{etmod}.1, 
$\varprojlim_nR''/p^{n+1}R''$ is a finite \'etale $W(A)$-algebra. On the other hand, we have 
$R''\cong\varprojlim_nR''/p^{n+1}R''$ by Lemma \ref{dual}.2, hence the assertion. 
\qed}

We recall a well-known result on $W(A)^\times/(W(A)^\times)^2$ in the 
case of $p=2$. 

\begin{lm}\label{Wunr}
Let $p=2$ and let $A$ be a perfect $\F_2$-algebra.  
\begin{enumerate}
\item We have the inclusion $1+8W(A)\subset (W(A)^\times)^2$. 
\item There exists a canonical short  exact sequence 
\begin{equation}\label{exseq}
0\to A/\wp(A)\to W(A)^\times/(W(A)^\times)^2\to A\to0
\end{equation}
where $\wp\colon A\to A$ denotes the map $x\mapsto x^2-x$. The first arrow sends $a\in A$ to the class of $1+4[a]\in W(A)^\times$ and 
the second arrow sends $1+2[a]$ to $a$.  
\item Let $\alpha$ be an element of $W(A)^\times$. 
The following are equivalent. 
\begin{enumerate}
\item The $W(A)[\frac{1}{2}]$-algebra $W(A)[\frac{1}{2}][u]/(u^2-\alpha)$ extends to a (necessarily 
unique) finite 
\'etale $W(A)$-algebra. 
\item For any generic point $\eta$ of $A$,
the extension $F_\eta[u]/(u^2-\bar{\alpha})$ where 
$F_\eta$ denotes the fraction field 
of $W(k(\eta))$ and $\bar{\alpha}$ is the image of $\alpha$ under the map $W(A)\to W(k(\eta))$ is unramified. 

\item The class of $\alpha$ in $W(A)^\times/(W(A)^\times)^2$ comes from $A/\wp(A)$ under the injection $A/\wp(A)\to W(A)^\times/(W(A)^\times)^2$ 
in (\ref{exseq}). 
\end{enumerate}
Suppose that these conditions are satisfied and set $\alpha\equiv1+4[b]$ in $W(A)^\times/(W(A)^\times)^2$. 
Let $R'$ be a finite \'etale $W(A)$-algebra such that 
$R'[\frac{1}{2}]\cong W(A)[\frac{1}{2}][u]/(u^2-\alpha)$. 
Then 
the quadratic character in $ H^1({\rm Spec}(A),\mathbb{Z}/2)$ 
given by the quadratic extension $R'\otimes_{W(A)}A$ of $A$ corresponds to $b$ via the Artin--Schreier isomorphism $A/\wp(A)\cong H^1({\rm Spec}(A),\mathbb{Z}/2)$. 
\end{enumerate}
\end{lm}
\proof{
The exponential $\exp(x)=1+x+\frac{x^2}{2!}+\cdots$ converges for $x\in 4 W(A)$ and the logarithm $\log(x)=
(x-1)-\frac{(x-1)^2}{2}+\frac{(x-1)^3}{3}-\cdots$ converges for $x\in 1+4W(A)$. They define group isomorphisms between 
$1+4W(A)$ and $4W(A)$, each of which is the inverse to the other. They also respect the filtrations $1+2^iW(A)\subset 1+4W(A)$ and 
$2^iW(A)\subset 4W(A)$ for $i\geq2$. In particular, we have 
$(1+4W(A))^2=1+8W(A)$, hence the assertion $1$. On the other hand, since $A$ is perfect of characteristic $2$, we have $([A^\times])^2=[A^\times]$. Thus, letting $M:=(1+2W(A))/
(1+8W(A))$, whose group law we write additively, we have $W(A)^\times/(W(A)^\times)^2\cong M/2M$. 

Consider a commutative diagram 
\begin{equation*}
\xymatrix{
0\ar[r]& A\ar[r]\ar[d]^2& M\ar[r]\ar[d]^2& A\ar[d]^2\ar[r]& 0\\
0\ar[r]& A\ar[r]& M\ar[r]& A\ar[r]& 0
}
\end{equation*}
where the left horizontal arrows $A\to M$ are $a\mapsto 1+4[a]$ and the right horizontal arrows $M\to A$ are 
$1+2[a]+4[b]\mapsto a$. By the snake lemma, we have 
\begin{equation*}
A\to A\to M/2M\to A\to0. 
\end{equation*}
It is straightforward to check that the connecting homomorphism $A\to A$ is $\wp$. The assertion $2$ follows. 

$3$. The uniqueness in $(a)$ follows from Lemma \ref{Wet}. 
The implication $(a)\Rightarrow(b)$ is clear. 

We show $(c)\Rightarrow(a)$. Write $\alpha=\beta^2(1+4[b])$ where  
$\beta\in W(A)^\times, b\in A$. Putting $\beta(1+2v)=u$, 
$W(A)[\frac{1}{2}][u]/(u^2-\alpha)$ is isomorphic to $W(A)[\frac{1}{2}][v]/(v^2+v-[b])$. This extends to a finite \'etale 
$W(A)$-algebra $W(A)[v]/(v^2+v-[b])$. The last assertion on the 
quadratic character 
follows from this presentation. 

We show $(b)\Rightarrow(c)$. Let $\alpha\in W(A)^\times$ be an element such that 
$F_\eta[u]/(u^2-\bar{\alpha})$ is unramified for any $\eta$. Write $\alpha=\beta^2(1+2[a]+4[b])$
 where $\beta\in W(A)^\times, a,b\in A$. 
We need to show that $a=0$. 
Suppose otherwise. By Lemma \ref{genred}, there exists a 
generic point $\eta$ such that the image $\bar{a}$ of $a$ in 
$W(k(\eta))$ is non-zero. 
This means that $2[a]+4[b]$ maps to a uniformizer $\pi$ in $W(k(\eta))$. 
Putting $u=\beta(w+1)$, we have 
$F_\eta[u]/(u^2-\alpha)\cong F_\eta[w]/(w^2+2w-\pi)$, which is a non-trivial totally ramified extension, a contradiction. 
\qed}

Using the above results, we prove some properties on $\mu_2$-torsors and 
$\mathbb{Z}/2$-coverings over Witt rings which will be needed in Section $4$. 
Let $A$ be a perfect $\F_2$-algebra and $W_3(A)$ be the 
ring of Witt vectors of length $3$. 
Write 
$S_2:={\rm Spec}(W_3(A))$ (we use the subscript $2$ in the left since we write 
$S_0$ for ${\rm Spec}(A)$ in Section $4$). 
Write $\mu_2$ for the fppf sheaf on schemes sending a scheme $T$ to the group 
\begin{equation*}
\{a\in\Gamma(T,\OO_T)\vert a^2=1\}. 
\end{equation*}
This fits into the short exact sequence of fppf sheaves 
\begin{equation}\label{sexmu}
0\to\mu_2\to\mathbb{G}_m\xrightarrow{a\mapsto a^2}\mathbb{G}_m\to0. 
\end{equation}
\begin{lm}\label{mapz/2mu_2}
Let $p=2$. Consider a map 
$\mathbb{Z}/2\to\mu_2$ of fppf sheaves sending $\mathbb{Z}/2\ni1\mapsto-1\in\mu_2$. 
\begin{enumerate}
\item Write 
$S_\infty:={\rm Spec}(W(A))$. Then the maps 
$H^1(S_\infty,\mathbb{Z}/2)\to H^1(S_2,\mathbb{Z}/2)$ and 
$H^1(S_\infty,\mu_2)\to H^1(S_2,\mu_2)$ are isomorphisms. 

\item Write $S_\infty[\frac{1}{2}]={\rm Spec}
(W(A)[\frac{1}{2}])=S_\infty\times_{\mathbb{Z}}\mathbb{Z}[\frac{1}{2}]$. The map $\mathbb{Z}/2\to\mu_2$ induces a commutative diagram 
of cohomology groups 
\begin{equation}\label{commmuz}
\xymatrix{
H^1(S_2,\mathbb{Z}/2)\ar[r]\ar[rd]&H^1(S_2,\mu_2)\ar[d]\\
&H^1(S_\infty[\frac{1}{2}],\mathbb{Z}/2). 
}
\end{equation}
The slant arrow is injective. Consequently, 
the horizontal map 
$H^1(S_2,\mathbb{Z}/2)\to H^1(S_2,\mu_2)$ is injective. 
\item 
Let $\{\eta_i\}_{i\in I}$ be the set of generic points of $A$. 
Set $S_{2,i}:={\rm Spec}(W_3(k(\eta_i)))$. 
Then the diagram 
\begin{equation*}
\xymatrix{
H^1(S_2,\mathbb{Z}/2)\ar[r]\ar[d]&H^1(S_2,
\mu_2)\ar[d]\\
\Pi_i H^1(S_{2,i},\mathbb{Z}/2)\ar[r]&\Pi_i
H^1(S_{2,i},\mu_2)
}
\end{equation*}
is cartesian. 

\item Further assume that $A$ is normal, i.e. 
$A$ is the finite product of normal domains. 
Then the vertical arrow $H^1(S_2,\mathbb{Z}/2)\to\prod_i
H^1(S_{2,i},\mathbb{Z}/2)$ in $3$ is injective. 

\end{enumerate}
\end{lm}
\proof{
$1$. The first map $H^1(S_\infty,\mathbb{Z}/2)\to H^1(S_2,\mathbb{Z}/2)$ is an isomorphism by Proposition \ref{etmod}.2. 

We show that the second one is also an isomorphism. The Kummer 
exact 
sequence (\ref{sexmu}) gives a commutative diagram 
\begin{equation}
\xymatrix{
0\ar[r]&W(A)^\times/(W(A)^\times)^2\ar[r]\ar[d]&
H^1(S_\infty,\mu_2)\ar[r]\ar[d]&H^1(S_\infty,\mathbb{G}_m)[2]
\ar[r]\ar[d]&0\\
0\ar[r]&W_3(A)^\times/(W_3(A)^\times)^2\ar[r]&
H^1(S_2,\mu_2)\ar[r]&H^1(S_2,\mathbb{G}_m)[2]
\ar[r]&0
}
\end{equation}
where the horizontal lines are exact and $G[2]$ for an abelian group $G$ denotes the $2$-torsion part. 
By Lemma \ref{Wunr}.1, the left vertical arrow is an isomorphism. 
Applying Lemma \ref{Gmisom} to $(R,\varpi)=(W(A),p^3)$, we know that 
the right vertical arrow is an isomorphism. 
Therefore, the assertion follows from the five lemma.

$2$. 
Since the map $\mathbb{Z}/2\to\mu_2$ is 
an isomorphism on $S_\infty[\frac{1}{2}]$, the map 
$H^1(S_\infty,\mathbb{Z}/2)\to H^1(S_\infty[\frac{1}{2}],\mathbb{Z}/2)
$ factors through $H^1(S_\infty,\mu_2)$. Then the diagram is obtained from this under the identifications given in $1$. 
The injectivity of the horizontal arrow follows from that 
of the slant one, which is a consequence of Lemma \ref{Wet}.

$3$. We may replace $S_2,S_{2,i}$ in the diagram by $S_\infty$ and $
S_{\infty,i}={\rm Spec}(W(k(\eta_i)))$ by $1$. 
We already know that the horizontal arrows are injective by $2$. 
Hence, it is enough to show that 
a $\mu_2$-torsor $S'$ on $S_\infty$ comes from a (unique) $\mathbb{Z}/2$-covering on $S_\infty$ provided that 
its restrictions to $S_{\infty,i}$ come from $\mathbb{Z}/2$-coverings. 

Let $R':=\Gamma(S',\OO_{S'})$, which is finite projective as a 
$W(A)$-module. We show that the normalization 
$R''$ of $R'$ in $R'[\frac{1}{2}]$ is finite \'etale over $W(A)$ (cf. Lemma \ref{Wet}). 
Then the action of $\mathbb{Z}/2$ on $R'[\frac{1}{2}]$ uniquely lifts to 
an action on $R''$, which gives an element in $H^1(S_\infty, 
\mathbb{Z}/2)$ mapping to $S'$. 

The exact sequence (\ref{sexmu}) gives an exact sequence
\begin{equation}\label{mugm}
0\to W(A)^\times/(W(A)^\times)^2\to
H^1(S_\infty,\mu_2)\to H^1(S_\infty,\mathbb{G}_m). 
\end{equation}
Let $M$ be an invertible $W(A)$-module which is the image of 
$S'$ via the right arrow in (\ref{mugm}). Take elements 
$g_1,\dots,g_m\in A$ such that $M\otimes_{W(A)}A_{g_i}$ are 
free and $g_1\dots, g_m$ generate the unit ideal of $A$. 
Then Nakayama's lemma implies that $M\otimes_{W(A)}W(A_{g_i})$ 
are free. By Lemma \ref{aaa}, we may replace $A$ by each 
$A_{g_i}$ to show that $R''$ is finite \'etale over $W(A)$. 
 Hence we may assume that $M$ itself is free. In this case, 
the exact sequence (\ref{mugm}) tells us that there exists an element 
$\alpha\in W(A)^\times$ such that 
$R'$ is isomorphic to $W(A)[u]/(u^2-\alpha)$. 
Then the assertion follows from the equivalence $(a)\iff(b)$ in Lemma 
\ref{Wunr}.3. 

$4$. We may assume that $A$ is a normal domain. Then the assertion 
follows from the surjectivity of the map 
$\pi_1^{ab}(\eta)\to\pi_1^{ab}({\rm Spec}(A))$ where $\eta$ is 
the generic point of $A$. 
\qed}

\section{Isolated singularities and symmetric bilinear forms}\label{singbil}
\subsection{Reminder on residue symbols}
In this and next subsections, we recall the notion of residue symbol and construct symmetric bilinear forms 
from isolated singular points. 

Residue symbol is introduced in \cite[III 9]{Res} at least for (locally) noetherian schemes. 
In the sequel, we will need the same construction for non-noetherian cases, in which the arguments in \cite{Res} also work. 
However, we will need explicit descriptions of some isomorphisms appearing in this context, in order to  
complete the proofs of the main results of 
this paper. For this reason, we record the arguments in loc.~cit.~in a manner which is suitable for our purposes, but content ourselves with explaining only what is necessary in this paper. Affine cases are treated in \cite{algGr} in detail.

For a ringed space $(X,\mathcal{A})$, ${\rm Qcoh}(\mathcal{A})$ denotes the category of quasi-coherent $\mathcal{A}$-modules. 
We also write ${\rm Qcoh}(X)$ for it if no confusions occur. 
\begin{df}(cf. \cite[III 6]{Res})
Let $Z,S$ be schemes. For a finite morphism $g\colon Z\to S$, define a functor 
\begin{equation*}
g^!\colon {\rm Qcoh}(S)\to {\rm Qcoh}(Z)
\end{equation*}
by 
\begin{equation*}
 \mathcal{F}\mapsto \OO_Z\otimes_{g^{-1}g_\ast\OO_Z}
g^{-1}\mathcal{H}om_{\OO_S}(g_\ast\OO_Z,\mathcal{F}),
\end{equation*}
 where $g^{-1}\mathcal{H}om_{\OO_S}(g_\ast\OO_Z,\mathcal{F})$ is 
regarded as a $g^{-1}g_\ast\OO_Z$-module via the left term of $\mathcal{H}om$. 
\end{df}
When $S={\rm Spec}(A),Z={\rm Spec}(B)$ are affine and $\mathcal{F}$ is associated with 
an $A$-module $M$, $g^!\mathcal{F}$ is the quasi-coherent sheaf associated with the 
$B$-module ${\rm Hom}_A(B,M)$. 
\begin{lm}\label{propg!}
Let $g\colon Z\to S$ be a finite morphism. 
\begin{enumerate}
\item The functor $g_\ast g^!(-)$ is canonically isomorphic to $\mathcal{H}om_{\OO_S}(g_\ast\OO_Z,-)$. 
\item $g^!$ is a right adjoint to $g_\ast$. When $S={\rm Spec}(A),Z={\rm Spec}(B)$ are affine and $\mathcal{F}$ is the associated quasi-coherent sheaf with an $A$-module $M$, the counit corresponds to the evaluation map 
${\rm Hom}_A(B,M)\to M$ at $1\in B$. 
\item Let $\mathcal{M}$ be a quasi-coherent subsheaf of $g^!\OO_S$. Suppose that the composition 
$g_\ast\mathcal{M}\to g_\ast g^!\OO_S\to \OO_S$ is zero where the latter map is the counit of 
the adjunction. Then $\mathcal{M}$ is zero as a quasi-coherent sheaf. 
\end{enumerate}
\end{lm}
\proof{
$1$. The map $\mathcal{H}om_{\OO_S}(g_\ast\OO_Z,\mathcal{F})\to g_\ast g^!\mathcal{F}$ 
 is the composition of 
\begin{equation*}
\mathcal{H}om_{\OO_S}(g_\ast\OO_Z,\mathcal{F})\to g_\ast g^{-1}\mathcal{H}om_{\OO_S}(g_\ast\OO_Z,\mathcal{F})\to
g_\ast(\OO_Z\otimes_{g^{-1}g_\ast\OO_Z} g^{-1}\mathcal{H}om_{\OO_S}(g_\ast\OO_Z,\mathcal{F})). 
\end{equation*}
To check that this is an isomorphism, we may assume that $S$ (hence also $Z$) is affine, in which case it follows from the description 
of $g^!$ for affine cases. 

$2$. Under the identification $g_\ast\colon {\rm Qcoh}(Z)\xrightarrow{\cong}{\rm Qcoh}(g_\ast\OO_Z)$, 
$g^!$ corresponds to $\mathcal{H}om_{\OO_S}(g_\ast\OO_Z,-)$ by $1$ and 
$g_\ast$ corresponds to the forgetful functor. Then the adjunction is well-known. 
The last assertion follows from these descriptions of $g_\ast, g^!$. 

$3$. Let $\mathcal{M}$ be a quasi-coherent subsheaf of $g^!\OO_S$ such that $g_\ast\mathcal{M}\to\OO_S$ is zero. 
Via the adjunction in $2$, the map $g_\ast\mathcal{M}\to\OO_S$ corresponds to the inclusion map $\mathcal{M}\to g^!\OO_S$. 
This is zero if and only if $\mathcal{M}$ itself is zero. 
\qed}

For the computation of $g^!$ in the sequel, it is convenient to record \cite[III 7.2]{Res} in the following manner. 

Let $X$ be a scheme and $\mathcal{E}$ be a locally free sheaf of constant rank $n$ on $X$. 
Let $s\in \Gamma(X,\mathcal{E})$ be a global section. 
Let $Z\subset X$ be the vanishing locus of $s$,~i.e.~the fibered product of the diagram 
\begin{equation*}
\xymatrix{
&X\ar[d]^-s\\X\ar[r]^-0&\mathbb{V}(\mathcal{E})
}
\end{equation*}
where $\mathbb{V}(\mathcal{E})$ is the spectrum of the symmetric $\OO_X$-algebra ${\rm Sym}^\bullet\mathcal{E}^{\vee}$ and $0\colon X\to
\mathbb{V}(\mathcal{E})$ is the $0$-section.
We assume that $s$ gives an $\OO_X$-regular sequence. Namely, taking a trivialization $\mathcal{E}\cong\mathcal{O}_X^{\oplus n}$ Zariski-locally 
by which $s$ corresponds to $(s_i)_i\in\Gamma(X,\OO_X^{\oplus n})$, the sequence $s_1,\dots,s_n$ is $\OO_X$-regular. 
Set $i\colon Z\to X$. 

Write ${\cal I}_Z\in{\cal O}_X$ for the sheaf of ideals defining $Z$. Then the map of 
${\cal O}_X$-algebras 
$s^\ast\colon {\rm Sym}^\bullet{\cal E}^\vee\to{\cal O}_X,{\cal E}^\vee\ni\varphi\mapsto\varphi(s)$ 
induces a surjection ${\rm Sym}^+{\cal E}^\vee\to{\cal I}_Z$ where 
${\rm Sym}^+{\cal E}^\vee$ is the augmentation ideal $\bigoplus_{m\geq1}{\rm Sym}^m{\cal E}^\vee$ 
of ${\rm Sym}^\bullet{\cal E}^\vee$. 
\begin{lm}\label{conE}
\begin{enumerate}
\item 
The surjection ${\rm Sym}^+{\cal E}^\vee\to{\cal I}_Z$ induces an isomorphism 
$i^\ast{\cal E}^\vee\to{\cal I}_Z/{\cal I}_Z^2$ of locally free ${\cal O}_Z$-modules. 
\item Let $({\cal E}_1,s_1)$ be another pair as $({\cal E},s)$ such that the vanishing locus of $s_1$ is $Z$. 
Suppose that an isomorphism $\alpha\colon {\cal E}\to{\cal E}_1$ which sends $s$ to $s_1$ is given. 
Then the diagram 
\begin{equation*}
\xymatrix{
i^\ast{\cal E}_1^\vee\ar[rr]^-{i^\ast\alpha^\vee}\ar[rd]_-{1.}&&i^\ast{\cal E}^\vee\ar[ld]^-{1.}\\
&{\cal I}_Z/{\cal I}_Z^2
}
\end{equation*}
is commutative. 
\end{enumerate}
\end{lm}
\proof{
$1$. 
As $s^\ast\colon{\rm Sym}^\bullet{\cal E}^\vee\to{\cal O}_X$ sends ${\rm Sym}^+{\cal E}^\vee$ onto 
${\cal I}_Z$, we have a surjection  
\begin{equation*}
{\cal E}^\vee\cong{\rm Sym}^+{\cal E}^\vee/
({\rm Sym}^+{\cal E}^\vee)^2\to{\cal I}_Z/{\cal I}_Z^2. 
\end{equation*}
Hence we get a surjection $i^\ast{\cal E}^\vee\to{\cal I}_Z/{\cal I}_Z^2$, which 
is an isomorphism as 
the both sides are locally free ${\cal O}_Z$-modules of the same rank. 

$2$. This is clear from the construction. 
\qed
}

In the situation as above, we define a natural isomorphism (\ref{mapcor}) 
for ${\cal F}\in{\rm Qcoh}(X)$ 
as follows. 

Consider a Koszul complex 
\begin{equation}\label{KosE}
K^\bullet :(\det\mathcal{E})^\vee\xrightarrow{{\rm id}\otimes(s\wedge)}(\det\mathcal{E})^\vee\otimes\mathcal{E}\xrightarrow{
{\rm id}\otimes(s\wedge)}
(\det\mathcal{E})^\vee\otimes\wedge^2\mathcal{E}\xrightarrow{{\rm id}\otimes(s\wedge)}
\cdots\xrightarrow{{\rm id}\otimes(s\wedge)}(\det\mathcal{E})
^\vee\otimes\det\mathcal{E}\cong\OO_X. 
\end{equation}
Here $\det\mathcal{E}$ denotes the highest exterior power 
$\wedge^{n}\mathcal{E}$. As the section 
$s$ gives an $\OO_X$-regular sequence, 
(\ref{KosE}) 
is a locally free resolution of $i_\ast \OO_Z$. Let ${\cal F}$ be a 
quasi-coherent ${\cal O}_X$-module. 
The resolution (\ref{KosE}) defines a quasi-isomorphism 
$R\mathcal{H}om_{\OO_X}(i_\ast\OO_Z,\mathcal{F})\cong{\cal H}om_{{\cal O}_X}(K^\bullet,{\cal F})$. By definition, the $m$-th component of the right-hand side is ${\cal H}om_{{\cal O}_X}(K^{-m},{\cal F})$ and the differential ${\cal H}om_{{\cal O}_X}(K^{-m},{\cal F})\to{\cal H}om_{{\cal O}_X}(K^{-m-1},{\cal F})$ is given by 
$\varphi\mapsto (-1)^{m+1}\varphi d_{K^{-m-1}}$ where we write $d_{K^{-m-1}}\colon K^{-m-1}\to
K^{-m}$. 
Then the canonical map 
\begin{equation}\label{mapperf}
{\cal F}\otimes_{{\cal O}_X}
{\cal H}om_{{\cal O}_X}(K^\bullet,{\cal O}_X)\to 
{\cal H}om_{{\cal O}_X}(K^\bullet,{\cal F})
\end{equation}
 is a morphism of complexes, hence an 
isomorphism of complexes as $K^\bullet$ is a bounded complex of locally finite free ${\cal O}_X$-modules. 

The complex ${\cal H}om_{{\cal O}_X}(K^\bullet,{\cal O}_X)$ is of the form
\begin{equation}\label{KosE*F}
\OO_X\to\det\mathcal{E}\otimes\wedge^{n-1}\mathcal{E}^\vee\to
\cdots\to\det\mathcal{E}
\end{equation}
where $\mathcal{O}_X$ is put on degree $0$. 
 Therefore, the quotient map 
$\det\mathcal{E}[-n]\to i_\ast i^\ast\det\mathcal{E}[-n]$ and the inverse of (\ref{mapperf}) give us 
a map 
\begin{align}\label{mapcor}
R\mathcal{H}om_{\OO_X}(i_\ast\OO_Z,\mathcal{F})
\xrightarrow{(\ref{mapperf})}&{\cal F}\otimes_{{\cal O}_X}{\cal H}om_{{\cal O}_X}(K^\bullet,{\cal O}_X)\\
\notag&\to
 \mathcal{F}\otimes_{\OO_X}^Li_\ast i^\ast\det\mathcal{E}[-n]
\xrightarrow{\cong} i_\ast(Li^\ast\mathcal{F}\otimes_{\OO_Z}i^\ast\det\mathcal{E})[-n]
\end{align}
in the derived category. Here we write $\otimes_{{\cal O}_Z}i^\ast \det{\cal E}$ instead of $\otimes^L_{{\cal O}_Z}i^\ast\det{\cal E}$ 
in the right-hand side since $i^\ast\det{\cal E}$ is a locally free ${\cal O}_Z$-module. 
\begin{lm}\label{resKos}
\begin{enumerate}
\item 
In the situation as above, The map (\ref{mapcor}) is a quasi-isomorphism. 
\item Let $({\cal E}_1,s_1)$ be as in Lemma \ref{conE}.2 with an isomorphism $\alpha\colon{\cal E}\to{\cal E}_1$ sending $s$ to $s_1$. Then 
the diagram 
\begin{equation*}
\xymatrix{
&R\mathcal{H}om_{\OO_X}(i_\ast\OO_Z,\mathcal{F})\ar[ld]_-{(\ref{mapcor})}\ar[rd]^-{(\ref{mapcor})}&\\
i_\ast(Li^\ast\mathcal{F}
\otimes_{\OO_Z}i^\ast\det{\cal E})[-n]\ar[rr]&&i_\ast(Li^\ast\mathcal{F}
\otimes_{\OO_Z}i^\ast\det{\cal E}_1)[-n]
}
\end{equation*}
is commutative. Here the horizontal arrow is induced from $\det\alpha$. 
\item{(\cite[III 7.3]{Res})} Let $j\colon W\to Y$ be a regular closed immersion  
of schemes of  codimension $n$. For a quasi-coherent ${\cal O}_Y$-module ${\cal F}$ which 
is $Lj^\ast$-acyclic, 
we have a natural quasi-isomorphism 
\begin{equation}\label{isomKos}
R\mathcal{H}om_{\OO_Y}(j_\ast\OO_W,\mathcal{F})\cong j_\ast(j^\ast\mathcal{F}
\otimes_{\OO_W}\omega_{W/Y})[-n]
\end{equation}
with the following properties. 
Here $\omega_{W/Y}$ is the determinant of the conormal sheaf. 
\begin{enumerate}

\item Suppose that we are given a locally free sheaf ${\cal E}$ of rank $n$ on $Y$ and 
a section $s$ of $\cal E$ whose vanishing locus is $W$. Then, we have a commutative diagram 
\begin{equation*}
\xymatrix{
&R\mathcal{H}om_{\OO_Y}(j_\ast\OO_W,\mathcal{F})\ar[ld]_-{(\ref{mapcor})}\ar[rd]^-{(\ref{isomKos})}&\\
j_\ast(j^\ast\mathcal{F}
\otimes_{\OO_W}j^\ast\det{\cal E})[-n]&&j_\ast(j^\ast\mathcal{F}
\otimes_{\OO_W}\omega_{W/Y})[-n]\ar[ll]
}
\end{equation*}
where the horizontal arrow is the one induced from the isomorphism $j^\ast{\cal E}^\vee\to 
{\cal I}_W/{\cal I}_W^2$ in Lemma \ref{conE}. 
\item For an open immersion $U\hookrightarrow Y$, the restriction of (\ref{isomKos}) to $U$ is 
equal to the one defined for ${\cal F}|_U$ on $U$. 
\item Further assume that ${\cal F}$ is ${\cal O}_Y$-flat. 
Let $f\colon Y'\to Y$ be a morphism of schemes such that $W':=W\times_YY'\to Y'$ is a 
regular immersion of codimension $n$. Set $j'\colon W'\to Y'$ and $f_W\colon W'\to W$. 
The diagram 
\begin{equation}\label{daigreg}
\xymatrix{
Lf^\ast R{\cal H}om_{{\cal O}_Y}(j_\ast{\cal O}_W,{\cal F})\ar[r]^-{(\ref{isomKos})}\ar[d]&
Lf^\ast j_\ast(j^\ast{\cal F}\otimes_{{\cal O}_W}\omega_{W/Y})[-n]\ar[d]\\
R{\cal H}om_{{\cal O}_{Y'}}(j'_\ast{\cal O}_{W'},f^\ast{\cal F})\ar[r]^-{(\ref{isomKos})}&
 j'_\ast(j'^\ast f^\ast{\cal F}\otimes_{{\cal O}_{W'}}\omega_{W'/Y'})[-n]
}
\end{equation}
is commutative. 
\end{enumerate}
\end{enumerate}
\end{lm}
The flatness assumption in $3$ is apparently superfluous, and indeed one can construct (\ref{isomKos}) 
which satisfies $(a),(b)$ for general quasi-coherent sheaves, which is done in \cite[III, 7]{Res}. 
However, it seems difficult to verify the property $(c)$ in this full generality, for which 
the author could not find a reference. 
\proof{
$1$. The only non-trivial part is that the middle arrow in (\ref{mapcor}) is a quasi-isomorphism. 
Since (\ref{KosE*F}) is a locally free resolution of 
$i_\ast i^\ast\det\mathcal{E}$ as $\OO_X$-modules, the assertion follows. 

$2$. The exterior powers of $\alpha$ induce an isomorphism between the Koszul complexes (\ref{KosE}) 
defined from $({\cal E},s),({\cal E}_1,s_1)$. The assertion is verified directly from the 
construction of (\ref{mapcor}). 

$3$. 
Let $\cal F$ be a quasi-coherent ${\cal O}_Y$-module which is $Lj^\ast$-acyclic. Locally on $Y$, 
we can take $({\cal E},s)$ as in $(a)$. Then, by $1$, we have a quasi-isomorphism 
$R\mathcal{H}om_{\OO_Y}(j_\ast\OO_W,\mathcal{F})\to j_\ast(j^\ast\mathcal{F}
\otimes_{\OO_W}j^\ast\det{\cal E})[-n].$ Composing $j^\ast\det{\cal E}\xrightarrow{\cong}
\omega_{W/Y}$ defined from Lemma \ref{conE}.1, we have a quasi-isomorphism 
$R\mathcal{H}om_{\OO_Y}(j_\ast\OO_W,\mathcal{F})\to j_\ast(j^\ast\mathcal{F}
\otimes_{\OO_W}\omega_{W/Y})[-n].$ By Lemma \ref{conE}.2 and 2, this isomorphism 
is independent of the choice of $({\cal E},s)$. Since the complexes under consideration are 
concentrated on one degree, the isomorphisms glue to an isomorphism on the whole of $Y$, which 
satisfies the properties $(a),(b)$ by the construction. 

We verify the property $(c)$. We show the commutativity of (\ref{daigreg}). By the 
assumption on $f$, the canonical map $Lf^\ast j_\ast\to j'_\ast Lf_W^\ast$ is a quasi-isomorphism. 
 Hence 
the right vertical arrow of (\ref{daigreg}) is a quasi-isomorphism. Also 
the horizontal arrows are quasi-isomorphisms. 
Consequently 
the complexes in the diagram are concentrated on degree $n$. 
Therefore we may replace $Y$ and $Y'$ 
by arbitrary open coverings to check the commutativity. Hence we may assume that 
we are given a pair $({\cal E},s)$ on $Y$ as in $(a)$. Then, the assertion follows from the 
compatibility of the construction of (\ref{mapcor}). 
\qed}

Let $S$ be a scheme and $f\colon X\to S$ be a smooth morphism of schemes purely of relative dimension $n$. 
Let $\mathcal{E}$ be a locally free $\OO_X$-module of rank $n$ on $X$ and $s\in\Gamma(X,\mathcal{E})$ be a global section with the 
vanishing locus $Z$. 
Let $i\colon Z\to X$ be the 
closed immersion and $g=f\circ i\colon Z\to S$ be the structure morphism. 
We apply Lemma \ref{resKos}.3 to describe $g^!$. 
\begin{lm}\label{cg!}
Assume that $Z$ is {\it quasi-finite} over $S$. 
\begin{enumerate}
\item The immersion $i\colon Z\to X$ is transversally regular of codimension $n$ relative to $S$ (\cite[(19.2.2)]{EGA4}). 
$Z$ is flat of finite presentation over $S$. 
\item Assume further that $Z$ is finite over $S$. Then there is a canonical isomorphism  
\begin{equation}\label{idg!}
\Psi_s\colon g^!\to g^\ast(-)\otimes i^\ast(\omega_{X/S}\otimes\det\mathcal{E}) 
\end{equation}
of functors. 
In particular, 
the counit ${\rm Tr}_g\colon g_\ast g^!\mathcal{O}_S\to\mathcal{O}_S$ is canonically identified with a map 
$g_\ast i^\ast(\omega_{X/S}\otimes\det\mathcal{E})\to
\mathcal{O}_S$. This map is denoted by ${\rm Res}[\frac{}{s}]$, or 
${\rm Res}_S[\frac{}{s}]$ if we want to specify the base scheme, and 
called the residue symbol (cf. \cite[III 9]{Res}). 
\item Let $\mathcal{E}_1$ be another locally free sheaf of rank $n$ on $X$ with a global section $s_1$ whose vanishing locus equals to 
$Z$. For an isomorphism $\alpha\colon \mathcal{E}\to\mathcal{E}_1$ which sends $s$ to $s_1$, the composition of 
\begin{equation*}
i^\ast(\omega_{X/S}\otimes\det\mathcal{E})\xrightarrow{\Psi_s^{-1}} g^!\OO_S\xrightarrow{
\Psi_{s_1}} i^\ast(\omega_{X/S}\otimes\det\mathcal{E}_1)
\end{equation*}
 equals to $i^\ast({\rm id}\otimes\det\alpha)$. 
\end{enumerate}
\end{lm}
\proof{
$1$. Since the assertion is local on $X$, we may assume that 
$S={\rm Spec}(A),X={\rm Spec}(B)$ are affine and that 
$\mathcal{E}\cong\OO_X^{\oplus n}$ is free. Let $s_j\in 
B$ be the $j$-th coordinate of $s$. We need to show that $(s_j)_j$ is a $B$-regular sequence around $Z$ and that 
$B/\sum_js_jB$ is $A$-flat. 

By a limit argument, we may assume that $A$ and $B$ are noetherian. In this case, \cite[Corollary to Theorem 22.5]{Mat} is applicable. Using this, we reduce the assertion to the case where $A$ is a field. Take 
$x\in Z\subset{\rm Spec}(B)$. By the assumption that $Z$ is finite over $S$, 
the local ring $B_x/\sum_js_jB_x$ is $0$-dimensional, which implies that $s_1,\dots,s_n$ is a system of parameters. 
As $B$ is regular, hence Cohen-Macaulay, $(s_j)_j$ is a $B$-regular sequence.

$2$. The proof goes similarly as in \cite{Res}. Since we need a concrete description for $g^!\OO_S$, we give a proof for 
this case. The general case can be deduced from this particular case, as is explained in \cite[5.4]{NeeG}. 

We define a canonical isomorphism $g^!\OO_S\cong i^\ast(\omega_{X/S}\otimes_{\OO_X}\det\mathcal{E})$. 
Consider a commutative diagram 
\begin{equation*}
\xymatrix{
Z\ar[rd]^{\rm id}\ar@{^{(}->}[r]^{\Gamma_i\ \ }&Z\times_SX\ar[d]^{{\rm pr}_Z}\ar[r]^{\ \ {\rm pr}_X}&X\ar[d]^f\\
& Z\ar[r]^g&S 
}
\end{equation*}
where the square is cartesian. 
As ${\rm pr}_Z^\ast g^!{\cal O}_Z$ is $L\Gamma_i^\ast$-acyclic, we can apply Lemma \ref{resKos}.3 to $(\Gamma_i, {\rm pr}_Z^\ast g^!{\cal O}_Z)$ to get 
an isomorphism 
\begin{equation}\label{isom}
R\mathcal{H}om_{\OO_{Z\times X}}(\Gamma_{i\ast}\OO_Z,{\rm pr}_Z^\ast g^!\OO_S)\cong 
\Gamma_{i\ast}(g^!\OO_S\otimes_{\OO_Z} i^\ast
\omega_{X/S}^{-1})[-n]. 
\end{equation}
Here we use a canonical isomorphism $\omega_{Z/Z\times_SX}\cong i^\ast\omega_{X/S}^{-1}$. 
From this, we compute 
\begin{align*}
\Gamma_{i\ast}(g^!\OO_S)&\cong R\mathcal{H}om_{\OO_{Z\times_S X}}(\Gamma_{i\ast}\OO_Z,{\rm pr}_Z^\ast g^!\OO_S\otimes{\rm pr}_X^\ast\omega_{X/S}[n])\\
&\cong R\mathcal{H}om_{\OO_{Z\times_S X}}(\Gamma_{i\ast}\OO_Z,
{\rm pr}_X^!\omega_{X/S}[n]). 
\end{align*}
Taking ${\rm pr}_{X\ast}$, we have 
\begin{align}
i_\ast g^!\OO_S&\cong {\rm pr}_{X\ast}R\mathcal{H}om_{\OO_{Z\times_S X}}(\Gamma_{i\ast}\OO_Z,{\rm pr}_X^!\omega_{X/S}[n])\notag\\
& \cong R\mathcal{H}om_{\OO_X}(i_\ast\OO_Z,\omega_{X/S}[n])\notag\\
&\cong i_\ast i^\ast(\omega_{X/S}\otimes\det\mathcal{E}) \label{isom2} 
\end{align}
where the second isomorphism is the derived version of the adjunctions and the last one is due to Lemma \ref{resKos}.1. 

$3$. The exterior powers of $\alpha$ induce an isomorphism of Koszul complexes (\ref{KosE}). The assertion follows from the construction in $2$ and Lemma \ref{resKos}.2. 
\qed}

We collect some of functorialities of the construction. 
\begin{lm}\label{bcvan}
Let $h\colon S'\to S$ be a morphism of schemes. 
Let $X',\mathcal{E}'$ be the pullbacks of $X,\mathcal{E}$ to $S'$. Let $s'\in\Gamma(X',\mathcal{E}')$ be the pullback of $s$. 
\begin{enumerate}
\item The vanishing locus $Z'$ of $s'$ is the pullback of $Z$ by $X'\to X$. 
\item Assume that $Z$ is finite over $S$. Then the coherent sheaf 
$g^!{\cal O}_S$ is ${\cal O}_Z$-flat. 
\item Assume that $Z$ is {\it finite} over $S$. 
Let $g'\colon Z'\to S'$ be the structure map and $h_Z\colon Z'\to Z$ be the projection. Let $i'\colon Z'\to X'$ be the immersion. 
The base change map $h_Z^\ast g^!\to g'^!h^\ast$ fits into a commutative diagram 
\begin{equation*}
\xymatrix{
h_Z^\ast g^!\ar[r]\ar[d]&g'^!h^\ast\ar[d]\\
h_Z^\ast (g^\ast(-)\otimes i^\ast(\omega_{X/S}\otimes\det\mathcal{E}))\ar[r]&
g'^\ast h^\ast(-)\otimes i'^\ast(\omega_{X'/S'}\otimes\det\mathcal{E'}). 
}
\end{equation*}
Here the vertical arrows are the ones constructed in Lemma \ref{cg!}.2. 
\item Assume that $Z$ is {\it finite} over $S$. The diagram 
\begin{equation*}
\xymatrix{
h^\ast g_\ast i^\ast(\omega_{X/S}\otimes\det\mathcal{E})\ar[r]^{\ \ \ \ \ h^\ast{\rm Res}[\frac{}{s}]}\ar[d]&\OO_{S'}\\
g'_\ast h_Z^\ast i^\ast(\omega_{X/S}\otimes\det\mathcal{E})\ar[r]&g'_\ast i'^\ast(\omega_{X'/S'}\otimes\det\mathcal{E}')\ar[u]_{{\rm Res}
[\frac{}{s'}]}
}
\end{equation*}
is commutative. 
\end{enumerate} 
\end{lm}
\proof{
The assertion $1$ is clear from the definition. 

To apply Lemma \ref{resKos}.3.$(c)$ to the assertion $3$, 
we need to prove $2$, which is verified as follows. 
 By a limit argument, we may assume that $S$ is noetherian. By the local criterion of flatness, 
we further reduce it to the case where $S$ is the spectrum of an algebraically closed field $k$. 
Shrinking $X$, we also assume that $Z$ is local. 
In this case, ${\cal O}_Z$ is a Gorenstein local ring. In particular, 
${\cal O}_Z$ has a unique non-zero minimal ideal. 
Dually, $g^!{\cal O}_Z={\cal H}om_k({\cal O}_Z,k)$ has a unique ${\cal O}_Z$-linear surjection 
$g^!{\cal O}_Z\to k$ up to scalar multiplication. This implies that $g^!{\cal O}_Z\otimes_{{\cal O}_Z}k$ is one-dimensional. 
Hence there exists an ${\cal O}_Z$-linear surjection ${\cal O}_Z\to g^!{\cal O}_Z$, which 
must be an isomorphism as they have the same lengths. The assertion $2$ follows. 

We show $3$. 
Lemma \ref{resKos}.3.$(c)$ ensures the commutativity of the first isomorphism in 
(\ref{isom2}) with base change. For the second one, it follows from the compatibility of the adjunction 
with base change, which can be verified directly in this case. The commutativity of 
the last one follows from the construction of (\ref{mapcor}).

The assertion $4$ follows from $3$ and the functoriality of the counits. 
\qed}
\begin{lm}\label{Trfet}
Let $f\colon X\to S'$ be a smooth morphism of schemes purely of relative dimension $n$ 
and let $\mathcal{E}$ be a locally free $\OO_X$-module of rank $n$. Let $s$ be a global section of $\mathcal
{E}$ whose vanishing locus $Z$ is finite over $S'$. 

Let $h\colon S'\to S$ be a finite \'etale morphism. Then, 
the residue symbol ${\rm Res}_S[\frac{}{s}]$ equals to ${\rm Tr}_{S'/S}\circ{\rm Res}_{S'}[\frac{}{s}]$. 
\end{lm}
\proof{
Regard $X$ as a smooth scheme over $S$ and consider a commutative diagram 
\begin{equation*}
\xymatrix{
Z\ar@{^{(}->}[r]\ar[rd]_{\rm id}&Z'\ar[r]\ar[d]^{h_Z}&X'\ar[r]\ar[d]^{h_X}&S'\ar[d]^h\\
&Z\ar[r]&X\ar[r]&S
}
\end{equation*}
where the squares are cartesian. Write $g\colon Z\to S$, $g'\colon Z'\to S'$. 
We have a commutative diagram 
\begin{equation*}
\xymatrix{
h_{Z\ast}g'^!\OO_{S'}\ar[r]^\cong\ar[d]_{{\rm Tr}_{g'}}&h_{Z\ast}h_Z^\ast g^!\OO_S\ar[r]^{{\rm Tr}_{Z'/Z}}&g^!\OO_S\ar[d]^{{\rm Tr}_g}\\
h_{\ast}\OO_{S'}\ar[rr]^{{\rm Tr}_{S'/S}}&&\OO_S. 
}
\end{equation*}
Applying Lemma \ref{bcvan}.3, 
the statement follows by restricting $g'^!\OO_{S'}$ to the open and closed subscheme $Z$ of $Z'$. 
\qed}

\begin{lm}\label{prodtr}
Let $X_j\ (j=1,2)$ be smooth $S$-schemes purely of relative dimensions $n_j$ and let $\mathcal{E}_j$ be locally free sheaves of rank $n_j$ on $X_j$. 
 Let $s_j$ be global sections of $\mathcal{E}_j$ whose vanishing loci $Z_j$ are finite over $S$. 
 Set $g_j\colon Z_j\to S$ and $i_j\colon Z_j\to X_j$. 
 Let 
 $X:=X_1\times_SX_2\xrightarrow{{\rm pr}_j}X_j$ be the projections from the fibered product. 
 
 Let $\mathcal{E}:=\mathcal{E}_1\boxtimes_{\OO_S}\mathcal{E}_2={\rm pr}_1^\ast\mathcal{E}_1\oplus 
 {\rm pr}_2^\ast\mathcal{E}_2$ be the external product on $X$. Write $s$ for ${\rm pr}_1^\ast s_1+{\rm pr}_2^\ast s_2\in\Gamma(X,\mathcal{E})$. 
 Then the vanishing locus of $s$ equals to $Z:=Z_1\times_SZ_2$. We write $i\colon Z\to X$ for the inclusion. 
 Under the identification 
 $i^\ast(\omega_{X/S}\otimes\det\mathcal{E})\cong 
 {\rm pr}^\ast_1i_1^\ast(\omega_{X_1/S}\otimes\det\mathcal{E}_1)\otimes {\rm pr}^\ast_2i_2^\ast(\omega_{X_2/S}\otimes\det\mathcal{E}_2)
 $, ${\rm Res}[\frac{}{s}]$ corresponds to ${\rm Res}[\frac{}{s_1}]\otimes{\rm Res}[\frac{}{s_2}]$. 
 \end{lm}
 \proof{ 
 The first assertion on the vanishing locus is clear, as the ideal sheaf of $Z$ is generated by those of $Z_1,Z_2$. 
 
 Let $g\colon Z\to S$. For the second assertion, we have a diagram 
 \begin{equation*}
 \xymatrix{
g_\ast g^!\OO_S\ar[dd]_{c}\ar[r]^{a}&g_\ast{\rm pr}_1^!g_1^!\OO_S\ar[r]^{b}
\ar[d]_{d}&g_{1\ast}g_1^!\OO_S\ar[dd]_h\\
&g_\ast( {\rm pr}^\ast_2i_2^\ast(\omega_{X_2/S}\otimes\det\mathcal{E}_2)
\otimes_{\OO_Z}{\rm pr}_1^\ast g_1^!\OO_S)\ar[d]_e\ar[ru]^f&\\
g_\ast i^\ast(\omega_{X/S}\otimes\det\mathcal{E})\ar[r]^{k\ \ \ \ \ \ \ \ \ \ \ \ \ \ \ \ \ \ \ \ \ }
&g_\ast( {\rm pr}^\ast_2i_2^\ast(\omega_{X_2/S}\otimes\det\mathcal{E}_2)\otimes_{\OO_Z}{\rm pr}^\ast_1i_1^\ast(\omega_{X_1/S}\otimes\det\mathcal{E}_1))
\ar[r]^{\ \ \ \ \ \ \ \ \ \ \ \ \ \ \ \ \ \ \ \ \ \ \ \ \ \ \ \ l}&\OO_S. 
 }
 \end{equation*}
Here $a,k$ are the canonical ones, $b,h$ are the counits of adjunctions, $c,d,e$ are the isomorphisms in (\ref{idg!}). 
$f$ is induced from ${\rm Res}[\frac{}{{\rm pr}_2^\ast s_2}]$. $l$ is ${\rm Res}[\frac{}{s_2}]\otimes{\rm Res}[\frac{}{s_1}]$. 
It is straightforward to check that it is commutative and that the map $g_\ast g^!\OO_S\to\OO_S$ from the upper left to 
the lower right is the counit. 
\qed}

At the end of this preliminary subsection, 
we give a way to compute 
the isomorphism (\ref{idg!}) in Lemma \ref{cg!}. 
To make the notation compatible with the one in the subsection 3.3, let us suppose that
 $X$ is a smooth $S$-scheme purely of relative dimension $n+1$. Let ${\cal E}$ be a locally 
free ${\cal O}_X$-module of rank $n+1$. Suppose that we are given a global section $s$ of 
${\cal E}$ whose vanishing locus $Z$ is finite over $S$. Further assume that $X$ is equipped with 
local parameters $t_0,\dots,t_n$. 

In this setting, 
we define a map $\Phi_s\colon g_\ast i^\ast(\omega_{X/S}\otimes_{{\cal O}_X}\det{\cal E})\to g_\ast g^!{\cal O}_S$ as follows. 
Let $K^\bullet(\underline{\xi})$ be the Koszul complex defined by $\xi_i=1\otimes t_i-
t_i\otimes1\in\OO_{Z\times_SX}$. We regard $\OO_{Z\times_S X}$ as a right $\OO_X$-module by the right component of the tensor product. 
Take a morphism of resolutions of ${\cal O}_Z$ by locally free 
$\OO_X$-modules 
\begin{equation}\label{morKos}
\xymatrix{
K^\bullet(\underline{\xi})\ar@{}[r]|:&\OO_{Z\times_S X}\ar[r]\ar[d]_{\eta}&\OO_{Z\times_S X}^{n+1}\ar[d]
\ar[r]&\cdots\ar[r]\ar[d]&\OO_{Z\times_S X}\ar[d]^{\zeta}\\
(\ref{KosE})\ar@{}[r]|:&(\det{\cal E})^\vee\ar[r]&(\det{\cal E})^\vee\otimes{\cal E}\ar[r]&\cdots\ar[r]&{\cal O}_X. 
}
\end{equation}
 Namely, it is a morphism of complexes such that the composition ${\cal O}_{Z\times_SX}\xrightarrow{\zeta}{\cal O}_X\to{\cal O}_Z$ coincides with the quotient map. 

We identify the $\OO_S$-sheaves 
$g_\ast i^\ast(\omega_{X/S}\otimes\det{\cal E})$, $ g_\ast{\cal H}om_{{\cal O}_Z}(i^\ast 
(\det{\cal E})^\vee,i^\ast\omega_{X/S})$ under the isomorphism  $g_\ast i^\ast(\omega_{X/S}\otimes\det{\cal E})\ni x\otimes y\mapsto(\varphi\mapsto\varphi(y)x)$. 
For an element $\bar{\alpha}\in g_\ast i^\ast(\omega_{X/S}\otimes\det{\cal E})$, take an 
${\cal O}_X$-linear map $\alpha\colon (\det{\cal E})^\vee\to\omega_{X/S}$ which is a lift of $\bar{\alpha}
\in  g_\ast i^\ast(\omega_{X/S}\otimes\det{\cal E})\cong g_\ast{\cal H}om_{{\cal O}_Z}(i^\ast 
(\det{\cal E})^\vee,i^\ast\omega_{X/S})$. 

Composing $\alpha$ and $\eta$ in (\ref{morKos}), we get a map $\alpha\circ\eta\colon 
{\cal O}_{Z\times_SX}\to\omega_{X/S}$. Trivializing $\omega_{X/S}$ by $1\mapsto\wedge_idt_i$, 
we identify this map with an element in ${\cal H}om_{{\cal O}_X}
({\cal O}_{Z\times_SX},{\cal O}_X)$, for which we use the same symbol $\alpha\circ\eta$.

Under the canonical isomorphism 
\begin{equation*}
\mathcal{H}om_{\OO_X}(\OO_{Z\times_S X},\OO_X)\cong 
\mathcal{H}om_{\OO_S}(\OO_{Z},\OO_S)\otimes_{\OO_S}\OO_X, 
\end{equation*}
suppose that the map $ \alpha\circ\eta$ in the left-hand side 
corresponds to $\sum_i\phi_i\otimes a_i$ in the right-hand side. 
Then we get an element $\Phi_s(\alpha)\in g_\ast g^!\OO_S={\cal H}om_{{\cal O}_S}({\cal O}_Z,{\cal O}_S)$ 
by setting $\Phi_s(\alpha)=\sum_i\bar{a}_i\phi_i=\sum_i\phi_i(\bar{a}_i-)$ 
  where $\bar{a}_i$ is the image 
 of $a_i$ in $\OO_Z$. 
 \begin{lm}\label{computePsi}
 The map $\Phi_s\colon g_\ast i^\ast(\omega_{X/S}\otimes\det{\cal E})\to g_\ast g^!{\cal O}_S$ 
 constructed above coincides with $g_\ast\Psi_s^{-1}$ in (\ref{idg!}). 
 \end{lm}
 \proof{
By Lemma \ref{resKos}.3, the isomorphism (\ref{isom}) in Lemma  \ref{cg!} is identified with 
 the one (\ref{mapcor}) constructed for $({\cal E},s)=({\cal O}_{Z\times_SX}^{\oplus (n+1)},(1\otimes t_i-t_i\otimes1)_i)$, 
 under the identification $\det{\cal E}^\vee={\cal O}_{Z\times_SX}\to
 \omega_{X/S},1\mapsto\wedge_idt_i$. Then the construction of $\Phi_s$ is just a restatement of the construction in Lemma 
 \ref{cg!}.2. 
 \qed
 }
 
\subsection{Bilinear forms of isolated singularities}\label{cbil}
In this subsection, we construct non-degenerate symmetric bilinear forms from isolated singularities. 

First let us give the notion of a family of isolated singularities in a manner which suits our purposes. 
\begin{df}\label{isoldef}
 {\rm A family of isolated singularities} is a commutative diagram of schemes
\begin{equation*}
\xymatrix{
Z\ar@{^{(}->}[r]^i\ar[rrd]_g&X\ar[rr]^f\ar[rd]&&C\ar[ld]\\
&&S&
}
\end{equation*}
where $X$ is smooth over $S$, $C$ is a smooth $S$-curve, and 
$Z$ is the singular locus of $f$ which is assumed to be quasi-finite over $S$. 
\end{df}

Let $S$ be a scheme and $X$ be a smooth scheme purely of relative dimension $n$ over $S$. Consider a commutative diagram of $S$-schemes 
\begin{equation*}
\xymatrix{
X\ar[rr]^f\ar[rd]&&C\ar[ld]\\
&S&
}
\end{equation*}
where $C$ is a smooth $S$-curve. The map 
$f^\ast\colon f^\ast\Omega^1_{C/S}\to\Omega^1_{X/S}$ 
defines a section of the locally free ${\cal O}_X$-module $\mathcal{H}om(f^\ast\Omega^1_{C/S},\Omega^1_{X/S})$, which we also denote by  $f^\ast$. 
Recall that the singular locus $Z$ of $f$ is the closed subscheme of $X$ defined by vanishing of $f^\ast$. 

We record Lemma \ref{cg!} again 
in the form which fits into the current situation. 
\begin{lm}\label{cg!now}
Assume that $Z$ is quasi-finite over $S$. 
\begin{enumerate}
\item $Z$ is flat of finite presentation over $S$. 
\item Assume further that $Z$ is finite over $S$. 
Then, the $!$-pull-back $g^!$ on quasi-coherent sheaves 
is canonically isomorphic to 
$g^\ast(-)\otimes i^\ast(\omega^{\otimes2}_{X/S}\otimes f^\ast\omega_{C/S}^{\otimes(-n)})$. 
\item Assume further that $Z$ is finite over $S$ and that 
$C$ admits an everywhere non-zero differential $\omega$. Then the isomorphism 
\begin{equation*}\label{idg!now}
g^!\cong g^\ast(-)\otimes i^\ast\omega_{X/S}^{\otimes2}  
\end{equation*}
which is given by applying Lemma \ref{cg!}.2 to $({\cal E},s)=(\Omega^1_{X/S},f^\ast\omega)$ coincides with 
the isomorphism in $2$ followed by the trivialization $f^\ast\omega_{X/S}^{\otimes(-n)}\cong{\cal O}_X, 
f^\ast\omega^{\otimes(-n)}\mapsto 1$. 
\end{enumerate}
\end{lm}
\proof{
For $1$, $2$, 
apply Lemma \ref{cg!} to $(\mathcal{E},s)=(\mathcal{H}om(f^\ast\Omega^1_{C/S},\Omega^1_{X/S}),f^\ast)$. 
The assertion $3$ follows from Lemma \ref{cg!}.3. 
\qed}

From now on, we assume that $Z$ is {\it finite} over $S$ and that an everywhere non-zero 
differential $\omega\in\Omega^1_{C/S}$ is given. 
In this case, the residue symbol 
\begin{equation*}
{\rm Res}[\frac{}{f^\ast\omega}]\colon g_\ast i^\ast\omega_{X/S}^{\otimes2}\to\OO_S
\end{equation*}
is given by the composition of $g_\ast i^\ast\omega_{X/S}^{\otimes2}\cong 
g_\ast g^!{\cal O}_S\xrightarrow{
{\rm Tr}_g}{\cal O}_S$ where the first isomorphism is the one in Lemma \ref{cg!now}.3. 

Set $\varphi_{f}:=g_\ast i^\ast\omega_{X/S}$. Let us write $B_{f,\omega}$ for the composition of 
\begin{equation*}
\varphi_{f}\otimes_{\OO_S}\varphi_{f}
\to g_\ast (i^\ast\omega_{X/S}
\otimes_{\OO_Z} i^\ast\omega_{X/S})=g_\ast i^\ast\omega_{X/S}^{\otimes2}\xrightarrow{{\rm Res}[\frac{}{f^\ast\omega}]}
\OO_S. 
\end{equation*}
This defines an $\OO_S$-bilinear form on the coherent sheaf $\varphi_{f}$ on $S$. 
\begin{df}\label{phiB}
Assume that $Z$ is finite over $S$.  
{\rm The bilinear form attached to} $(f,\omega)$ is the pair $(\varphi_{f},B_{f,\omega})$ constructed above. 
We also write $(\varphi_{f,S},B_{f,\omega,S})$ for it to specify the base scheme. 
\end{df}

\begin{lm}\label{bpofbil}
Let the notation be as above. Assume that $Z$ is $S$-finite. 
\begin{enumerate}
\item The $\OO_S$-module $\varphi_{f}$ is locally free with the rank equal to ${\rm rk}_{\OO_S}\OO_Z$. 
\item $B_{f,\omega}$ is symmetric and $\OO_Z$-equivariant,~i.e.~we have $B_{f,\omega}(\alpha x,y)=B_{f,\omega}(x,\alpha y)$ for 
$\alpha\in g_\ast\OO_Z$ and $x,y\in\varphi_{f}$. 
\item $B_{f,\omega}$ is non-degenerate, i.e. the induced map $\varphi_{f}\to\varphi_{f}^\vee$ is an isomorphism. 
\item For a unit $\alpha\in\Gamma(C,\OO_C)$, we have $B_{f,\alpha\omega}(-,-)=B_{f,\omega}(-,\alpha^{n}-)$, where $n$ is the 
relative dimension of $X/S$.  
\item Let $h\colon S'\to S$ be a morphism of schemes. Let $f',\omega'$ be the pullbacks of $f,\omega$ by $h$. Then we have a 
canonical isomorphism 
\begin{equation*}
h^\ast(\varphi_{f},B_{f,\omega})\cong(\varphi_{f'},B_{f',\omega'}). 
\end{equation*}
\end{enumerate}
\end{lm}
\proof{

$1$. It follows since $i^\ast\omega_{X/S}$ is an invertible $\OO_Z$-module and $g\colon Z\to S$ is finite flat of finite presentation 
(Lemma \ref{cg!now}.1). 

$2$. Note that, since $i^\ast\omega_{X/S}$ is of rank $1$, the map 
\begin{equation*}
i^\ast\omega_{X/S}\otimes_{\OO_Z}i^\ast \omega_{X/S}\to i^\ast\omega_{X/S}\otimes_{\OO_Z} i^\ast\omega_{X/S}
\end{equation*}
 defined by switching the components 
$x\otimes y\mapsto y\otimes x$ equals to the identity. The assertion follows from this observation. 

$3$. Let 
$x\in\varphi_{f}$ be an element in the kernel 
of $\varphi_{f}\to\varphi_{f}^\vee$. Let $\mathcal{M}' :=\OO_Z\cdot x$ be 
the quasi-coherent subsheaf of $i^\ast\omega_{X/S}$ generated by $x$ and let $\mathcal{M}:=\mathcal{M}'\otimes_{\OO_Z}
i^\ast\omega_{X/S}$, which is a quasi-coherent subsheaf of $i^\ast\omega_{X/S}^{\otimes2}\cong 
g^!\OO_S$. By the assumption, 
we have ${\rm Res}[\frac{x\otimes y}{f^\ast\omega}]=
B_{f,\omega}(x,y)=0$ for any element $y\in i^\ast\omega_{X/S}$.  Since $x\otimes y$ generates $g_\ast(
\mathcal{M}\otimes_{\OO_Z} i^\ast\omega_{X/S})$, 
the composition 
\begin{equation*}
g_\ast(\mathcal{M}\otimes_{\OO_Z} i^\ast\omega_{X/S})\to g_\ast g^!\OO_S\xrightarrow{{\rm Tr}_g}
\OO_S
\end{equation*} 
is zero. Thus, by Lemma \ref{propg!}.3, we know that $\mathcal{M}'\otimes_{\OO_Z}
i^\ast\omega_{X/S}$ is zero, which implies $x=0$ since $i^\ast\omega_{X/S}$ is invertible. 

Therefore we know that the map 
$\varphi_{f}\to\varphi_{f}^\vee$ induced by $B_{f,\omega}$ 
is injective. Since this holds after any base change $S'\to S$, the map is an isomorphism. 

$4$. It follows from Lemma \ref{cg!now}.3. 

$5$. It follows from Lemmas \ref{bcvan}.3, 4. 
\qed}

\begin{lm}\label{etres}
Consider a family of isolated singularities 
\begin{equation*}
\xymatrix{
X\ar[rr]^f\ar[rd]&&C\ar[ld]\\
&S'
}
\end{equation*}
with the singular locus $Z$ finite over $S'$. Let $h\colon S'\to S$ be a finite \'etale morphism. Then the 
bilinear form $(\varphi_{f,S},B_{f,\omega,S})$ is canonically isomorphic to 
$(h_\ast\varphi_{f,S'},{\rm Tr}_{S'/S}\circ B_{f,\omega,S'})$. 
\end{lm}
\proof{
It follows from Lemma \ref{Trfet}. 
\qed}

\begin{lm}\label{prodphi}
Let $X_j\ (j=1,2)$ be smooth $S$-schemes and let $f_j\colon X_j\to \A^1_S$ be $S$-morphisms whose singular loci $Z_j$ are 
finite over $S$. Let $X:=X_1\times_SX_2$ and let $f\colon X\to\A^1_S$ be the map $(x,y)\mapsto f_1(x)+f_2(y)$. 
Then the singular locus of $f$ is $Z_1\times_SZ_2$. We have a canonical isomorphism 
\begin{equation*}
(\varphi_{f},B_{f,\omega})\cong (\varphi_{f_1},B_{f_1,\omega})\otimes_{\OO_S}
(\varphi_{f_2},B_{f_2,\omega})
\end{equation*}
for an everywhere non-zero differential $\omega$. 
\end{lm}
\proof{
This follows from Lemma \ref{prodtr}. Indeed, by the projections ${\rm pr}_j\colon X\to X_j$, we have canonically 
$\Omega^1_{X/S}\cong{\rm pr}_1^\ast\Omega^1_{X_1/S}\oplus{\rm pr}_2^\ast\Omega^1_{X_2/S}$. Under this 
identification, the pullback 
$f^\ast\colon \Omega^1_{\A^1_S/S}\to\Omega^1_{X/S}$ equals to $f_1^\ast+f_2^\ast$. The assertion follows. 
\qed}

Next we recall the notion of Milnor number and its continuity. 

Let $S$ be a scheme and $X$ be a smooth $S$-scheme. consider an $S$-morphism  
\begin{equation*}
f\colon X\to C
\end{equation*}
to a smooth $S$-curve. Let $Z$ be an open and closed subscheme of the 
singular locus of $f$ which is finite over $S$. 
Recall then that $Z$ is finite locally free over $S$ (Lemma \ref{cg!now}.1). 
\begin{df}\label{Mildef}
 Let the notation be as above. 
 We call the rank of $\OO_{Z}$ as locally free $\OO_S$-module {\rm the Milnor number} of $f$ along $Z$ and write $\mu(f,Z)$ for it. 
 When $Z$ stands for the underlying set of such an open and closed subscheme of the singular locus, we also write 
 $\mu(f,Z)$ for it, as this convention does not lead to confusions. 
\end{df}
When $S$ is the spectrum of a field and $Z$ consists of only one point $x$, 
we also write $\mu(f,x)$ for 
$\mu(f,Z)$. 

\begin{pr}\label{Milconti}
Consider a family of isolated singularities (Definition \ref{isoldef}.2) 
\begin{equation*}
\xymatrix{
Z\ar@{^{(}->}[r]^i\ar[rrd]_g&X\ar[rr]^f\ar[rd]&&C\ar[ld]\\
&&S.&
}
\end{equation*}
Let $\bar{z}$ be a geometric point of $Z$ and let 
$\bar{s}$ be the geometric point of $S$ induced from $\bar{z}$ by composing $g$. 
Let 
$\bar{t}$ be a geometric point of $S$ which specializes to $\bar{s}$.
Then we have 
\begin{equation}\label{valmil}
\mu(f_{\bar{s}},\bar{z})=\sum_{z'\in Z_{(\bar{z})}\times_{S_{(\bar{s})}\bar{t}}}\mu(f_{\bar{t}},z'). 
\end{equation}
\end{pr}
For the notation $S_{(\bar{s})},f_{\bar{s}}$, etc., see the end of the introduction. 
\proof{
We may replace $S$ by the strict henselization $S_{(\bar{s})}$. Then we replace $X$ by an open neighborhood around (the image of) 
$\bar{z}$ so that $Z$ equals to $Z_{(\bar{z})}$. In this case, $\OO_Z$ is a finite free 
$\OO_S$-module and the left (resp. right) hand side in (\ref{valmil}) is the rank of $\OO_Z\otimes_{\OO_S}k(\bar{s})$ 
(resp. $\OO_Z\otimes_{\OO_S}k(\bar{t})$). 
\qed}

\subsection{Ordinary quadratic singularities}
In this subsection, we recall the notion of {\it ordinary quadratic singularities} and 
compute various invariants of such singularities, which is necessary for the proofs 
of the main theorems. 

\begin{df}(\cite[1.2]{PL}) \label{quadsing}
\begin{enumerate}
\item Let $Y$ be a scheme of finite type over a field $k$. Let $y\in Y$ be a closed point. 
\begin{enumerate}
\item When $k$ is algebraically closed, we say that $y$ is {\rm an ordinary (resp.~non-degenerate) quadratic point} if the completion $\hat{\OO}_{Y,y}$ is $k$-isomorphic 
to $k[[x_0,\dots,x_{n}]]/(f)$ where $f$ starts in degree $2$ and the homogeneous part of degree $2$ is an ordinary (resp.~non-degenerate)  quadratic form 
\cite[1.1]{Quad}. 
\item In general, we say that $y$ is {\rm an ordinary (resp.~non-degenerate) quadratic point} if the following holds: 
Let $Y\times_{k}\bar{k}$ be the base change to $\bar{k}$, where $\bar{k}$ is an algebraic closure of $k$. For some (hence all) $\bar{y}\in Y\times_{k}\bar{k}$ which projects to $y$, $\bar{y}$ is an ordinary (resp.~non-degenerate) quadratic point in the sense of $(a)$. 
\end{enumerate}
\item Let $f\colon X\to T$ be a flat morphism of schemes of finite type. We say that a point $x\in X$ is  {\rm an ordinary 
(resp.~non-degenerate) quadratic point} if so is 
$x$ in the fiber $X\times_Tf(x)$. 
\end{enumerate} 
\end{df}
Let $y\in Y$ be an ordinary quadratic point on a $k$-scheme $Y$ of finite type. Then $y$ is non-degenerate if and only if 
${\rm char}(k)$ is odd or $Y$ has odd dimension at $y$ \cite[1.2.2]{PL}. 

First we recall results on deformations of ordinary quadratic singularities (\cite{PL}). 
Let $(R,\mathfrak{m}_R)$ be a henselian local ring with the residue field $k$ and $S={\rm Spec}(R)$ be its affine spectrum. 
We write $s$ for the closed point. 
Let $X$ be a smooth $S$-scheme of relative dimension $n+1$ and let 
$f\colon X\to C$ be an $S$-morphism to a smooth $S$-curve $C$ with the singular locus $Z$ finite over $S$. We 
assume that $Z$ is {\it local}. Thus $Z$ has only one closed point,  which we denote by $x$. 

To state the lemma, let us introduce the following notation. 
Let $(A,\mathfrak{m}_A)$ be a henselian local ring with a local homomorphism $R\to A$. An element $a\in \mathfrak{m}_A$ 
gives a unique morphism of local $R$-algebras $R\{t\}\to A$, $t\mapsto a$, where $R\{t\}$ is the henselization of $R[t]$ at 
$(\mathfrak{m}_R,t)$. 
We write $f(a)$ for the image of $f\in R\{t\}$ under this map. 

\begin{lm}(\cite[Proposition 1.3.1]{PL})\label{classquad}
Let the notation be as above. We assume that $k(x)$ is isomorphic to $k$. 
\begin{enumerate}
\item Suppose that $k$ is of odd characteristic or that $n+1$ is even.  
Then there is an $R$-isomorphism 
$R\{t\}\cong \OO_{C,(f(x))}$ of henselizations such that 
the henselization of $\OO_{X,x}$ is isomorphic to 
\begin{equation*}
R\{t,t_0,\dots,t_n\}/(Q-t)
\end{equation*}
as algebras over $\OO_{C,(f(x))}=R\{t\}$. Here $Q\in R[t_0,\dots,t_n]$ is a non-degenerate quadratic form. 
The Milnor number $\mu(f_s,x)$ of the closed fiber $f_s\colon X_s\to C_s$ equals to $1$. 
The singular locus $Z$ equals to the closed subscheme of 
$X_{(x)}$ defined by $t_0=\cdots=t_n=0$. In particular, $Z\to S$ is 
isomorphic. 

\item Suppose that $k$ is of characteristic $2$ and that $n+1$ is odd. Then, if we choose appropriately an 
 $R$-isomorphism 
$R\{t\}\cong \OO_{C,(f(x))}$, an element $b(t)\in R\{t\}$ in the maximal ideal whose reduction $\bar{b}(t)$ in $k\{t\}$ is non-zero, 
and a non-degenerate quadratic form $Q\in R[t_1,\dots,t_n]$ in variables $t_1,\dots,t_n$, 
the henselization $\OO_{X,(x)}$ is $R\{t\}$-isomorphic to 
\begin{equation}\label{presinit}
\OO_{X,(x)}\cong R\{t,t_0,t_1,\dots,t_n\}/(Q+t_0^2+b(t)t_0-t). 
\end{equation}
The Milnor number $\mu(f_s,x)$ of the special fiber equals to $2{\rm ord}_t(\bar{b}(t))$, twice the normalized valuation 
of the reduction $\bar{b}(t)$. 

Further if $R$ is over $\F_2$ and if $\bar{b}(t)\in k\{t\}$ is a 
uniformizer, the map $Z\to S$ is finite flat of finite presentation of degree $2$ and is 
a universal homeomorphism. 
\item In the situation $2$, the $R\{t\}$-algebra (\ref{presinit}) is isomorphic to the henselization of 
\begin{equation}\label{userep}
R\{t_1,\dots,t_n\}\otimes_RR\{u\}
\end{equation}
at $(\mathfrak{m}_R,u,t_1,\dots,t_n)$ which is regarded as an $R\{t\}$-algebra via the map 
\begin{equation*}
t\mapsto Q+\tilde{t}:=Q\otimes1+1\otimes\tilde{t}. 
\end{equation*}
Here $\tilde{t}$ is an element of $R\{u\}$ which satisfies the relation 
$u^2+b(\tilde{t})u-\tilde{t}=0$. 
\end{enumerate}
\end{lm}
\proof{
For $1$, 
the assertion on the presentation of $\OO_{X,(x)}$ is given in \cite[1.3.1]{PL}. Note that $b$ in \cite[1.3.1(i)]{PL} 
can be taken as $b=t$ as in $1$ because the fiber $X_s$ is regular. 
For $2$, the presentation is given in \cite[1.3.1(ii)]{PL}. 
We can take $b,c$ in {\it loc.~cit.}~as in 
$2$ because $X_s$ is regular and $x$ is an isolated singular point. 
The computation on the Milnor numbers are given in 
\cite[1.13]{Mil} or can be verified directly from the definition in this case. 

We verify the assertions on $Z$ in each case $1$ and $2$. For $1$, the assertion follows since 
$Q$ is non-degenerate, hence the coefficients of $dQ$ span the ideal 
$(t_0,\dots,t_n)$. For $2$, assume that $R$ is an $\F_2$-algebra. Then we have 
\begin{equation*}
d(Q+t_0^2+b(t)t_0-t)=dQ+b(t)dt_0-(1-b'(t)t_0)dt
\end{equation*}
where $b'(t)=\frac{db}{dt}$. Hence we have 
\begin{equation*}
\OO_Z=\OO_{X,(x)}/(t_1,\dots,t_n,b(t))\cong 
R\{t,t_0\}/(t_0^2+b(t)t_0-t,b(t))\cong R\{t_0\}/(b(t_0^2)). 
\end{equation*}
As $R$ is of characteristic $2$, we have $b(t_0^2)=b(t_0)^2$. By the assumption that $\bar{b}(t)$ is a uniformizer of 
$k\{t\}$, the map $R\{b\}\to R\{t_0\}$ sending 
$b$ to $b(t_0)$ is an isomorphism. The assertion follows.

We  show $3$. 
Let $A$ be the henselization of (\ref{userep}) at $(\mathfrak{m}_R,t_1,\dots,t_n,u)$. 
We show that the polynomial $Q(x)+x_0^2+b(Q(t)+\tilde{t})x_0-(Q(t)+\tilde{t})\in A[x_0,\dots,x_n]$ has a solution 
$(u',t_1',\dots,t_n')\in A^{n+1}$ which equals to $(u,t_1,\dots,t_n)$ mod $(u,t_1,\dots,t_n)^2+\mathfrak{m}_RA$. Such a solution 
gives an \'etale map $A\to R\{t,x_0,\dots,x_n\}/(Q(x)+x_0^2+b(t)x_0-t)$ of henselian $R\{t\}$-algebras, 
which is an isomorphism. 

We find such a solution applying \cite[5.11]{App} to $y^\circ=(u,t_1,\dots,t_n),\mathfrak{a}=Au$, and 
$f=Q(x)+x_0^2+b(Q(t)+\tilde{t})x_0-(Q(t)+\tilde{t})$ 
with the notation given there. 
We have 
\begin{align*}
Q(t)+u^2+b(Q(t)+\tilde{t})u-(Q(t)+\tilde{t})&=u^2+b(Q(t)+\tilde{t})u-\tilde{t}\\
&\equiv u^2+b(\tilde{t})u-\tilde{t}=0\hspace{5mm}{\rm mod}\hspace{1mm}(t_1,\dots,t_n)^2u. 
\end{align*}
As $Q$ is non-degenerate, the ideal $\Delta$ in {\it loc.~cit.}~equals to 
$(t_1,\dots,t_n,2u+b(Q(t)+\tilde{t}))$. 
Thus, by \cite[5.11]{App}, we find a solution $(u',t_1',\dots,t_n')$ as desired. 
\qed}

We compute the bilinear form $(\varphi_{f},B_{f,\omega})$ for certain quadratic 
singularities in mixed characteristic. 
\begin{pr}\label{compphi}
\begin{enumerate}
\item Let $R$ be a ring and  
$Q\in R[t_0,\dots,t_n]$ be a non-degenerate quadratic form. Let 
\begin{equation*}
f\colon X:={\rm Spec}(R[t_0,\dots,t_n])\to\A^1_{R}={\rm Spec}(R[t])
\end{equation*}
be the $R$-morphism defined by $t\mapsto Q$. Then the singular locus $Z$ of $f$ equals to 
$\{t_0,\dots,t_n=0\}\subset X$ as a subscheme. The bilinear form $(\varphi_{f},B_{f,\omega})$ 
for a differential $\omega$ on $\A^1_{R}$ which is non-zero along $\{t=0\}$ is isomorphic to 
$(R,B)$ where $B\colon R\otimes R\to R$ is defined by $1\otimes 1\to \frac{1}{{\rm disc} Q}\cdot(\frac{\omega}{dt}|_{t=0})^{n+1}$. 
\item Let $k$ be a perfect field of characteristic $2$. Let $R=W_{m+1}(k)=W(k)/p^{m+1}$ be the ring of Witt vectors of length $m+1$. 
Let $A$ be an \'etale $R[t]$-algebra such that $A/(t)$ is isomorphic to $R$. 
Let $Q\in R[t_1,\dots,t_n]$ be a 
non-degenerate quadratic form (hence $n$ is assumed even) 
and $b\in A$ be an element such that its image in $A/p$ is non-zero and is contained in the ideal $t\cdot A/p$. Consider 
\begin{equation*}
f\colon X':={\rm Spec}(A[t_0,t_1,\dots,t_n]/(Q+t_0^2+bt_0-t))
\to\A^1_R={\rm Spec}(R[t]). 
\end{equation*}
Let $x\in X'$ be the $k$-rational point defined by $\{p,t_0,\dots,t_n=0\}$. 
Then $X'$ is $R$-smooth around $x$, and $x$ is an isolated singular point with respect to $f$. 

Let $X$ be an open neighborhood of $x$ 
in $X'$ so that $X$ is $R$-smooth and the singular locus $Z$ of $f$ 
 consists of only $x$. 
Further assume that the Milnor number $\mu(f\otimes_Rk,x)$ equals to $2$. Then the discriminant of the bilinear form $(\varphi_{f},B_{f,dt})$ equals to 
$-1$ modulo $(W_{m+1}(k)^\times)^2$. 
\end{enumerate}

\end{pr}

In the situation $2$, we will prove the equality ${\rm disc}B_{f,dt}\equiv (-1)^{\frac{(n+1)\mu(f,x)}{2}}$ without the assumption on the Milnor number  
as a special case of Proposition \ref{arfodd}. 
\proof{
Put $S={\rm Spec}(R)$. We apply Lemma \ref{computePsi} to 
$({\cal E},s)=(\Omega^1_{X/S},f^\ast \omega)$ in 
each case. We follow 
the notations given there. In particular, we write $\xi_i$ for $1\otimes t_i-t_i\otimes1\in
{\cal O}_{Z\times_SX}$. 
 
 $1$. 
 In this case, $Z$ equals to the $0$-section $S\hookrightarrow X=\A^{n+1}_S$ and $\xi_i=1\otimes t_i$ in $\OO_{Z\times_SX}$. 
 We trivialize $\Omega^1_{X/S}$ by $dt_0\dots,dt_n$ and identify the complex (\ref{KosE}) defined 
 for $s=f^\ast dt$ with the Koszul complex $K^\bullet(\underline{f})$ defined by $\frac{\partial Q}{\partial t_0},\dots,\frac{\partial Q}{\partial t_n}$. 
 Let 
 $F\colon \OO_X^{n+1}\to\OO_X^{n+1}$ be the $\OO_X$-linear map defined by the hessian 
 $(\frac{\partial^2Q}{\partial t_i\partial t_j})_{i,j}$. Its exterior powers give an isomorphism 
 $K(F)^\bullet\colon K^\bullet(\underline{\xi})\to K^\bullet(\underline{f})$. The map $\frac{1}{{\rm disc}Q}K^\bullet(F)$ then gives a morphism of resolutions of ${\cal O}_Z$. 
Using this, we can identify the residue symbol ${\rm Res}[\frac{}{f^\ast dt}]$ with 
 $\OO_S\xrightarrow{\frac{1}{{\rm disc}Q}}\OO_S$. The assertion for $\omega$ comes from that for $dt$ and Lemma \ref{bpofbil}.4. 
 
 $2$.  
 First we treat the case $n=0$. Hence 
 we assume that $X'={\rm Spec}(A[u]/(u^2+bu-t))$. 
 We compute 
 \begin{equation*}
 d(u^2+bu-t)=(2u+b)du-(1-ub')dt
 \end{equation*}
 where $b'=\frac{db}{dt}$. As $u,t$ are local parameters on ${\rm Spec}(A[u])$, $X'$ is smooth at $x$. Around $x\in X'$, we have 
 \begin{equation}\label{diffq}
 f^\ast dt=\frac{2u+b}{1-ub'}du. 
 \end{equation}
 Hence, in the special fiber, we have $f^\ast dt=\frac{b}{1-ub'}du$.  
Thus there exists an open neighborhood $X$ around $x$ such that 
$f|_X$ has only one singular point at $x$. 
 
In the special fiber, $X_s\to\A^1_k$ is totally wildly ramified over $\{t=0\}$ and the kernel of the map 
$\OO_{X_s,x}\to\OO_{Z_s,x}$ is generated by $\bar{b}$ by (\ref{diffq}). Thus 
$\OO_Z\otimes_Rk$ has a basis $1,u,\dots,u^{2{\rm ord}_{t}(\bar{b})-1}$ as a $k$-vector space. 
By Nakayama's lemma, $\OO_Z$ itself has a basis $1,u,\dots,u^{2{\rm ord}_{t}(\bar{b})-1}$ as a free $R$-module. 

Assume further that $\mu(f_s,x)=2$. We identify $\OO_Z$ with the free $R$-module $R\oplus Ru$. 
In order to compute $(\varphi_{f},B_{f,dt})$, 
we express several elements in $\OO_Z$ as linear combinations of 
$1$ and $u$, which are summarized from $1$ to $5$ below. 
As 
$\bar{b}$ is a uniformizer of $k\{t\}$, 
the map $W_2(k)[b]/(b^2)\to R\{t\}/(4,b^2)$ is an isomorphism. 
Thus 
we find a relation of the form 
\begin{equation}\label{t,b}
t\equiv2c+ab\hspace{3mm}{\rm mod.}\hspace{1mm}(4,b^2)
\end{equation}
 where $c\in R$ and $a\in R^\times$. 
 Up to replacing $\A^1_R$ by an \'etale neighborhood $C$ 
 around $0$, 
 we assume that we have a relation as (\ref{t,b}) in $\OO_C$.

We express $u^2\in\OO_Z$ as the form  
 $u^2=m+nu$ where $m,n\in R$. Then the ratio $\frac{u^2-m-nu}{2u+b}$ is defined around $x$ in $X$. 
We compute various elements in $\OO_Z/4\OO_Z$ as follows. Here $\equiv$ indicates the 
congruence modulo $4$. 
\begin{enumerate}
\item[(1)] $b= -2u$ in ${\cal O}_Z$. 
\item[(2)] $t\equiv 2c-2au$ in ${\cal O}_Z/4{\cal O}_Z$. 
\item[(3)] $b'\equiv \frac{1}{a}$ in ${\cal O}_Z/4{\cal O}_Z$. 
\item[(4)] $m\equiv 2c$, $n\equiv2a$ in $R/4R$. 
\item[(5)] $\frac{u^2-m-nu}{2u+b}\equiv -a+u$ in ${\cal O}_Z/4{\cal O}_Z$. 
\end{enumerate}
The first one follows since $2u+b$ is a generator of the ideal sheaf defining $\OO_Z$. 
The second one follows from $(1)$ and (\ref{t,b}). 
Differentiating (\ref{t,b}) by $t$, we have 
\begin{equation*}
1\equiv ab'\hspace{3mm}{\rm mod.}\hspace{1mm}(4,2bb',b^2).  
\end{equation*}
As we have $b=-2u$ in $\OO_Z$, $(3)$ holds in $\OO_Z/4\OO_Z$. 

We have 
\begin{align*}
(u-a)(2u+b)-u^2&=u^2+bu-2au-ab\\&=t-2au-ab
\equiv 2c-2au\hspace{3mm}{\rm mod.}\hspace{1mm}(4,b^2)
\end{align*}
by (\ref{t,b}). Since we have $(4,b^2)=(4,(2u+b)^2)$, the assertions 
$(4)$ and $(5)$ follow.

To compute $(\varphi_{f},B_{f,dt})$, 
we construct a morphism of complexes as (\ref{morKos}). 
Define $\zeta\colon \OO_Z\otimes_{\OO_S}\OO_X\to\OO_X$ to be the $\OO_X$-linear 
map sending $1\otimes 1\mapsto 1,u\otimes1\mapsto u$. 
We need to find an $\OO_X$-linear map $\eta\colon \OO_Z\otimes_{\OO_S}\OO_X\to 
\OO_X$  which makes the diagram
\begin{equation*}
\xymatrix@M=8pt{
\OO_Z\otimes_{\OO_S}\OO_X\ar[r]^{1\otimes u-u\otimes1}\ar[d]_{\eta}&\OO_Z\otimes_{\OO_S}\OO_X\ar[d]_{\zeta}\\
\OO_X\ar[r]^{\frac{2u+b}{1-ub'}}&\OO_X
}
\end{equation*}
commutative. The image of $1\otimes1$ by $\nu\colon\OO_Z\otimes_{\OO_S}\OO_X\xrightarrow{1\otimes u-u\otimes1}
\OO_Z\otimes_{\OO_S}\OO_X\xrightarrow{\zeta}\OO_X$ is zero. Thus we can set $\eta(1\otimes 1):=0$. 
The image $\nu(u\otimes1)$ equals to $\zeta(u\otimes u-(m+nu)\otimes1)
=u^2-m-nu$. As this can be divided by $2u+b$, we can define 
$\eta(u\otimes1)$ to be $\frac{u^2-m-nu}{2u+b}(1-ub')$.  Set $\beta:=
\frac{u^2-m-nu}{2u+b}(1-ub')$. 
Under the identification 
\begin{equation*}
\mathcal{H}om_{\OO_X}(\OO_{Z\times X},\OO_X)\cong 
\mathcal{H}om_{\OO_S}(\OO_{Z},\OO_S)\otimes_{\OO_S}\OO_X, 
\end{equation*}
the map $\alpha\eta$ with $\alpha\in \OO_X$ corresponds to 
\begin{equation*}
u^\ast\otimes \beta\alpha
\end{equation*}
where $u^\ast\colon\OO_Z\to\OO_S$ maps $
1\otimes1\mapsto0, u\otimes1\mapsto 1$. Therefore, 
we can identify the bilinear form $B_{f,dt}$ with 
\begin{equation*}
B\colon \OO_Z\times\OO_Z\to\OO_S,\hspace{1mm}(r_1,r_2)\mapsto 
u^\ast(\beta r_1r_2). 
\end{equation*}
From now on, we compute everything modulo $4$. 
By $3,4,5$ above, we have 
\begin{equation*}
\beta\equiv (-a+u)(1-\frac{u}{a})=2u-a-\frac{u^2}{a}\equiv 
2u-a-\frac{2c}{a}-2u=-a-\frac{2c}{a}. 
\end{equation*}
Thus we have 
\begin{equation*}
B(r_1,r_2)\equiv -(a+\frac{2c}{a})u^\ast(r_1r_2). 
\end{equation*}
Therefore, we have 
 $B(1,1)\equiv 0,
B(u,u)\in2{\cal O}_S$. 
Hence we have 
\begin{equation*}
{\rm disc}B_{f,dt}=B(1,1)B(u,u)-B(1,u)^2
\equiv -B(1,u)^2\hspace{3mm}{\rm mod}\hspace{1mm}8. 
\end{equation*}
This completes the proof, as $B$ is non-degenerate and $1+8W_m(k)$ is contained in $(W_m(k)^\times)^2$ by Lemma \ref{Wunr}. 

For general $n\geq0$, we reduce it to 
the case $n=0$. By Lemma \ref{classquad}.3, the henselization 
 $\OO_{X',(x)}$ is $R\{t\}$-isomorphic to $R\{t_1,\dots,t_n\}\otimes_R
 R\{u\}$ by the map sending $t\mapsto Q\otimes1+1\otimes\tilde{t}$ where $\tilde{t}\in R\{u\}$ satisfies the relation 
 $u^2+b(\tilde{t})u-\tilde{t}=0$. As the Milnor number of 
 ${\rm Spec}(A[u]/(u^2+bu-t))\to\A^1_R$ at $(t,u)=(0,0)$ is $2$,  Lemma \ref{prodphi} and $1$ imply the assertion. 
\qed}
\begin{cor}\label{Arfquad}
In the situation in $1$ of Proposition \ref{compphi}, assume that 
$R=W_3(k)=W(k)/8W(k)$ for a perfect field $k$ of characteristic $2$. 
Then $(-1)^{\frac{n+1}{2}}$ times  the discriminant of $B_{f,dt}$  is contained in the image of 
$1+4W_3(k)\to W_3(k)^\times/(W_3(k)^\times)^2$ and equals to $1+4[{\rm Arf}(Q\otimes_Rk)]$ where 
${\rm Arf}(-)$ denotes the Arf invariant of a non-degenerate quadratic form. 
\end{cor}
\proof{

The discriminant of $B_{f,dt}$ equals to 
$\frac{1}{{\rm disc}(Q)}\equiv{\rm disc}(Q)$ mod $(W_3(k)^\times)^2$. Then the assertion follows from \cite[1.4]{FQE}. 
Indeed, let $Q'$ be a quadratic form on a free $W(k)$-module 
which is a lift of $Q$. It is proved in loc.~cit.~that 
the discriminant of $Q'$ equals to $(-1)^{\frac{n+1}{2}}(1+4[{\rm Arf}(
Q'\otimes_{W(k)}k)])$ in $W(k)^\times/(W(k)^\times)^2$. 
Since we have $W(k)^\times/(W(k)^\times)^2\cong 
W_3(k)^\times/(W_3(k)^\times)^2$ by Lemma \ref{Wunr}.1, the assertion is verified. 
\qed}

We give a corresponding result of Proposition \ref{compphi} 
for local epsilon factor. The following is used in the proof of Theorem 
\ref{mainMil}. 

Let $p$ be a prime number. Fix a prime number 
$\ell\neq p$ and take an algebraic closure $\overline{\Ql}$ of the 
$\ell$-adic field $\Ql$. 
We fix a non-trivial 
additive character $\psi\colon\F_p\to\overline{\Ql}^\times$. 
For a finite extension $k/\F_p$, we write $\psi_k$ for the composition 
$k\xrightarrow{{\rm Tr}_{k/\F_p}}\F_p\xrightarrow{\psi}\overline{\Ql}^\times$. We use these characters to define local epsilon factors 
in equal-characteristic case, as explained in \cite[(3.1.5.4)]{Lau}. 
We follow the notation for local epsilon factor in loc. cit. 

Let $f\colon X\to\mathbb{A}^1_k$ be a $k$-morphism of smooth 
$k$-schemes with a $k$-rational isolated singular point $x\in X$. 
We write $R\Phi_f(\Ql)_x$ for the vanishing cycles complex 
supported at $x$. This is a bounded complex of finite dimensional 
$\Ql$-vector spaces with $G_\eta$-action, where $G_\eta$ is the 
absolute Galois group of generic point $\eta$ of the henselization $\A^1_{k,(f(x))}$. Hence the local epsilon factor 
$\varepsilon_0(\A^1_{k,(f(x))},R\Phi_f(\Ql)_x,dt)$ is defined. 

\begin{lm}\label{epex}
Let $k$ be a finite extension of $\F_p$ with $q$ elements. 
\begin{enumerate}
\item Assume that $p$ is odd. For $a\in k^\times$, define 
\begin{equation*}
f\colon\mathbb{A}^1_k\to\mathbb{A}^1_k
\end{equation*}
by $t\mapsto at^2$. 
Then we have 
\begin{equation*}
\varepsilon_0(\A^1_{k,(0)},R\Phi_{f}(\Ql)_0,dt)=(\frac{-a}{k})
\sum_{x\in k}\psi_k(x^2). 
\end{equation*}

\item Assume that $p=2$. Let $a\in k$ and consider a map 
\begin{equation*}
f\colon \A^2_k\to\A^1_k
\end{equation*}
defined by $(x,y)\mapsto x^2+xy+ay^2$. Then the product 
\begin{equation*}
-\varepsilon_0(\A^1_{k,(0)},R\Phi_{f}(\Ql)_0,dt)q
\end{equation*}
is $\pm1$. This is $1$ if and only if $a\in\wp(k)$. 
\item Assume that $p=2$. Let $t\colon C\to\A^1_k$ 
be an \'etale morphism of smooth $k$-curves such that $C\times_{\A^1_k}0$ equals to a $k$-rational point, which we also denote by $0$. 
 Let $b\in\Gamma(C,\OO_C)$ be an element which gives a uniformizer at $0$. Let 
\begin{equation*}
f\colon X:={\rm Spec}(\OO_C[u]/(u^2+bu-t))\to C. 
\end{equation*}
Then $X$ is smooth around $\{u=0\}$ and 
we have 
\begin{equation*}
\varepsilon_0(C_{(0)},R\Phi_{f}(\Ql)_0,dt)=q. 
\end{equation*}
\end{enumerate}
\end{lm}
\proof{
Although they are well-known, we include a proof for completeness. 
Set $F:=k((t))$ and let $\OO_F$ be the ring of integers in $F$. 

$1$. The quadratic extension $E:=k((t))[u]/(u^2-at)$ of $F$ corresponds to a tame character $\chi$ 
of order $2$ of the absolute Galois group of $F$. 
Then the vanishing cycles complex $R\Phi_f(\Ql)_0$ is isomorphic to $\chi$. We have (\cite[5.10]{Del}) 
\begin{align*}
\varepsilon_0(\A^1_{k,(0)},\chi,dt)&=\chi(t)\sum_{x\in k^\ast}\chi(x)\psi_k({\rm Res}(x\frac{dt}{t}))\\
&=\chi(-a)\sum_{x\in k^\ast}\chi(x)\psi_{k}(x). 
\end{align*}
As $E/F$ is a totally tamely ramified quadratic extension, $\chi|_{k^\ast}$ is the unique non-trivial character of order $2$,~i.e.~it equals to the Legendre symbol $(\frac{}{k})$. The assertion follows. 

$2$. Let $X:={\rm Proj}({\cal O}_F[x,y,z]/(x^2+xy+ay^2-tz^2))$. This is regarded as a closed subscheme of 
$\Proj^2_{\OO_F}$. By the jacobian criterion, the structure map 
$\bar{f}\colon X\to{\rm Spec}(\OO_F)$ is smooth except at the point $(x:y:z)=(0:0:1)$ in the special fiber. By 
the local acyclicity of smooth morphism and the proper base change theorem, we have a distinguished triangle 
\begin{equation*}
R\Gamma(X_{\bar{s}},\Ql)\to R\Gamma(X_{\bar{\eta}},\Ql)\to R\Phi_{\bar{f}}(\Ql)_0\to 
\end{equation*}
of $G_F$-representations. As $X_{\bar{\eta}}$ is isomorphic to a projective line, the cohomology groups are 
$H^0(X_{\bar{\eta}},\Ql)=\Ql,H^2(X_{\bar{\eta}},\Ql)\cong\Ql(-1)$ and $0$ in the other degrees. 
The geometric special fiber $X_{\bar{s}}$ is the union of two projective lines glued at the origins. 
The action of $G_k$ on the irreducible components is non-trivial if and only if the special fiber 
$X_s$ is irreducible, which happens if and only if $a$ is non-zero in $k/\wp(k)$. 
Let $\chi_a$ be the correspnding character of $G_k$ via 
Artin--Schreier theory. 
We have $R\Phi_{\bar{f}}(\Ql)_0\cong \chi_a\otimes\Ql(-1)[-1]$. Since we have 
$R\Phi_{\bar{f}}(\Ql)_0\cong R\Phi_{f}(\Ql)_0$, the assertion follows in this case. 

$3$. We are treating the quadratic extension $E:=F[u]/(u^2+bu-t)$. Let $\chi$ be the character of order $2$ corresponding to it. 
The Artin conductor of $\chi$ is $2$. Take $\gamma\in F$ with ${\rm ord}_F(\gamma)=2$. Let 
$dx$ be the Haar measure on the additive topological group $F$ normalized as $\OO_F$ has volume $1$. 
We have (cf. \cite[5.8]{Del})
\begin{align*}
\varepsilon_0(\OO_F,\chi,dt)&=\int_{\gamma^{-1}\OO_F^\times}\chi(x)\psi_k({\rm Res}(xdt))dx\\
&=q^2\int_{\OO_F^\times}\chi(\gamma x)\psi_k({\rm Res}(\frac{xdt}{\gamma}))dx\\
&=\chi(\gamma)q^2\sum_{a\in k^\times}\sum_{b\in(1+(t))/(1+(t)^2)}
\int_{1+(t)^2}\chi(abx)\psi_k({\rm Res}(\frac{abxdt}{\gamma}))dx\\
&=\chi(\gamma)q^2\sum_{a,b}\chi(ab)\int_{(t)^2}\psi_k({\rm Res}(\frac{abdt}{\gamma}))\psi_k({\rm Res}(\frac{abxdt}{\gamma}))dx. \\
\end{align*}
Note that the term $\psi_k({\rm Res}(\frac{abxdt}{\gamma}))=1$ since $x\in (t)^2$. We proceed 
\begin{align*}
\chi(\gamma)q^2\sum_{a,b}\chi(ab)\int_{(t)^2}\psi_k({\rm Res}(\frac{abdt}{\gamma}))\psi_k({\rm Res}(\frac{abxdt}{\gamma}))dx&=
\chi(\gamma)\sum_{a,b}\chi(ab)\psi_k({\rm Res}(\frac{abdt}{\gamma})). 
\end{align*}
Take $\alpha\in k^\times$ so that $\chi(1+z)=\psi_k({\rm Res}(\alpha\frac{zdt}{\gamma}))$ for $z\in(t)/(t)^2$. 
We have 
\begin{align*}
\chi(\gamma)\sum_{a,b}\chi(ab)\psi_k({\rm Res}(\frac{abdt}{\gamma}))&=
\chi(\gamma)\sum_{a\in k^\times,z\in(t)/(t)^2}\chi(a(1+z))\psi_k({\rm Res}(\frac{a(1+z)dt}{\gamma}))\\
&=\chi(\gamma)\sum_{a\in k^\times}\chi(a)\psi_k({\rm Res}(\frac{adt}{\gamma}))
\sum_{z\in(t)/(t)^2}\psi_k({\rm Res}(\frac{(\alpha+a) zdt}{\gamma})). 
\end{align*}
The last sum $\sum_{z\in(t)/(t)^2}\psi_k({\rm Res}(\frac{(\alpha+a) zdt}{\gamma}))$ is $0$ unless $a=\alpha$. 
Thus we obtain 
\begin{equation*}
\varepsilon_0(\OO_F,\chi,dt)=\chi(\gamma)\chi(\alpha)\psi_k({\rm Res}(\frac{\alpha dt}{\gamma}))q. 
\end{equation*}
From now on, we take $\gamma$ to be $t^2$. For this choice, the terms $\chi(\gamma),\psi_k({\rm Res}(\frac{\alpha dt}{\gamma}))$ are trivial. 
The term $\chi(\alpha)$ is also trivial since the orders of elements in $k^\times$ are odd. Thus the assertion follows. 
\qed}

\subsection{Deforming isolated singularities}

Here we explain that any function with an isolated singularity of equal characteristic 
can be embedded into a family of isolated singularities which generically contains ordinary quadratic singularities (cf. \cite[2.5]{Mil}). 
\begin{df}\label{deffuniso}
Let $T$ be a henselian trait. 
{\rm A function with an isolated singularity} is a map of schemes $f\colon X\to T$ which can be obtained as the henselization 
of a $T$-scheme $X'$ flat of finite type at a closed point $x\in X'$ with the following properties. 
$X'$ is regular, $x$ maps to the closed point of $T$, and $X'\setminus \{x\}$ is smooth over $T$. 
In particular, $X$ is regular and henselian local. 
\end{df}
\begin{lm}\label{deform}
Let $k$ be a field. Let $f\colon X\to{\rm Spec}(\OO)$ be a function with an isolated singularity (Definition \ref{deffuniso}) such that $\OO$ is isomorphic to 
the henselization $\OO_{\A^1_{k,(0)}}$. Let $x\in X$ be the singular point. Assume that the residue field $k(x)$ is isomorphic to $k$. 
Then there exists a family of isolated singularities (Definition \ref{isoldef}) 
\begin{equation}\label{familydef}
 \xymatrix{
Y\ar[rr]^{\tilde{f}}\ar[rd]_\pi&&\A^1_S\ar[ld]\\
 &S
 }
 \end{equation}
 such that 
\begin{enumerate}
\item $S$ is a smooth connected $k$-curve. The   
singular locus $Z$ is {\rm finite} over $S$. 
\item There exists a $k$-rational point $s\in S(k)$ such that 
\begin{enumerate}
\item $Z\times_Ss$ consists of only one point $x_0$. The point $x_0$ maps to the origin of $\A^1_{k(s)}=\A^1_S\times_Ss$. 
\item There is a $k$-homomorphism $\OO_{\A^1_{k(s),(0)}}\to\OO$ which 
preserves uniformizers such that the henselization of $Y\times_{\A^1_S}\OO$ 
at $x_0$ is $\OO$-isomorphic to $X$. 
\end{enumerate}
\item There exists an open dense subset $U\subset S$ 
which has a section $\iota\colon U\to Z$ such that, for any $u\in U$, 
$\iota(u)$ is an ordinary quadratic point of $\tilde{f}$ 
with Milnor number $1$ or $2$. 
\item $Y,S$ are affine. 
\end{enumerate}
\end{lm}
\proof{
Let $\dim X=n+1$. 
Let $J^n(X/\OO)$ be the jacobian ideal. Namely, this is the annihilator 
of $\Omega_{X/\OO}^{n+1}$. As $X$ is $\OO$-smooth outside $x$, 
$J^n(X/S)$ 
contains $\mathfrak{m}_x^N$, a power of 
the  ideal sheaf $\mathfrak{m}_x$ of $x$. 

Fix a uniformizer $t\in \OO$, hence identifications $\OO\cong k\{t\}$ and $\hat{\OO}\cong k[[t]]$. 
Since $\OO_{X,x}$ is regular, there is a presentation as a $k[[t]]$-algebra 
\begin{equation*}
\hat{\OO}_{X,x}\cong k[[t,t_0,\dots,t_{n}]]/(h). 
\end{equation*}
Let $M$ be an integer $\geq 2N+1$. Set 
$h_1$ to be the degree $\leq M$ part of $h$. As $J^n(X/\OO)\supset
\mathfrak{m}_x^N$, 
\cite[5.11]{App} implies that the natural isomorphism 
$\OO_{X,x}/\mathfrak{m}_x^M\cong k[t,t_0,\dots,t_{n}]/(h_1,(t,t_0,\dots,t_{n})^M)$ extends to an $\OO$-isomorphism 
\begin{equation*}
\OO_{X,(x)}\cong k\{t,t_0,\dots,t_{n}\}/(h_1). 
\end{equation*}

We take a polynomial $\Tilde{Q}\in k[t,t_0,\dots,t_n]$ of the following form depending on which cases we are dealing with: 
\begin{enumerate}
\item If ${\rm char} (k)$ is odd or $n+1$ is even, set $\Tilde{Q}=Q+t$ for some non-degenerate quadratic form 
$Q$ over $k$ in variables $t_0,\dots,t_n$. 
\item Otherwise, set $\Tilde{Q}=Q+t_0^2+tt_0+t$ for some non-degenerate quadratic form $Q$ in variables 
$t_1,\dots,t_n$ (note that $n$ is even in this case). 
\end{enumerate}

We construct a ``homotopy'' which connects $h_1$ and $\Tilde{Q}$, 
adding an extra coordinate $a$. 
We consider a commutative diagram 
\begin{equation*}
\xymatrix{
Y'\ar[rr]^-{\tilde{f}}\ar[rd]_-\pi&&\mathbb{A}^1_S={\rm Spec}(k[a,t])\ar[ld]\\
&S:={\rm Spec}(k[a])
}
\end{equation*}
where $Y'$ is defined for now by $Y':={\rm Spec}(k[a,t,t_0,\dots,t_n]/(h_1+a(\Tilde{Q}-h_1)))$. 
Let $x_0=(0,0,\dots,0),x_1=(1,0,\dots,0)$ be points of $Y'$. 
Since $Y'$ is Cohen--Macaulay and the fibers of $\pi$ over $a=0,1$ are of codimension $1$ in $Y'$, $\pi$ is flat over $a=0,1$. 
By the local criterions of smoothness and flatness, we know that $\pi$ is smooth and $\tilde{f}$ is flat at $x_0,x_1$. Shrinking $Y'$ 
around $x_0,x_1$, we may assume that $\pi$ is smooth and that $\tilde{f}$ is flat. 

Let $Z$ be the singular locus of $\tilde{f}$. Then $x_0,x_1$ are isolated in the fibers $Z\times_{S}\{a=0\}$ and 
$Z\times_S\{a=1\}$ respectively. Replacing $Y'$ by an open neighborhood 
around $x_0,x_1$ if necessary, we may assume that $Z$ is quasi-finite (by \cite[(13.1.3)]{EGA43}), hence flat, over $S$ and that $Z\times_{S}\{a=0\}$ consists of only one point $x_0$. We write $Y$ for such an open neighborhood. 

Let $Z_0$ be the closed subscheme of $Y$ defined by $t,t_0,\dots,t_n=0$. This defines a section of 
$Z\to S$ 
 which passes through $x_0$ and $x_1$. 
Note that $Z_0$ is an 
irreducible component of $Z$ as $Z$ is $S$-finite. 

Note that $x_1$ is an ordinary quadratic point of $\tilde{f}$ with 
the Milnor number $1$ or $2$. Therefore 
Lemma \ref{classquad} can be applicable to $Y_{(x_1)}\to S_{(1)}$. 
In particular, 
 the map $Z\to S$ is universally injective around $x_1$. 
 On the other hand, $Z_0$ gives a section of $Z\to S$. Hence 
 there exists an open neighborhood $U'$ of $x_1$ in $Z$ which is 
 contained in $Z_0$. 
 By \cite[1.3.4]{PL}, we may assume that 
 $U'$ only consists of ordinary quadratic points. 
Applying Proposition \ref{Milconti} to $U'$, we know that 
the Milnor number at {\it any} point in $U'$ equals to 
 $1$ or $2$, accordingly to the division into cases made above. 

After replacing $S$ by an irreducible \'etale neighborhood around $a=0$ and then shrinking $Y$ by an open neighborhood around 
$x_0$, we can make $S$, $Y$ affine, and $Z$ finite over $S$. As $Z_0$ is irreducible and intersects with $Z$ at $x_0$, 
 $U'$ remains non-empty after this replacement. 
 Then such $Y,\pi,\tilde{f}$ satisfy the conditions; we can take $U$ 
 in the condition $3$ to be the image of $U'$ by the map $Z\to S$. 
\qed}

This lemma gives us a following well-known result in characteristic $2$. 

\begin{pr}\label{tdeven}
Let $k$ be a perfect field of chracteristic $2$. 
Let $X$ be a smooth $k$-scheme of odd dimension $n$. 
Let $f\colon X\to C$ be a $k$-morphism to a smooth curve which has an isolated singular point at $x\in X$. Then 
$\mu(f,x)$ is even. 
\end{pr}
\proof{
Taking the base change to an algebraic closure, we may assume that 
$k$ is algebraically closed, hence $x$ is $k$-rational. 

Applying Lemma \ref{deform} to the map of the henselizations 
$X_{(x)}\to C_{(f(x))}$, we find a 
family (\ref{familydef}) as in the lemma. 
We follow the notation there. 
By the condition $2$ in the lemma, we need to show that the $\OO_S$-rank of $\OO_Z$ is even. 
Let $\bar{t}$ be a geometric generic point of $S$. 
Then the fiber $Z_{\bar{t}}$ consists of an ordinary quadratic point with the Milnor number $1$ or $2$. 
However, by 
Lemma \ref{classquad}, the Milnor number must be $2$ in this case. 
Therefore $Z_{\bar{t}}$ consists of points $z_1,\dots,z_m$ with $\mu(\tilde{f}_{\bar{t}},z_1)=2$. 
By the continuity (Proposition \ref{Milconti}), we have 
\begin{equation*}
\mu(f,x)=\sum_i\mu(\tilde{f}_{\bar{t}},z_i)=2+\sum_{i\neq1}\mu(\tilde{f}_t,{\bar{t}}_i). 
\end{equation*}
In particular, the Milnor numbers $\mu(\tilde{f}_{\bar{t}},z_i)$ for $i=2,\dots,m$ is strictly less than $\mu(f,x)$ 
and we conclude by induction on $\mu(f,x)$. 
\qed}

\section{$\mathbb{Z}/2$-coverings associated with isolated singularities}
\subsection{Discriminants of symmetric bilinear forms}

We recall that one can define the discriminants of non-degenerate symmetric bilinear forms as $\mu_2$-torsors. 
Let $S$ be a scheme. Let $(\varphi,B)$ be a pair of a locally free $\OO_S$-module $\varphi$ of finite rank and a non-degenerate 
symmetric bilinear form $B$ on $\varphi$. 
The form $B$ induces an isomorphism $\varphi\to
\varphi^\vee$. Taking the determinant, 
it induces an isomorphism 
$\det\varphi\to\det\varphi^\vee$, hence 
an isomorphism $\det\varphi^{\otimes2}\to\OO_S$, which 
is denoted by $\det B$.

Let $({\rm Sch}/S)$ denote the category of $S$-schemes. 
For an $S$-scheme $T$, we use the same symbol $\det B$ 
for the base change $(\det\varphi\otimes_{\OO_S}\OO_T)^{
\otimes2}\to\OO_T$ of $\det B$. 
{\it A presheaf} on $S$ means a functor $({\rm Sch}/S)^{\rm op}\to 
({\rm Set})$ to the category of sets. 

For a non-degenerate symmetric bilinear form $(\varphi,B)$ and an integer $N$, define 
${\rm disc}_NB$ to be the presheaf on $S$ sending an $S$-scheme $T$ to the set 
\begin{equation}\label{discdef}
{\rm disc}_NB(T):=\{x\in\det\varphi\otimes_{\OO_S}\OO_T\vert
\det B(x,x)=(-1)^N\}. 
\end{equation}

Let $\mu_2$ be the presheaf on $S$ sending 
an $S$-scheme $T$ to the group 
\begin{equation*}
\{a\in \Gamma(T,\OO_T)|a^2=1\}. 
\end{equation*}
This functor is represented by ${\rm Spec}(\OO_S[t]/(t^2-1))$ and is isomorphic to the kernel of $\mathbb{G}_{m,S}\xrightarrow{a\mapsto a^2}\mathbb{G}_{m,S}$. We consider the action of $\mu_2$ on 
${\rm disc}_NB$ defined by 
\begin{equation*}
\mu_2(T)\times {\rm disc}_NB(T)\ni(a,x)\mapsto ax\in {\rm disc}_NB(T). 
\end{equation*}
\begin{lm}\label{rep}
\begin{enumerate}
\item Assume that $\det\varphi$ is monogenic. 
Set $\alpha=\det B(x,x)$ for a basis $x\in\det\varphi$. 
Then 
${\rm disc}_NB$ is represented by ${\rm Spec}(\OO_S[u]/(u^2-(-1)^N\alpha))$. 
\item 
The presheaf ${\rm disc}_NB$ is representable by an $S$-scheme 
of finite flat of finite presentation. 
The action of $\mu_2$ gives ${\rm disc}_NB$ a structure of 
a $\mu_2$-torsor over $S$ in the fppf topology. 
\item Let $h\colon S'\to S$ be a morphism of schemes. Let $(\varphi',B')$ be the pullback of $(\varphi,B)$ by $h$. Then 
we have ${\rm disc}_NB'\cong {\rm disc}_NB\times_SS'$. 
\end{enumerate}
\end{lm}
\proof{
$1$.  
The isomorphism $\OO_S\to\det\varphi,1\mapsto x$ 
identifies $\det B$ with 
$\OO_S\otimes\OO_S\to\OO_S, a\otimes b\mapsto ab\alpha$. 
This identification gives the identification of ${\rm disc}_NB$ with the presheaf 
\begin{equation*}
 T\mapsto\{a\in\Gamma(T,\OO_T)|a^2=(-1)^N\alpha^{-1}\}, 
\end{equation*}
which is represented by ${\rm Spec}(\OO_S[u]/(u^2-(-1)^N\alpha))$. 

$2$. 
As ${\rm disc}_NB$ is a sheaf in the Zariski topology, we can work 
Zariski-locally on $S$. Thus we may assume that 
$\det\varphi$ is free of rank $1$. 
Then the assertion follows from 1. 

$3$. It is obvious from the definition. 
\qed}

The Kummer short exact sequence 
\begin{equation*}
0\to\mu_2\to\mathbb{G}_{m}\xrightarrow{a\mapsto a^2}
\mathbb{G}_{m}\to0
\end{equation*}
gives an injection $\Gamma(S,\mathbb{G}_m)/(\Gamma(S,\mathbb{G}_m))^2\to H^1(S,\mu_2)$. 
Lemma \ref{rep}.1 means that, if $\det\varphi$ has a global basis $x$, the $\mu_2$-torsor ${\rm disc}_NB$ is equal to the image of $\det B(x,x)$ under this injection. 

\begin{df}\label{mu2tor}
\begin{enumerate}
\item 
For a non-degenerate symmetric bilinear form $(\varphi,B)$ and $N\in\mathbb{Z}$, we call 
the presheaf ${\rm disc}_NB$ defined by (\ref{discdef}) the $N${\rm -signed discriminant} of 
$(\varphi,B)$. When $N=0$, we omit $N$ in the symbol and write 
${\rm disc}B$ for it, which is simply called the discriminant of $(\varphi,B)$. 
\item When $\det\varphi$ is monogenic, we also use the same symbol ${\rm disc}_NB$ (resp. 
${\rm disc}B$) for the 
element in $\Gamma(S,\mathbb{G}_m)/(\Gamma(S,\mathbb{G}_m))^2$ mapping to 
${\rm disc}_NB$ (resp. ${\rm disc}B$) under the injection $\Gamma(S,\mathbb{G}_m)/(\Gamma(S,\mathbb{G}_m))^2\to H^1(S,\mu_2)$. 
\end{enumerate}
\end{df}

\subsection{The case where $2$ is invertible}\label{2inv}

Let $S$ be a scheme. Let 
$X$ be a smooth $S$-scheme purely of relative dimension $n$ and $f\colon X\to C$ be an $S$-morphism 
to a smooth $S$-curve $C$. Assume that the singular locus $Z$ is finite over $S$. We also fix 
an everywhere non-zero differential $\omega\in\Omega^1_{C/S}$.

Assume that $2$ is invertible in $S$. In this case, $\mu_2$ is canonically isomorphic to the constant sheaf $\mathbb{Z}/2$ and 
we are allowed to give the following definition. 
\begin{df}\label{dtodd}
For an integer $N$, define the element 
\begin{equation*}
\rho_{f,\omega,N}^{\rm g}\in H^1(S,\mathbb{Z}/2)
\end{equation*}
to be the cohomology class corresponding to the $N$-signed discriminant 
${\rm disc}_NB_{f,\omega}$ (
Definition \ref{mu2tor}) under the identification 
$H^1(S,\mathbb{Z}/2)\cong H^1(S,\mu_2)$. 
\end{df}
When $S$ is noetherian and connected, $H^1(S,\mathbb{Z}/2)$ is identified with 
the group of continuous group homomorphisms 
$\pi_1^{ab}(S)\to\{\pm 1\}$. In this case we also 
use the same symbol $\rho_{f,\omega,N}^{\rm g}$ for the 
corresponding group homomorphism. 

\begin{lm}
Let the notation be as above. Assume that $2$ is invertible in $S$. 
\begin{enumerate}
\item For a morphism $h\colon S'\to S$ of schemes, 
we have 
\begin{equation*}
h^\ast\rho_{f,\omega,N}^{\rm g}=\rho_{f',\omega',N}^{\rm g}\hspace{2mm}{\rm in}
\hspace{1mm}
H^1(S',\mathbb{Z}/2). 
\end{equation*}
Here $f'$ and $\omega'$ are the pullbacks of $f$ and $\omega$ to 
$S'$. 
\item Assume that $\det\varphi_{f}$ 
is monogenic. 
Put $\alpha:= \det B_{f,\omega}(x,x)$ for a basis $x$. Then $\rho^{\rm g}_{f,\omega}$ 
is the character of the square root of $(-1)^N\alpha$. 
\end{enumerate}
\end{lm}
\proof{
$1$. This follows from the isomorphism of functors 
${\rm disc}_NB_{f',\omega'}\cong {\rm disc}_NB_{f,\omega}\times_SS'$ by Lemma \ref{rep}.3. 

$2$. This is a restatement of Lemma 
\ref{rep}.1 in this particular case. 
\qed}

\subsection{The case of characteristic 2}

In this subsection, we construct finite \'etale $\mathbb{Z}/2$-coverings from isolated singularities in characteristic $2$, which 
can be regarded as an {\it \'etale} analogue of the construction in the subsection 
\ref{2inv}.

As we consider a lift to Witt rings, we slightly change the notation as follows. Let $S_0$ be an $\F_2$-scheme, 
$X_0$ be a smooth $S_0$-scheme, 
and $f_0\colon X_0\to\A^1_{S_0}$ be an $S_0$-morphism. 
The singular locus $Z_0\subset X_0$ of $f_0$ is assumed to be finite 
over $S_0$. Write $t$ for the standard coordinate of $\A^1_{S_0}$. 
From these data, we construct a $\mathbb{Z}/2$-covering over $S_0$ in the case where $S_0,X_0$ are affine or 
$S_0$ is the spectrum of a field. 

First we treat the case where $S_0$ and $X_0$ are affine. Let $S^{{\rm perf}}_0$ be the perfection of $S_0$,~i.e.~the projective limit of the 
absolute Frobeniuses $\Phi\colon S_0\to S_0$
\begin{equation*}
\cdots\to S_0\xrightarrow{\Phi} S_0\xrightarrow{\Phi} S_0
\end{equation*}
indexed by non-negative integers. 
Note that the $0$-th projection $\pi_0\colon S^{{\rm perf}}_0\to S_0$ 
induces an equivalence between their \'etale topoi. Therefore we can replace $S_0$ by $S^{{\rm perf}}_0$ to construct finite \'etale coverings. 

Thus we further assume that $S_0$ is perfect from now on. 
Let $A:=\Gamma(S_0,\OO_{S_0})$. 
Let $S_n:={\rm Spec}(W_{n+1}(A))=
{\rm Spec}(W(A)/p^{n+1})$ be the spectrum of the ring of 
Witt vectors of length $n+1$. 

We collect some basic terminologies of lifts which are needed. 
\begin{df}\label{liftdf}
Let $m\geq n\geq0$ be integers. 
\begin{enumerate}
\item Let $X_n$ be a smooth $S_n$-scheme. 
 {\rm A lift} of $X_n$ to $S_m$ is a pair $(X_m,\iota_m)$ where $X_m$ is a smooth $S_m$-scheme together with an $S_n$-isomorphism 
$\iota_m\colon X_n\to X_m\times_{S_m}S_n$. 

\item Let $X_n,Y_n$ be smooth $S_n$-schemes and $f_n\colon X_n\to Y_n$ be an $S_n$-morphism. 
Let $(X_m,\iota_m),(Y_m,\kappa_m)$ be lifts of $X_n,Y_n$ to $S_m$. An $S_m$-morphism $f_m\colon X_m\to Y_m$ is 
called {\rm a lift} of $f_n$ 
if we have $f_m\circ\iota_m=\kappa_m\circ f_n$. 
\end{enumerate}
\end{df}

\begin{lm}(\cite{lisext})\label{lift}
Let $m\geq n\geq0$ be integers. Let $X_n$ be a smooth $S_n$-scheme. Assume that $X_n,S_n$ are affine. 
\begin{enumerate}
\item Let $(X_m,\iota_m)$ be a lift to 
$S_m$. 
Let $Y_n$ be another affine smooth $S_n$-scheme with a lift $(Y_m,\kappa_m)$ to $S_m$. Let $f_n\colon X_n\to Y_n$ be an $S_n$-morphism. Then  
there exists at least one lift of $f_n$ to $S_m$. If $f_n$ is an isomorphism, any 
lift is an isomorphism. 
\item There exists at least one lift of $X_n$ to $S_m$. Such lifts are isomorphic to each other. 
\end{enumerate}
\end{lm}
\proof{
$1$. 
The existence of a lift of $f_n$ follows from the smoothness of $Y_m$.  
The second assertion is a special case of \cite[4.2]{lisext}. 

$2$. This is \cite[6.7]{lisext}: since $X_n$ is affine, the obstructions to 
the existences of lifts always vanish. The latter assertion follows from $1$. 
\qed}

Let us go back to our situation. Hence 
$X_0$ is an affine smooth $S_0$-scheme and $f_0\colon 
X_0\to\A^1_{S_0}$ is an $S_0$-morphism with the singular locus $Z_0$ finite over $S_0$. By Lemma \ref{lift}, we take and fix a lift 
$X_2$ of $X_0$ to $S_2$ and a lift $f_2\colon 
X_2\to\A^1_{S_2}$ of $f_0$. 
Let $Z_2$ be the singular locus of $f_2$. By Lemma \ref{bcvan}.1, 
we know that $Z_2\times_{S_2}S_0$ is isomorphic to $Z_0$. In particular, $Z_2$ is finite over $S_2$. Applying the construction in 
Definition \ref{phiB}, we have a symmetric bilinear form 
$(\varphi_{f_2},B_{f_2,dt})$ on $W_3(A)$. 

Following is the main theorem in this section. 
\begin{thm}\label{z/2ext}
Let the notation be as above. Set 
$N=\frac{{\rm dim}(X_0/S_0)\mu(f_0,Z_0)}{2}$, which is an integer by 
Proposition \ref{tdeven}. 
\begin{enumerate}
\item{\rm (\textbf{\'etaleness})} The class  ${\rm disc}_NB_{f_2,dt}\in H^1(S_2,\mu_2)$ belongs to the image of 
the injective map 
$H^1(S_2,\mathbb{Z}/2)\to H^1(S_2,\mu_2)$ in Lemma \ref{mapz/2mu_2}.2. We write $\rho^{\rm g}_{f_2,dt}$ for the class in 
$H^1(S_2,\mathbb{Z}/2)$ mapping to ${\rm disc}_NB_{f_2,dt}$. 
\item{\rm (\textbf{compatibility with base change})} Let $h\colon S'_0\to S_0$ be a morphism of affine perfect schemes. Then we have 
\begin{equation*}
h^\ast\rho^{{\rm g}}_{f_2,dt}=\rho^{{\rm g}}_{f_2',dt}\hspace{3mm}{\rm in}\hspace{1mm}
H^1(S_2',\mathbb{Z}/2).  
\end{equation*}
Here $f_2'$ denotes the pullback of $f_2$. 
\item{\rm (\textbf{independence})} Assume that $A$ is normal. Then $\rho^{\rm g}_{f_2,dt}$ is independent of the choice of lifts $X_2,f_2$. 
\end{enumerate}
\end{thm}
\begin{rmk}\label{Wliftform}
Suppose that $X_2$ lifts to a formal ${\rm Spf}(W(A))$-scheme 
$\cal X$ formally smooth formally of finite type and that $f_2$ lifts 
to a morphism ${\bf f}\colon {\cal X}\to\hat{\A}^1_{W(A)}$ of formal  
$W(A)$-schemes where $\hat{\A}^1_{W(A)}$ is the formal 
$p$-adic completion of $\A^1_{W(A)}$. Such lifts can be taken using Lemma \ref{lift} inductively. For $n\geq0$, set 
$X_n={\cal X}\otimes_{W(A)}W_{n+1}(A)$ and 
$f_n={\bf f}\otimes_{W(A)}W_{n+1}(A)$. We have a projective system 
$\{(\varphi_{f_n},B_{f_n,dt})\}_n$ of 
bilinear forms on finite projective $W_{n+1}(A)$-modules. 
Its limit $(\varphi,B):=\varprojlim_n(\varphi_{f_n},B_{f_n,dt})$ gives a non-degenerate symmetric bilinear form on $W(A)$. 
Consider the $N$-signed discriminant 
${\rm disc}_NB$ 
defined in Definition \ref{mu2tor}.1 for $N=\frac{{\rm dim}(X_0/S_0)\mu(f_0,Z_0)}{2}$. 
Then, by Lemmas \ref{Wet}, \ref{mapz/2mu_2}.1,2, 
the statement $1$ in the theorem 
is equivalent to saying that the normalization of 
${\rm Spec}(W(A))$ in ${\rm disc}_NB\otimes_{W(A)}W(A)[\frac{1}{2}]$ is finite \'etale over $W(A)$. Indeed, the image of $\rho^{\rm g}_{f_2,dt}$ 
in $H^1(S_\infty[\frac{1}{2}],\mathbb{Z}/2)$ in the notation of Lemma 
\ref{mapz/2mu_2} is equal to ${\rm disc}_NB\otimes_{W(A)}W(A)[\frac{1}{2}]$. 
\end{rmk}

We generalize the definition of $\rho^{\rm g}_{f_2,dt}$ to 
an affine normal $\F_2$-scheme $S$ as follows. 
Here we say that
 an affine 
scheme $S$ is normal if $\Gamma(S,\OO_S)$ is the finite product of 
normal domains. 

Let $S$ be an affine normal $\F_2$-scheme. 
Let $f\colon X\to\A^1_{S}$ be a morphism 
of smooth affine $S$-schemes with the singular locus $Z$ finite over $S$. We write $S^{\rm perf}$ for the perfection 
of $S$. Applying Theorem \ref{z/2ext} to the base change of $f$ to $S^{\rm perf}$, we obtain a quadratic character 
$\rho^{\rm g}_{f_2,dt}$ of $S^{\rm perf}$, which is independent of the choice of lifts $f_2$ by Theorem \ref{z/2ext}.3.  
Since the \'etale topoi of $S$ and $S^{\rm perf}$ are canonically equivalent, the following definition is allowed. 
\begin{df}\label{char2}
Let $S$ be an affine normal $\F_2$-scheme. 
We define $\rho_{f,dt}^{\rm g}\in H^1(S,\mathbb{Z}/2)$ to be the cohomology class corresponding to $\rho^{\rm g}_{f_2,dt}$ for any lift $f_2$ of $f\times_SS^{{\rm perf}}$, via 
the equivalence of the \'etale topoi of $S$ and $S^{\rm perf}$. 
\end{df}
Recall that the Artin--Schreier exact sequence gives an isomorphism of groups 
\begin{equation}\label{ASiso}
A/\wp(A)\to H^1({\rm Spec}(A),\mathbb{Z}/2)
\end{equation}
 for an $\F_2$-algebra $A$. Using this, 
we can define a generalization of Arf invariant as follows. 
\begin{df}\label{arfex}
Let $S$ be an affine normal $\F_2$-scheme. 
Let $A$ denote $\Gamma(S,\OO_S)$. 
Let $X$ be a smooth $S$-scheme purely of dimension $n$. Let $f\colon X\to \A^1_S$ be an $S$-morphism whose 
singular locus $Z$ is finite over $S$. 
We define ${\rm Arf}(f,Z)$ to be the element in $A/\wp(A)$ which corresponds to 
$\rho_{f,dt}^{\rm g}$ in Definition \ref{char2} via the isomorphism (\ref{ASiso}). 
\end{df}
This terminology is justified by Corollary \ref{Arfquad}: in this corollary, we see that the invariant ${\rm Arf}(Q,0)$ so obtained from 
a non-degenerate quadratic form $Q$ over a perfect field of 
characteristic $2$ coincides with the Arf invariant 
of $Q$ in the usual sense. 

With an extra assumption that $\det\varphi_{f}$ is {\it monogenic}, we have a following explicit definition of ${\rm Arf}(f,Z)$ 
(compare Corollary \ref{Arfquad}). 
\begin{pr}\label{conarf}
Let the notation be as in Definition \ref{arfex}. Further assume that $A$ is perfect and that $\det\varphi_{f}$ is monogenic. 
Take a lift $f_2$ of $f$ to $W_3(A)$. Then, $\det\varphi_{f_2}$ is monogenic and we have the discriminant ${\rm disc}B_{f_2,dt}$ 
in $W_3(A)^\times/(W_3(A)^\times)^2$ (Definition \ref{mu2tor}.2). Then we have an equality 
\begin{equation*}
(-1)^N{\rm disc}B_{f_2,dt}\equiv 1+4[{\rm Arf}(f,Z)]
\end{equation*}
in $W_3(A)^\times/(W_3(A)^\times)^2$. Here $N=\frac{n\mu(f,Z)}{2}$. 
\end{pr}
The equality uniquely characterizes ${\rm Arf}(f,Z)$ by Lemma \ref{Wunr}.2. 
\proof{
We follow the notations in Lemma \ref{mapz/2mu_2}. 
Let $\alpha\in W(A)^\times$ be an element in the class 
${\rm disc}_NB_{f_2,dt}=(-1)^N{\rm disc}B_{f_2,dt}
\in 
W(A)^\times/(W(A)^\times)^2=W_3(A)^\times/(W_3(A)^\times)^2\subset H^1(S_2,\mu_2)$. Consider the sequence of maps 
\begin{equation*}
H^1(S_2,\mathbb{Z}/2)\to H^1(S_2,\mu_2)\to 
H^1(S_\infty[\frac{1}{2}],\mathbb{Z}/2)
\end{equation*}
in Lemma \ref{mapz/2mu_2}.2. As its composition $H^1(S_2,\mathbb{Z}/2)\to H^1(S_\infty[\frac{1}{2}],
\mathbb{Z}/2)$ can be identified with the canonical one $H^1(S_\infty,\mathbb{Z}/2)\to H^1(S_\infty[\frac{1}{2}],
\mathbb{Z}/2)$ by Lemma \ref{mapz/2mu_2}, the statement $1$ in Theorem \ref{z/2ext} implies that $W(A)[\frac{1}{2}][u]/(u^2-\alpha)$ extends to a 
finite \'etale $W(A)$-algebra, which is equivalent to saying that 
$\alpha$ is congruent to an element of the form $1+4[a]$ modulo $(W(A)^\times)^2$ by Lemma \ref{Wunr}.3. 
The equality $a= {\rm Arf}(f,Z)$ in $A/\wp(A)$ follows from Lemma \ref{Wunr}.3. 
\qed}

In the rest of this subsection, we 
explain how to reduce Theorem \ref{z/2ext} to the case of perfect fields. 
The compatibility $2$ follows from Lemma \ref{rep}.3. 
For the assertion $1$, we may replace $A$ by the residue field at 
each generic point by Lemma \ref{mapz/2mu_2}.3. 
The assertion $3$ in Theorem \ref{z/2ext} can 
be checked at the generic point $\eta$ of $A$ by Lemma \ref{mapz/2mu_2}.4.   

In this way, we reduce Theorem \ref{z/2ext} to the case where $A=k$ is a perfect field, which we treat in the next subsection. 

\subsection{The case of a field of characteristic $2$}

Let $k$ be a perfect field of characteristic $2$. Let $X_0$ be a smooth $k$-scheme purely of dimension $n$ and $f_0\colon X_0\to\A^1_k$ be a 
$k$-morphism. Assume that the singular locus of $f_0$ consists of {\it one point} $x$. Shrinking $X_0$ around 
$x$, we assume that $X_0$ is affine. 

Let $X_2$ be a lift of $X_0$ to $W_3(k)$ and $f_2\colon X_2\to\A^1_{W_3(k)}$ be a lift of $f_0$ 
(Definition \ref{liftdf}). 
Theorem \ref{z/2ext} follows from the following. 
\begin{thm}\label{z/2extk}
Set $N=\frac{n\mu(f_0,x)}{2}$. 
\begin{enumerate}
\item $(-1)^N{\rm disc}B_{f_2,dt}$ belongs to the image of $1+4W_3(k)\to W_3(k)^\times/(W_3(k)^\times)^2$. 
\item $(-1)^N{\rm disc}B_{f_2,dt}$ is independent of the choice of lifts $X_2$ and $f_2$. 
\end{enumerate}
\end{thm}
\proof{

First we treat the case where $x$ is $k$-rational. 
 In this case, 
we prove the assertion by induction on the Milnor number $\mu(f_0,x)$. If $\mu(f_0,x)=0$, there is nothing to prove as $f_0$ is smooth at $x$. 

Suppose that $\mu(f_0,x)>0$. Take a family 
\begin{equation*}
 \xymatrix{
Y\ar[rr]^{\tilde{f}}\ar[rd]_{\pi}&&\A^1_{S}\ar[ld]\\
 &S
 }
 \end{equation*}
which satisfies the conditions in Lemma \ref{deform} for 
$(X\to{\rm Spec}(\OO))=(X_{0,(x)}\to\A^1_{k,(f_0(x))})$. We follow the 
notation given there. 
In particular, we are given a $k$-rational point $s\in S$ such that the henselizations of the fibers $(Y\times_{S}s)_{(x_0)}\to 
\A^1_{k,(0)}$ is isomorphic to the map of the henselizations $f_0\colon X_{0,(x)}\to\A^1_{k,(f_0(x))}$. 
Let $Z$ be the singular locus of $\tilde{f}$. 

Let $S_0$ be the perfection of $S$. Note that $s\to S$ uniquely lifts to $s\to S_0$. 
Let us put the symbol $(-)_0$ to indicate base change by $S_0\to S$. Write $A$ for $\Gamma(S_0,\OO_{S_0})$. 
Take a lift $\tilde{f}_2\colon Y_2\to \A^1_{W_3(A)}$ to 
$W_3(A)$ of $Y_0$ and $\tilde{f}_0$.

First we show that the class ${\rm disc}_NB_{\tilde{f}_2,dt}\in H^1(S_2,\mu_2)$ is contained in the image of 
$H^1(S_2,\mathbb{Z}/2)\to H^1(S_2,\mu_2)$ and that 
it is independent of the choice of $\tilde{f}_2$. 
To do this, we may replace $S_0$ by its generic point $\eta$ by 
Lemmas \ref{mapz/2mu_2}.3,4. 
Hence we assume that $A$ is a perfect field. By the condition $3$ in Lemma \ref{deform}, $Z_0\times_{S_0}\eta$ 
consists of points $z_1,\dots,z_m$ where $z_1$ is an ordinary quadratic 
point with $\mu(\tilde{f}_0\times_{S_0}\eta,z_1)=1,2$. For the points 
$z_2,\dots,z_m$, we can apply the induction hypothesis. For 
$z_1$, the assertion is proved in Corollary \ref{Arfquad} 
if $\mu(\tilde{f}_0\times_{S_0}\eta,z_1)=1$. If 
$\mu(\tilde{f}_0\times_{S_0}\eta,z_1)=2$, it is proved in Proposition  \ref{compphi}.2. 

We prove the theorem in the case where $x$ is $k$-rational. The base changes $X_2:=Y_2\times_{W_3(A)}W_3(k)$, $ f_2:=\tilde{f}_2\times_{W_3(A)}W_3(k)$ 
by $s\to S_0$ 
give lifts of 
$X_0$, $f_0$ respectively. 
The assertion $1$ is already proved for these lifts. Thus it suffices to show that, 
for other lifts $X_2',f_2'$ of $X_0,f_0$, we have the equality 
${\rm disc}_NB_{f'_2,dt}={\rm disc}_NB_{f_2,dt}$. First, by Lemma \ref{lift}.2, we may replace $X_2'$ by $X_2$. 
Then, since the reduction of $f_2,f_2'$ to $k$ are the same, there exists a function $g\colon X_2\to\A^1_{W_3(k)}$ such that 
$f_2'=f_2+2g$. Take a map $\tilde{g}\colon Y_2\to\A^1_{W_3(A)}$ whose base change by $W_3(A)\to W_3(k)$ 
equals to $g$. Then $\tilde{f}_2+2\tilde{g}$ gives another lift of $\tilde{f}_0$. The assertion follows since we already know 
${\rm disc}_NB_{\tilde{f}_2,dt}={\rm disc}_NB_{\tilde{f}_2+2\tilde{g},dt}$. 

We reduce the case where $x$ is not necessarily $k$-rational to the case treated above. First we show that ${\rm disc}_NB_{f_2,dt}\in W_3(k)^\times/(W_3(k)^\times)^2$ is independent of the choice of $f_2$. 
Replacing $f_0$ and $f_2$ by $f_0\times_kk(x)$ and $f_2\times_{W_3(k)}W_3(k(x))$, we may assume that the structure map 
$X_2\to{\rm Spec}(W_3(k))$ factors through ${\rm Spec}(W_3(k(x)))$. 
Then the assertion follows from the case over $k(x)$ and Lemma \ref{trbil} below, as we have 
${\rm disc}_NB_{f_2,dt}=(-1)^N{\rm disc}B_{f_2,dt}$. 

We show that ${\rm disc}_NB_{f_2,dt}$ belongs 
to the image of $1+4W_3(k)$. This is equivalent to showing that 
$F[u]/(u^2-\alpha)$ is unramified where $F=W(k)[\frac{1}{2}]$. 
Hence we may take the base changes to an algebraic closure and reduce it to the case where $x$ splits into rational points, which case 
is already treated. 
\qed}

\begin{lm}\label{trbil}
Let $R$ be a commutative ring. Let $R'$ be an $R$-algebra which is 
finite free as an $R$-module. 
Let $(V,B')$ be a finite free $R'$-module $V$ with an $R'$-bilinear form $B'$. Let $B$ be the composition 
\begin{equation*}
V\times V\xrightarrow{B'}R'\xrightarrow{{\rm Tr}_{{R'/R}}}R, 
\end{equation*}
which we view as an $R$-bilinear form on the finite free $R$-module $V$. Then the discriminant of $(V,B)$ equals to 
${\rm disc}(R'/R)^{{\rm rk}_{R'}V}N_{R'/R}({\rm disc}(V,B'))$ in $R^\times/(R^\times)^2$. Here ${\rm disc}(R'/R)$ denotes the discriminant of the form ${\rm Tr}_{{R'/R}}$. 
\end{lm}
\proof{
See \cite[Section 2, Proposition 9]{bour}. 
\qed}

We note that the Arf invariant is {\it trivial} if the variety is of odd dimension, as the following proposition shows. 
\begin{pr}\label{arfodd}
Let $k$ be a perfect field of characteristic $2$. 
Let $f\colon X\to\A^1_k$ be a morphism of smooth $k$-schemes with an isolated singularity $x\in X$. Assume that 
${\rm dim}X$ is odd. Then we have ${\rm Arf}(f,x)=0$. 
\end{pr}
\proof{
We prove the assertion by induction on $\mu(f,x)$. 

First we reduce it to the case where $x$ is $k$-rational. Replacing $X$ by $X\times_kk(x)$ and $x$ by the image of the diagonal $x\to X\times_kk(x)$, we may 
assume that $X\to{\rm Spec}(k)$ factors through $k(x)$. Then the assertion for $(X,x)$ over $k$ follows from that for $(X\times_kk(x),x)$ over $k(x)$ by Lemma \ref{trbil}. Note that, as the Milnor number is even in this case, the part ${\rm disc}(R'/R)^{{\rm rk}_{R'}V}$ in the lemma is trivial.

Take a family (\ref{familydef}) as in Lemma \ref{deform}. We follow the notation given there. 
Then, by Definition \ref{char2}, we have a character 
$\rho^{\rm g}_{\tilde{f},dt}\colon \pi_1^{ab}(S)\to\{\pm1\}$ such that, for a 
morphism $h\colon{\rm Spec}(k')\to S$ from the spectrum of a perfect field, the composition 
$\rho^{\rm g}_{\tilde{f},dt}\circ h_\ast$ corresponds to $\sum_{y\in Z_{k'}}{\rm Arf}(\tilde{f}_{k'},y)$ via the Artin--Schreier theory. 

We take such a $k'$ as the perfection of the generic point of $S$. By the condition $3$ in Lemma \ref{deform}, 
$Z_{k'}$ consists of $z_1,\dots,z_m$ where $z_1$ is ordinary quadratic of the Milnor number $1$ or $2$. 
However, as ${\rm dim}X$ is odd, it cannot be $1$ and we can apply 
Proposition \ref{compphi}.2 to $z_1$. 
For $z_2,\dots,z_m$, the Milnor numbers are strictly less than $\mu(f,x)$ as we have 
$\mu(f,x)=\sum_{i}\mu(\tilde{f}_{k'},z_i)$. Hence the induction hypothesis is applied. 
\qed}

Combining this proposition with Theorem \ref{mainMil} in the next section, we see that the local epsilon factors of isolated singularities in odd dimension are always trivial. This triviality is also deduced from  
the continuity of local epsilon factors (Theorem \ref{dta}) together with 
Lemmas \ref{deform}, \ref{epex}, without referring to the Arf invariants.

\section{Milnor  formula for local epsilon factors}\label{milep}

Throughout this section, we fix a perfect field $k$ 
of characteristic $p>0$ and a prime number $\ell$ different from $p$. We fix an algebraic closure $\overline{\Ql}$ of $\Ql$. 
We also fix a non-trivial additive character $\psi\colon\mathbb{F}_p\to\overline{\Ql}^\times$, which 
we use to define local epsilon factors as in \cite[(3.1.5.4)]{Lau}.  

We give some notation which is used in this section. 
\begin{itemize}
\item 
For a field $k$, we write $G_k$ for its absolute Galois group. 
\item Let $T$ be a scheme and $t\in T$ be a point. 
For a morphism of schemes $z\to t$ which comes from a finite 
separable field extension of $k$, we write 
$T_{(z)}$ for the unramified extension of the henselization $T_{(t)}$ 
 whose residue field is $k(z)$. 
\end{itemize}

\subsection{Preliminaries on local epsilon factors}

Let $T$ be a henselian trait which is isomorphic to the henselization of $\mathbb{A}^1_{k}$ at the origin. 
Let $s$ and $\eta$ be the closed point and the generic point of $T$. 
Fix a non-zero rational differential $\omega\in\Omega^1_{k(\eta)}$. 
In \cite{Y3}, \cite{geomep}, Yasuda and Guignard independently give 
generalizations of the
 theory of local epsilon factors 
to general perfect residue fields. 
They attach, in a canonical way, a character 
$\varepsilon_{0,\bar{k}}(T,V,\omega)\colon G_k^{ab}\to\overline{\Ql}^\times$ of the absolute Galois group of the residue field, to a 
finite dimensional $\Ql$-representation $V$ of $G_{k(\eta)}$. Their theories coincide by \cite[8.3]{Y2} (for finite field cases) and \cite[11.8]{geomep}. In this paper, 
we choose the settings and notation for local epsilon factors 
similar to \cite{geomep}, as it fits our purposes well. 
An explanation to translate the results in \cite{Y3} to 
our settings is given in \cite[3.1]{contiep}.

We recall the relation of their theories with the classical one. 
When $k$ is finite, the Galois group $G_k$ is topologically generated by the geometric Frobenius ${\rm Frob}_k$. The value $(-1)^{{\rm dimtot} V}\varepsilon_{0,\bar{k}}(T,V,\omega)({\rm Frob}_k)$ 
equals to the classical local epsilon factor in \cite{Del}, \cite[(3.1.5.4)]{Lau}, 
where ${\rm dimtot}$ denotes the sum of 
the Swan conductor and the actual dimension of $V$. 
This equality follows from the cohomological interpretation by Laumon 
\cite[(3.5.1.1)]{Lau} and \cite[8.3]{Y2}, \cite[11.8]{geomep}. 

We usually identify in the sequel the absolute Galois group of $\eta$ and that of the completion with respect to the discrete valuation. 
The identification depends on the choice of the embedding of an algebraic closure of $k(\eta)$ into that of the completion. 
When we use this identification, we fix one of such embeddings. 
Since such an isomorphism is unique up to the conjugation, 
local epsilon factors do not depend on its choice. 

To recall the continuity of local epsilon factors in \cite{contiep}, 
we give necessary notation. 
Let $T$ be a henselian trait with the residue field $k$ as above. Recall that, for a finite 
separable field extension $k(z)/k$, we write 
$T_{(z)}$ for the unramified extension of the henselization of $T$ whose residue field is $k(z)$. We also write $\eta_z$ for the generic point of $T_{(z)}$.

Let $X$ be a $T$-scheme of finite type and let 
$f\colon X\to T$ be the structure morphism. 
For a constructible complex $\mathcal{F}$ of $\Ql$-sheaves on $X$, the vanishing cycles complex 
$R\Phi_f(\mathcal{F})$ is defined to be a constructible complex 
on the geometric special fiber $X_{\bar{s}}$. Let 
$Z\subset X$ be the closed subset outside of which 
$f$ is 
universally locally acyclic relatively to $\cal F$. Then the complex $R\Phi_f(\mathcal{F})$ 
is supported on $Z_{\bar{s}}$ and it admits an action of the absolute Galois group $G_{\eta}$ equivariantly to the action on $Z_{\bar{s}}$. 
Therefore, for a closed point $z\in Z$, 
the restriction of $R\Phi_f(\mathcal{F})$ to $z\times_s\bar{s}\subset 
Z_{\bar{s}}$ gives a bounded complex of finite dimensional 
$\Ql$-representations of $G_{\eta_z}$, which 
we denote by $R\Phi_f(\mathcal{F})_z$. 

Let $S$ be a noetherian $\mathbb{F}_p$-scheme. Let 
\begin{equation*}
\xymatrix{
Z\ar@{^{(}->}[r]&X\ar[rr]^f\ar[rd]^\pi&&\A^1_S\ar[ld]\\
&&S}
\end{equation*}
be a family of isolated singularities. Namely, 
\begin{enumerate}
\item $Z$ is the singular locus of $f$ which we assume finite over $S$. 
\item $\pi$ is smooth. 
\end{enumerate}
For a finite separable extension $k'/k$ of fields, we denote by ${\rm tr}_{k'/k}\colon G_k^{ab}
\to G_{k'}^{ab}$ the transfer morphism induced  
by the inclusion $G_{k'}\hookrightarrow G_k$. The determinant character of the 
induced representation ${\rm Ind}_{G_{k'}}^{G_k}
1_{G_{k'}}$ of the trivial representation is denoted by $\delta_{k'/k}$. 
The following is an arithmetic analogue to Proposition \ref{Milconti}. 
\begin{thm}(A special case of \cite[4.8]{contiep})\label{dta}
Let the notation and assumption be as above. 
Assume that $S$ is connected. From these data 
given above, we can 
attach a character
\begin{equation*}
\rho_{f,dt}^{{\rm a}}\colon \pi_1(S)^{ ab}\to\overline{\Ql}^\times
\end{equation*}
in such a way that the following hold. 
\begin{enumerate}
\item The formation of $\rho_{f,dt}^{\rm a}$ commutes with base change 
$S'\to S$. 
\item When $S$ is the spectrum of a perfect field $k$, we have 
\begin{equation*}
\rho_{f,dt}^{\rm a}=\prod_{z\in Z}\delta_{k(z)/k}^{{\rm dimtot}(R\Phi_f(\Ql)_z)}
\cdot\varepsilon_{0,\bar{k}}(\A^1_{k,(z)},R\Phi_f(\Ql)_z,dt)\circ{\rm tr}_{k(z)/k}. 
\end{equation*}
\end{enumerate}
When $S$ is noetherian normal, the character $\rho_{f,dt}^{\rm a}$ is uniquely determined by these properties. 
\end{thm}

We also recall the compatibility with additive convolution. 
 Let $f_1\colon X_1\to\mathbb{A}^1_k$ 
and $f_2\colon X_2\to\mathbb{A}^1_k$ be $k$-morphisms of finite type. 
For each $i=1,2$, we assume that $X_i$ is smooth over $k$ and that 
isolated singular points $x_i\in X_i$ of $f_i$ are given. We further 
assume that $x_i$ are $k$-rational and that the image $f_i(x_i)$ equal 
to the origin. 
Let $X:=X_1\times_kX_2$ and denote by 
$f\colon X\to\mathbb{A}^1_k$ the map sending 
$(x,y)\mapsto f_1(x)+f_2(y)$. 
\begin{lm}(\cite[3.18]{charep})\label{epconv}
Let the notation be as above. Let $x:=(x_1,x_2)\in X$ be the $k$-rational 
point over $x_1$ and $x_2$.
Then we have the equality 
\begin{align*}
\varepsilon_{0,\bar{k}}(\A^1_{k,(0)},R\Phi_{f}&
(\Ql)_x,dt)^{-1}=\\
&\varepsilon_{0,\bar{k}}(\A^1_{k,(0)},R\Phi_{f_1}(\Ql)_{x_1},dt)^{{\rm dimtot} R\Phi_{f_2}
(\Ql)_{x_2}}\cdot\varepsilon_{0,\bar{k}}(\A^1_{k,(0)},
R\Phi_{f_2}(\Ql)_{x_2},dt)^{{\rm dimtot} R\Phi_{f_1}(\Ql)_{x_1}}
\end{align*}
as characters of $G_k$. 
\end{lm}

\subsection{Main result}

To state the main theorem, we introduce some notations. 

We are fixing a non-trivial character $\psi\colon\mathbb{F}_p\to\overline{\Ql}^\times$. For a finite extension $k/\mathbb{F}_p$, 
we write $\psi_k$ for $\psi\circ{\rm Tr}_{k/\mathbb{F}_p}$. Let $\wp\colon X=\A^1_{\mathbb{F}_p}\to\A^1_{\mathbb{F}_p}$ 
be the morphism defined by $t\mapsto t^p-t$. This is the Artin--Schreier covering, whose Galois group 
${\rm Aut}(X/\A^1_{\mathbb{F}_p})$ is canonically 
isomorphic to $\mathbb{F}_p$. 
The composition $\pi_1(\A^1_{\mathbb{F}_p})^{ab}\to{\rm Aut}(X/\A^1_{\mathbb{F}_p})\cong
\mathbb{F}_p\xrightarrow{\psi^{-1}}\overline{\Ql}^\times$ 
gives a smooth $\overline{\Ql}$-sheaf 
$\mathcal{L}_\psi$ on $\A^1_{\mathbb{F}_p}$. Let $k$ be a finite extension of $\mathbb{F}_p$ and let 
$a\in\A^1_{\mathbb{F}_p}(k)\cong k$ be a $k$-rational point. Then the geometric Frobenius of 
$G_k$ acts on the stalk $\mathcal{L}_{\psi,\bar{a}}$ by the multiplication of $\psi_k(a)$. 

To the additive character $\psi_k$, we attach the quadratic Gauss sum 
$\tau_{\psi,k}$ by the following 
\begin{equation*}
\tau_{\psi,k}:=-\sum_{a\in k}\psi_k(a^2)=-\sum_{a\in k^\times}\bigl(\frac{a}{k}\bigr)\psi_k(a). 
\end{equation*}
Here the symbol 
$(\frac{}{k})$ denotes the Legendre symbol. Note that 
its square $\tau_{\psi,k}^2$ equals to 
$(\frac{-1}{k})q$ where $q$ is the cardinality of $k$. 

\begin{lm}\label{quadchar}
Assume that $p$ is odd. Let 
$f\colon\A^1_{\mathbb{F}_p}\to\A^1_{\mathbb{F}_p}$ be a morphism 
of schemes defined by 
$t\mapsto t^2$. Then 
the \'etale cohomology group $H^i(\A^1_{\overline{\mathbb{F}_p}},f^\ast\mathcal{L}_\psi)$ 
vanishes except in degree $1$ and $H^1$ is a one dimensional $\overline{\Ql}$-vector space 
on which the geometric Frobenius of any finite extension $k/\mathbb{F}_p$ 
acts by the multiplication of $\tau_{\psi,k}$. 
\end{lm}
For a connected noetherian $\mathbb{F}_p$-scheme $S$, we write $\rho_\psi$ for the composition of the map 
$\pi_1^{ab}(S)\to\pi_1^{ab}({\rm Spec}(\mathbb{F}_p))$ and the character 
$\pi_1^{ab}({\rm Spec}(\mathbb{F}_p))\to\overline{\Ql}^\times$ corresponding to 
$H^1(\A^1_{\overline{\mathbb{F}_p}},f^\ast\mathcal{L}_\psi)$. 
\proof{
Let $G\subset {\rm Aut}(X_{\overline{\mathbb{F}_p}}/\mathbb{A}^1_{\overline{\mathbb{F}_p}})$ be the image of 
$\pi_1^{ab}(\mathbb{A}^1_{\overline{\mathbb{F}_p}})\xrightarrow{f_\ast}\pi_1^{ab}(\mathbb{A}^1_{\overline{\mathbb{F}_p}})\to 
{\rm Aut}(X_{\overline{\mathbb{F}_p}}/\mathbb{A}^1_{\overline{\mathbb{F}_p}})$. The quotient ${\rm Aut}(X_{\overline{\mathbb{F}_p}}/\mathbb{A}^1_{\overline{\mathbb{F}_p}})/G$ is a $2$-group since ${\rm Coker}(f_\ast)\cong\mathbb{Z}/2$ surjects onto it. Since it is also a quotient of 
${\rm Aut}(X_{\overline{\mathbb{F}_p}}/\mathbb{A}^1_{\overline{\mathbb{F}_p}})\cong\mathbb{F}_p$, we have $G={\rm Aut}(X_{\overline{\mathbb{F}_p}}/\mathbb{A}^1_{\overline{\mathbb{F}_p}})$. This shows that $H^0$ vanishes. 
For higher $H^i$, it vanishes in degree $\geq2$ as $\mathbb{A}^1_{\overline{\mathbb{F}_p}}$ is affine 
and $1$-dimensional. 

The assertion on the dimension of $H^1$ comes from the Grothendieck--Ogg--Shafarevich formula applied 
to $f^\ast\mathcal{L}_\psi$, as the Swan conductor of the 
covering $X\times_{\A^1_{\mathbb{F}_p}f}\mathbb{A}^1_{\mathbb{F}_p}\xrightarrow{{\rm pr}_2}\A^1_{\mathbb{F}_p}$ at 
the infinity 
is $2$. The last assertion follows from the Lefschetz trace formula. 
\qed}

Here is the main theorem of this article. 

\begin{thm}\label{mainMil}
Let $S$ be a noetherian affine connected normal $\mathbb{F}_p$-scheme and let 
\begin{equation*}
\xymatrix{
Z\ar@{^{(}->}[r]&X\ar[rr]^f\ar[rd]_\pi&&\A^1_S\ar[ld]\\
&&S}
\end{equation*}
be a commutative diagram of $S$-schemes of finite type where $\pi$ is smooth purely of relative 
dimension $n$ and $Z$ is the singular locus of $f$. We assume that 
$Z$ is finite over $S$. 

In this case, the character $\rho_{f,dt}^{\rm a}\colon 
\pi_1(S)^{ab}\to\overline{\Ql}^\times$ is given in 
Theorem \ref{dta}. 
Let $\mu(f,Z)$ be the Milnor number along $Z$ (Definition \ref{Mildef}). 
\begin{enumerate}
\item Assume that $p$ is odd. Let $\rho_\psi$ be as in Lemma \ref{quadchar}. 
Then we have an equality 
\begin{equation*}
\rho_{f,dt}^{\rm a}=\rho_{f,-2dt,0}^{\rm g}\cdot \rho_\psi^{(-1)^{n+1}n\mu(f,Z)}
\end{equation*}
of characters $\pi_1(S)^{ab}\to\overline{\Ql}^\times$. 
For the definition of $\rho^{\rm g}_{f,-2dt,0}$, see 
Definition \ref{dtodd}.
\item Assume that $p=2$. Let $\chi_{\rm cyc}$ denotes  the  character of the $\ell$-adic Tate twist, so that 
the value at the geometric Frobenius equals to $q^{-1}$ if 
$S$ is the spectrum of a finite field with $q$ elements. We have 
\begin{equation*}
\rho_{f,dt}^{\rm a}=\rho_{f,dt}^{\rm g}\cdot\chi_{\rm cyc}^{\frac{(-1)^{n}n\mu(f,Z)}{2}}. 
\end{equation*}
For the definition of $\rho_{f,dt}^{\rm g}$, see Definition 
\ref{char2}. 
\end{enumerate}
\end{thm}

As a corollary, we have the following result, as promised in the introduction. Recall our convention: 
$\varepsilon_0$ denotes the classical local epsilon factors and we have 
$\varepsilon_{0,\bar{k}}({\rm Frob}_k)=(-1)^{{\rm dimtot}}\varepsilon_0$. 
\begin{cor}\label{Milcor}
Let $k$ be a finite field. Let $X$ be a smooth $k$-scheme of dimension $n$ and let $f\colon X\to \A^1_k$ be 
a $k$-morphism with an isolated singular point $x\in X$. 
We write $\A^1_{k,(x)}$ for the unramified extension of 
$\A^1_{k,(f(x))}$ with the residue field $k(x)$. 
\begin{enumerate}
\item Assume that $k$ is of odd characteristic. Then we have 
\begin{align*}
(-1)^{[k(x):k]{\rm dimtot}R\Phi_f(\Ql)_x}\varepsilon_0(\A^1_{k,(x)},R\Phi_f(\Ql)_x,dt)&=\Bigl(\frac{{\rm disc}B_{f,-2dt}}{k}\Bigr)\cdot
\tau_{\psi,k}^{(-1)^{n+1}n\mu(f,x)}\\
&=\Bigl(\frac{(-2)^{n\mu(f,x)}{\rm disc}B_{f,dt}}{k}\Bigr)\cdot\tau_{\psi,k}^{(-1)^{n+1}n\mu(f,x)}. 
\end{align*}
\item Assume that $k$ is of characteristic $2$. Write $q$ for the cardinality of $k$. 
Then the ratio 
\begin{equation*}
(-1)^{[k(x):k]{\rm dimtot}R\Phi_f(\Ql)_x}\varepsilon_0(\A^1_{k,(x)},R\Phi_f(\Ql)_x,dt)/q^{\frac{(-1)^{n+1}n\mu(f,x)}{2}}
\end{equation*}
is $\pm1$. This is $1$ if and only if ${\rm Arf}(f,x)\in \wp(k)$. 
\end{enumerate}
\end{cor}
\proof{
The statements are the special case of Theorem \ref{mainMil} where $S={\rm Spec}(k)$. The equality ${\rm disc}B_{f,-2dt}=(-2)^{n\mu(f,x)}{\rm disc}B_{f,dt}$ 
follows from Lemma \ref{bpofbil}.4. 
\qed}

Actually the proof of Theorem \ref{mainMil} is reduced to the finite field case. We prepare an auxiliary lemma, by which  
we can freely assume that a singular point in consideration is $k$-rational. Before doing so, we recall a well-known result on trace forms. For a finite 
separable field extension $L/F$, write $\delta_{L/F}$ for the 
determinant character of the 
induced representation ${\rm Ind}_{G_{L}}^{G_F}
1_{G_{L}}$ of the trivial representation. 
\begin{lm}\label{discdelta}
Let $F$ be a field in which $2$ is invertible. For a finite separable field 
extension $L/F$, write ${\rm disc}(L/F)$ for the discriminant 
of the non-degenerate quadratic form $(L,{\rm Tr}_{L/F}(x^2))$ over $F$. Then the quadratic character 
$\delta_{L/F}$ is defined by the square root of ${\rm disc}(L/F)$. 
\end{lm}
\proof{
Write $\mu_2=\{1,-1\}$ for the group of order $2$. Let $n=[L:F]$. 
Take a separable closure $F_s$ of $F$. The 
$n$-th symmetric group $\mathfrak{S}_n$ is naturally identified with 
the automorphism group of the \'etale $F_s$-algebra $E:=F_s^{\oplus n}$. As 
$L$ is an $F$-form of $E$, it gives an element 
$x$ in the Galois cohomology $H^1(F_s/F,\mathfrak{S}_n)$ (\cite[IX, Section 2]{Corp}). The image 
of $x$ under the map $H^1(F_s/F,\mathfrak{S}_n)\to H^1(F_s/F,\mu_2)$ defined by the signature is nothing but the 
character $\delta_{L/F}$. On the other hand, the trace form 
$(L,{\rm Tr}_{L/F}(x^2))$ is an $F$-form of 
$(E,{\rm Tr}_{E/F_s}(x^2))$. As the latter is equal to the standard quadratic form $x_1^2+\cdots+x_n^2$, $(L,{\rm Tr}_{L/F}(x^2))$ gives an element $y$ in $H^1(F_s/F,O_n)$, where $O_n$ denotes the orthogonal group. 
The image of $y$ under the map $H^1(F_s/F,O_n)\to 
H^1(F_s/F,\mu_2)$ defined by the determinant is equal to the character of the square root of 
${\rm disc}(L/F)$. Then the assertion follows from the observation that 
$x\in H^1(F_s/F,\mathfrak{S}_n)$ maps to $y\in 
H^1(F_s/F,O_n)$ via the map defined by the map 
$\mathfrak{S}_n\cong{\rm Aut}(E/F_s)\to O_n$. 
\qed}

\begin{lm}\label{k-ratred}
Consider the situation in Corollary \ref{Milcor}. Let $f_{x}\colon X_{k(x)}=X\times_kk(x)\to\A^1_{k(x)}$ be the base change of $f$ by 
$k\to k(x)$. We lift $x$ to a closed point of $X_{k(x)}$ by the diagonal map $x\to X_{k(x)}=X\times_kx$. Then, 
the statements in Corollary \ref{Milcor} for $(f,x)$ over $k$ are equivalent to the statements for $(f_x,x)$ over $k(x)$. 
\end{lm}
\proof{
We have $\mu(f,x)=\mu(f_x,x)[k(x):k]$ and $\tau_{\psi,k}^{[k(x):k]}=\tau_{\psi,k(x)}$ when $k$ is of odd characteristic (cf. Lemma \ref{quadchar}). As the local epsilon factors does not depend on the 
base fields, it remains to check the equality of signs. 

If the characteristic is odd, it remains to show that 
\begin{equation*}
(-1)^{([k(x):k]-1)\mu(f_x,x)}\Bigl(\frac{{\rm disc}B_{f_x,-2dt,k(x)}}{k(x)}\Bigr)=\Bigl(\frac{{\rm disc}B_{f,-2dt,k}}{k}\Bigr). 
\end{equation*}
In the left-hand side, we have $(-1)^{[k(x):k]-1}=\bigl(\frac{{\rm disc}(k(x)/k)}{k}\bigr)$ by Lemma \ref{discdelta}. 
Then the assertion follows once we verify the equality 
\begin{equation}\label{NTR}
{\rm disc}(k(x)/k)^{\mu(f_x,x)}N_{k(x)/k}({\rm disc}B_{f_x,-2dt,k(x)})=
{\rm disc}B_{f,-2dt,k}\hspace{4mm}{\rm in}\hspace{2mm} k^\times/(k^\times)^2. 
\end{equation}
If the characteristic is $2$, lifting $f$ to $W(k)$, the assertion is reduced to the similar equality as (\ref{NTR}) with 
$k$ replaced with the fraction field of $W(k)$ and $-2dt$ replaced with $dt$. 
By Lemma \ref{etres}, the equalities are consequences of Lemma \ref{trbil}. 
\qed}

(Proof of Theorem \ref{mainMil})

{
As the formations of $\rho_{f,dt}^{\rm a}, 
\rho^{\rm g}_{f,-2dt,0},$ and $\rho_{f,dt}^{\rm g}$ commute with base change, we may assume that 
$S$ is of finite type over $\mathbb{F}_p$. By Cebotarev density theorem, it suffices to show the 
equality of the values of the characters at the geometric Frobenius of every closed point. Again by the commutativity of the 
formations, we reduce it to the case where $S$ is the spectrum of a finite field $k$. We also assume that $Z$ consists of a 
$k$-rational point $x$ by Lemma \ref{k-ratred} and that $x$ maps to 
the origin of $\A^1_k$. 

We prove the statement by induction on $\mu(f,x)$. 
Applying Lemma \ref{deform} to $(f,x)$, we find a commutative diagram of smooth $k$-schemes 
\begin{equation*}
\xymatrix{
Y\ar[rr]^{\tilde{f}}\ar[rd]&&\A^1_{S'}\ar[ld]\\
&S'
}
\end{equation*}
such that 
\begin{enumerate}
\item $Y\to S'$ is smooth. The singular locus $Z'$ of $\tilde{f}$ is finite over $S'$. 
\item There is a $k$-rational point $s\in S'(k)$ such that $Z'_s=Z'\times_{S'}s$ consists of one point 
$x_0$ mapping to the origin by $\tilde{f}$ and the map of the henselizations of the fibers 
$\tilde{f}_s\colon Y_{s,(x_0)}\to\A^1_{k,(0)}$ is isomorphic to 
$X_{(x)}\to \A^1_{k,(0)}$. 
\item There is an open dense subset $U'\subset S'$ such that, for any closed point $t\in U'$, 
the fiber $Z'_t$ has a $k$-rational ordinary quadratic point $z_1$ with $\mu(\tilde{f}_t,z_1)=1,2$. 
\end{enumerate}
By Cebotarev density, it is enough to show the equlity on the geometric Frobeniuses at points on $U'$. 
Let $t\in U'$ be a closed point and let $z_1,\dots,z_m$ be the points in $Z'_t$ where $z_1$ is such an ordinary quadratic 
singularity. For $z_i\hspace{1mm}(i=2,\dots,m)$, we apply the induction hypothesis on Milnor number. 

It remains to treat the case where $x$ is $k$-rational and ordinary quadratic with $\mu(f,x)=1,2$. 
By Lemma \ref{classquad}, the henselization $\OO_{X,(x)}$ is isomorphic to the following as $k\{t\}$-algebras. 
\begin{enumerate}
\item If $p$ is odd or $n$ is even, \begin{equation*}
\OO_{X,(x)}\cong k\{t,t_0,\dots,t_{n-1}\}/(Q-t)
\end{equation*} for some non-degenerate quadratic form $Q\in k[t_0,\dots,t_{n-1}]$. 
\item If $p=2$ and $n$ is odd, 
\begin{equation*}
\OO_{X,(x)}\cong k\{t_1,\dots,t_{n-1}\}\otimes k\{u\}
\end{equation*}
which is regard as $k\{t\}$-algebra by the map 
$ t\mapsto Q\otimes1+1\otimes\tilde{t}$ where $Q\in k[t_1,\dots,t_{n-1}]$ is a non-degenerate quadratic form and $\tilde{t}$ satisfies 
$u^2+b(\tilde{t})u-\tilde{t}=0$ for some uniformizer $b\in k\{t\}$. 
\end{enumerate}
Up to changing the coordinates, we can take $Q$ to be the form $a_0t_0^2+\cdots+a_{n-1}t_{n-1}^2\hspace{3mm}(a_i\in k^\times)$ 
if $p$ is odd, and $(t_0^2+t_0t_1+a_0t_1^2)+(t_2^2+t_2t_3+a_2t_3^2)+\cdots\hspace{3mm}(a_i\in k)$ if $p=2$. 
By the compatibility with additive convolution (Lemmas \ref{prodphi}, \ref{epconv}), we may assume one of the following: 
\begin{enumerate}
\item $\OO_{X,(x)}\cong k\{t,t_0\}/(at_0^2-t)$ for $a\in k^\times$. 
\item $\OO_{X,(x)}\cong k\{t,t_0,t_1\}/(t_0^2+t_0t_1+at_1^2-t)$ for $a\in k$. 
\item $\OO_{X,(x)}\cong k\{t,u\}/(u^2+bu-t)$ for a uniformizer $b\in k\{t\}$. 
\end{enumerate}
Each case is treated in Proposition \ref{compphi} and Lemma \ref{epex}. 
The proof is completed. 
\qed}

\subsection{Example: discriminants of homogeneous polynomials}

At the end of this paper, we compute the discriminant of the bilinear form associated with a homogeneous polynomial $F$ in 
terms of its divided discriminant \cite{hypdet}. 

To start with, we recall a relation between conditions on smoothness of a homogeneous polynomial. 
\begin{lm}\label{condhom}
 Let $F\in k[T_0,\dots,T_{n+1}]$ be a homogeneous polynomial of degree $d>1$ with coefficients in a field $k$. Consider 
the following conditions. 
\begin{enumerate}
\item The map $F\colon\A^{n+2}_k\to\A^1_k$ has an isolated singular point at the origin. 
\item The map $F\colon\A^{n+2}_k\to\A^1_k$ is smooth outside the origin. 
\item The hypersurface $Y:=\{F=0\}\subset\mathbb{P}^{n+1}_k$ is smooth. 
\end{enumerate}
We have $1\iff2\Longrightarrow3$. If $d$ is invertible in $k$, they are equivalent. 
\end{lm}
\proof{
We may assume that $k$ is algebraically closed. 
The implication $2\Rightarrow1$ is obvious. For the other direction, suppose that $F$ has a singular point $x\in\A^{n+2}_k$ outside the origin. 
Since $F$ is homogeneous, the map $\A^1_{k}\to\A^{n+2}_k,t\mapsto tx$ factors through the singular locus. Hence, the origin is not 
isolated. 

By the jacobian criterion, the condition $2$ is equivalent to saying that the polynomials $\frac{\partial F}{\partial T_i}$ do not have 
a non-zero common root. On the other hand, $3$ is equivalent to saying that those polynomials together with $F$ do not have 
a non-zero common root, hence the implication $2\Rightarrow3$. The last assertion follows from Euler's identity $F=\frac{1}{d}\sum_iT_i\frac{\partial F}{\partial T_i}$. 
\qed
}

 Let $n$ be an integer $\geq0$. 
For a multi-index $I=(i_0,\dots,i_{n+1})\in 
\mathbb{N}^{n+2}$, write $T^I:=\prod_jT_j^{i_j}$ and $|I|:=\sum_ji_j$, as usual. 
Consider an affine space $S={\rm Spec}(\mathbb{Z}[\{C_I\}_I])$ with indeterminates $C_I$ indexed by the multi-indices $I$ with $|I|=d$, which 
we view as the {\it moduli space} 
of homogeneous polynomials of degree $d$ with the universal homogeneous polynomial $F=\sum_IC_IT^I$. Its divided discriminant 
${\rm disc}_d(F)\in \mathbb{Z}[\{ C_I\}]$ 
is defined in \cite[Section 2]{hypdet}. It is an irreducible homogeneous polynomial and fits into the equality  
\begin{equation}\label{resdiv}
{\rm Res}(\frac{\partial F}{\partial T_0},\dots, \frac{\partial F}{\partial T_{n+1}})=
d^{a(n,d)}{\rm disc}_d(F)
\end{equation}
where ${\rm Res}(\frac{\partial F}{\partial T_0},\dots, \frac{\partial F}{\partial T_{n+1}})$ is the resultant \cite[13.1.A]{discres} of the 
partial derivatives and 
$a(n,d)=\frac{(d-1)^{n+2}-(-1)^{n+2}}{d}$. The integer $d^{a(n,d)}$ is the greatest common divisor of the coefficients of the resultant. 

For a homogeneous polynomial $F$ of degree $d$ with coefficients in a ring $R$, we define the divided discriminant 
${\rm disc}_d(F)\in R$ to be the specialization by the map $\mathbb{Z}[\{C_I\}]\to R$ defining $F$. 
As is explained 
in \cite{hypdet}, the discriminant ${\rm disc}_d(F)\in R$ is invertible if and only if the hypersurface $Y:=\{F=0\}$ in 
$\Proj^{n+1}_R={\rm Proj}(R[T_0,\dots,T_{n+1}])$ 
is smooth. On the other hand, the resultant ${\rm Res}(\frac{\partial F}{\partial T_0},\dots, \frac{\partial F}{\partial T_{n+1}})\in R$ 
is invertible if and only if the map $F\colon \mathbb{A}^{n+2}_R\to\mathbb{A}^1_R$ is smooth outside the origin \cite[13.1.A]{discres}. 
By (\ref{resdiv}), the latter condition is equivalent to 
the former one only when $d=2$ as $a(n,d)$ is zero only in this case. When $d>2$, it also requires that 
$d$ is invertible in $R$. 
Taking this into account, the open subscheme 
$\Tilde{U}:={\rm Spec}(\mathbb{Z}[\{C_I\}_I,1/{\rm disc}_d(F)])$ is regarded as the moduli of homogeneous polynomials of degree $d$ whose hypersurfaces are smooth, whereas the open subscheme $\Tilde{U}_{{d}}:=\Tilde{U}\times_{\mathbb{Z}}\mathbb{Z}[1/d]$ 
(when $d>2$) or $\Tilde{U}$ (when $d=2$) is considered as the moduli of homogeneous polynomials with isolated singularities at 
the origin. To simplify the situation, let us focus on 
$\Tilde{U}_{2d}:=\Tilde{U}\times_{\mathbb{Z}}\mathbb{Z}[1/
2d]$; we also invert $2$ so that 
the square roots of invertible elements define \'etale $\mathbb{Z}/2$-coverings.

Set $\Tilde{U}_{{2d}}=\Tilde{U}\times_{\mathbb{Z}}\mathbb{Z}[1/2d]$ and 
$R_{{2d}}=\mathbb{Z}[\{C_I\},1/2d,1/{\rm disc}_d(F)]=\Gamma(\Tilde{U}_{{2d}},\OO_{\Tilde{U}_{{2d}}})$. 
Since $\mathbb{Z}[\{C_I\}]$ is a 
unique factorization domain 
and ${\rm disc}_d(F)$ is irreducible, 
the $\F_2$-vector space $R_{{2d}}^\times/(R_{{2d}}^\times)^2$ is spanned freely by 
the classes of $-1,{\rm disc}_d(F),\ell$ (where 
$\ell$ runs through the prime divisors of $2d$). For an element 
$a\in R_{2d}^\times$, we write $[a]$ for the quadratic character of the square root of $a$ in $H^1(\Tilde{U}_{2d},\mathbb{Z}/2)$. 
Since the map $a\mapsto[a]$ gives an isomorphism $R_{{2d}}^\times/(R_{{2d}}^\times)^2\to H^1(\Tilde{U}_{2d},\mathbb{Z}/2)$, 
 the latter has a basis $[-1],[{\rm disc}_d(F)],[\ell]$ as $\F_2$-vector space.

As the map $F\colon \mathbb{A}^{n+2}\to\mathbb{A}^1$ on $\Tilde{U}_{2d}$ is smooth outside the origin, its singular locus 
is finite over 
$\Tilde{U}_{{2d}}$. Therefore we have a non-degenerate symmetric bilinear form 
$(\varphi_F,B_{F,dt})$ on $R_{{2d}}$ and its disciriminant ${\rm disc}B_{F,dt}\in R_{{2d}}^\times/(R_{{2d}}^\times)^2
\cong H^1(\Tilde{U}_{2d},\mathbb{Z}/2)$. 
The purpose of this subsection is to express ${\rm disc}B_{F,dt}$ via the basis $[-1],[{\rm disc}_d(F)],[\ell]$. 

The quadratic 
character corresponding to $[{\rm disc}_d(F)]$ relates with the determinants of the $\ell$-adic cohomologies of middle degree 
of the hypersurfaces $Y$ when $n$ is even (\cite[Theorem 3.5]{hypdet}) and those of double covers of $\mathbb{P}^{n+1}$ ramified along $Y$ when $n$ is odd and $d$ is even (\cite[Theorem 2.3]{double}). 

We recall the result in \cite{hypdet}, as it is necessary 
for our proof. For a proper smooth $k$-variety $Z$ of even dimension $m$, 
the cup product induces a symmetric perfect pairing on $H^m(Z_{\bar{k}},\Ql(\frac{m}{2}))$. Hence its determinant character is of order at most $2$.  The same construction is also carried out on a proper smooth scheme 
over a general scheme by using $\ell$-adic sheaves for a prime $\ell$ invertible in the scheme. 

When $Z=Y$ is a hypersurface, we have 
\begin{thm}(\cite[Theorem 3.5]{hypdet})\label{hypdeet}
 
Let $Y\subset\Proj^{n+1}_{\Tilde{U}_{2d}}$ be the universal smooth hypersurface defined by the equation $F=0$ over $\Tilde{U}_{2d}$. 
 Suppose that $n$ is even. Then the quadratic character $\det R^nf_\ast\Ql(\frac{n}{2})$, where 
$f\colon Y\to\Tilde{U}_{2d}$ is the structure map and $\ell$ is a prime divisor of $2d$, is independent of 
$\ell$ and 
 is defined by the square root of $\varepsilon(n,d)\cdot{\rm disc}_d(F)$ with a sign 
$\varepsilon(n,d)=\pm1$. The sign is $(-1)^{\frac{d-1}{2}}$ when $d$ is odd and $(-1)^{\frac{d}{2}\frac{n+2}{2}}$ when $d$ is even. 
\end{thm}

The following is the main result in this subsection. 

\begin{pr}\label{dischomog}
Let $\Tilde{U}_{{2d}}={\rm Spec}(\mathbb{Z}[\{C_I\},1/2d,1/{\rm disc}_d(F)])={\rm Spec}(R_{{2d}})$. Let $F=\sum_IC_IT^I$ be the universal homogeneous polynomial 
on $\Tilde{U}_{{2d}}$ and let ${\rm disc}B_{F,dt}\in R_{{2d}}^\times/(R_{{2d}}^\times)^2\cong H^1(\Tilde{U}_{2d},\mathbb{Z}/2)$ be the discriminant of 
$(\varphi_F,B_{F,dt})$. Then we have the following equalities in $
 H^1(\Tilde{U}_{2d},\mathbb{Z}/2)$. 

\begin{enumerate}
\item 
Suppose that $d$ is even. 
\begin{enumerate}
\item When $n$ is odd, we have 
\begin{equation*}
{\rm disc}B_{F,dt}=[{\rm disc}_d(F)]+\frac{d-2}{2}[-1]+[d]. 
\end{equation*}
\item When $n$ is even, we have 
\begin{equation*}
{\rm disc}B_{F,dt}=[{\rm disc}_d(F)]. 
\end{equation*}
\end{enumerate}
\item Suppose that $d$ is odd. 
\begin{enumerate}
\item When $n$ is odd, we have 
\begin{equation*}
{\rm disc}B_{F,dt}=0.
\end{equation*}
\item When $n$ is even, we have 
\begin{equation*}
{\rm disc}B_{F,dt}=[d]+[{\rm disc}_d(F)]. 
\end{equation*}
\end{enumerate}
\end{enumerate}
\end{pr}

We record a corresponding result in characteristic $2$ as a corollary, 
in the following form. 
\begin{cor}
Let $k$ be a perfect field of characteristic $2$. Let $F\in k[T_0,\dots,T_{n+1}]$ be a homogeneous polynomial 
of odd degree $d\geq3$. Assume that $n$ is even and
 that $F\colon\A^{n+2}_{k}\to\A^1_{k}$ has an isolated singularity at the origin. Write $Y\subset\Proj^{n+1}_{k}$ for the 
 smooth hypersurface defined by $F=0$. 
 Then the quadratic character
\begin{equation}\label{sign2}
\det(H^n(Y_{\overline{k}},\Ql(\frac{n}{2})))
\end{equation}
corresponds to ${\rm Arf}(F,0)+\frac{d^2-1}{8}\in k/\wp(k)$ via Artin-Schreier theory. 
\end{cor}
\proof{
Let $W=W(k)$ be the Witt ring with $k$-coefficients and $K_0$ be its fraction field. Choose a homogeneous 
polynomial $\Tilde{F}\in W[T_0,\dots,T_{n+1}]$ which is a lift of $F$ to $W$. 
Since ${\rm Res}(\frac{\partial \Tilde{F}}{\partial T_0},\dots, \frac{\partial \Tilde{F}}{\partial T_{n+1}})\in W$ is invertible in $k$, itself is invertible in 
$W$. Thus the map 
$\Tilde{F}\colon \A^{n+2}_W\to\A^1_W$ is smooth outside the origin. Specializing the equality in  
Proposition \ref{dischomog}.2, we have  
\begin{equation}\label{qqe}
{\rm disc}B_{\Tilde{F},dt}=((-1)^{\frac{d-1}{2}}d)\cdot((-1)^{\frac{d-1}{2}}{\rm disc}_d(\Tilde{F}))
\end{equation}
in $K_0^\times/(K_0^\times)^2$. Consider the map $k/\wp(k)\to K_0^\times/(K_0^\times)^2,a\mapsto1+4[a]$. 
By Theorem \ref{hypdeet} and the proper base change theorem, $(-1)^{\frac{d-1}{2}}{\rm disc}_d(\Tilde{F})$ is the image 
of the element in $k/\wp(k)$ corresponding to $\det(H^n(Y_{\overline{k}},\Ql(\frac{n}{2})))$ (cf. Lemma \ref{Wunr}). 
Since $d^2-(-1)^{\frac{d-1}{2}}\cdot2d+1$ is a multiple of $16$ for any odd integer $d$, we have 
 $(-1)^{\frac{d-1}{2}}d\equiv1+4[\frac{d^2-1}{8}]$ in $K_0^\times/(K_0^\times)^2$. 
On the other hand, we have 
\begin{equation*}
{\rm disc}B_{\Tilde{F},dt}\equiv(-1)^N(1+4[{\rm Arf}(F,0)])
\end{equation*}
by Proposition \ref{conarf} where $N=\frac{(n+2)\mu(F,0)}{2}$. 
Note that  $\mu(F,0)$ is equal to $(d-1)^{n+2}$. This follows from Lemma \ref{discferm}.3 below and the continuity of Milnor numbers (Proposition \ref{Milconti}). 
Hence the discriminant ${\rm disc}B_{\Tilde{F},dt}$ is equal to $1+4[{\rm Arf}(F,0)]$ in $K_0\times/(K_0^\times)^2$. 
Since the map $k/\wp(k)\to K_0^\times/(K_0^\times)^2,a\mapsto1+4[a]$ is injective, the assertion follows. 
\qed}

First, considering homogeneous polynomials of Fermat type, we show that the proposition is true when $d$ is even and 
that ${\rm disc}B _{F,dt}$ is either $0$ or $[d]+[{\rm disc}_d(F)]$ when $d$ is odd. 
\begin{lm}\label{discferm}
Set $R_{a_\bullet,2d}=\mathbb{Z}[a_0^{\pm1},\dots,a_{n+1}^{\pm1},
1/2d]$. This ring admits a universal polynomial 
$F_{a_\bullet}:=\sum_ia_iT_i^d$. 
\begin{enumerate}
\item 
The divided discriminant of $F_{a_\bullet}$ is equal to 
\begin{equation*}
d^{(n+2)(d-1)^{n+1}-a(n,d)}\cdot(a_0\cdots a_{n+1})^{(d-1)^{n+1}}. 
\end{equation*}
\item The discriminant of $(\varphi_{F_{a_\bullet}},B_{F_{a_\bullet},dt})$ is equal to 
\begin{equation*}
(-1)^{\frac{(d-2)(d-1)^{n+2}(n+2)}{2}}\cdot d^{(d-1)^{n+2}(n+2)}\cdot (a_0\cdots a_{n+1})^{(d-1)^{n+2}}. 
\end{equation*}
\item The Milnor number $\mu(F_{a_\bullet},0)$ is equal to 
$(d-1)^{n+2}$. 
\end{enumerate}
\end{lm}
\proof{
$1$. 
By (\ref{resdiv}), we have 
\begin{equation*}
{\rm disc}_d(F_{a_\bullet})=d^{-a(n,d)}{\rm Res}((a_0d)T_0^{d-1},\dots, 
(a_{n+1}d)T_{n+1}^{d-1}). 
\end{equation*}
The assertion follows since the resultant is homogeneous in each  coefficient $a_i$ of degree $(d-1)^{n+1}$ and is normalized by 
${\rm Res}(T_0^{d-1},\dots,T_{n+1}^{d-1})=1$. 

$2,3$. We compute $(\varphi_{F_{a_\bullet}},B_{F_{a_\bullet},dt})$ explicitly as follows. 
For $i=0,\dots,n+1,$ set $F_i:=a_iT^d_i$. 
By Lemma \ref{prodphi}, we have an isomorphism 
\begin{equation*}
(\varphi_{F_{a_\bullet}},B_{F_{a_\bullet},dt})\cong\bigotimes_{i=0}^{n+1}(\varphi_{F_{i}},B_{F_{i},dt}). 
\end{equation*}
Hence the assertion is reduced to the case $F=aT^d$ with one indeterminate. In this case, $\varphi_{F}$ is isomorphic to 
$R_{a_\bullet,2d}[T]/(T^{d-1})$ and the bilinear form $B_{F,dt}$ is represented by $(\frac{1}{ad}\delta_{i,d-2-j})_{0\leq i,j\leq d-2}$ 
with respect to the basis $1,T,\dots,T^{d-2}$, where $\delta_{i,d-2-j}$ is $1$ if $i=d-2-j$ and $0$ otherwise. 
Its discriminant is equal to $\epsilon(ad)^{d-1}$ where $\epsilon$ is the determinant of the matrix $(\delta_{i,d-2-j})_{0\leq i,j\leq d-2}$, 
which is $(-1)^{\frac{(d-1)(d-2)}{2}}$. 
The assertion $2$ follows from this computation. The assertion $3$ 
follows as $R_{a_\bullet,2d}[T]/(T^{d-1})$ has a rank $d-1$. 
\qed}

Set 
$\Tilde{U}_{a_\bullet,2d}={\rm Spec}(R_{a_\bullet,2d})$. The universal polynomial 
$F_{a_\bullet}=\sum_ia_iT_i^d$ 
induces a map $\Tilde{U}_{a_\bullet,2d}\to\Tilde{U}_{2d}$, and a map 
\begin{equation}\label{mapferm}
H^1(\Tilde{U}_{2d},\mathbb{Z}/2)\to H^1(\Tilde{U}_{a_\bullet,2d},\mathbb{Z}/2). 
\end{equation}
The quadratic characters 
$[-1],[a_0],\dots,[a_{n+1}],[\ell]$ (for $\ell\mid2d$) 
 defined by the square roots of $-1,a_i,\ell\in R_{a_\bullet,2d}^\times$ form a basis of the $\F_2$-vector space $H^1(\Tilde{U}_{a_\bullet,2d},\mathbb{Z}/2)$.  Lemma \ref{discferm} gives the following. 
 
 \begin{lm}\label{lmhomog}
 Let the situation be as in Proposition \ref{dischomog}. 
 \begin{enumerate}
 \item When $d$ is even, the assertion $1$ in Proposition \ref{dischomog} holds. 
 \item When $d$ is odd, ${\rm disc}B_{F,dt}$ is equal to either $0$ or $[d]+[{\rm disc}_d(F)]$. 
 \end{enumerate}
 \end{lm}
\proof{
In the sequel, we frequently use Lemmas \ref{discferm}.1,2 to compute the images of ${\rm disc}B_{F,dt},
[{\rm disc}_d(F)]$ by the map (\ref{mapferm}). 

$1$. 
Suppose that $d$ is even. The image of $[{\rm disc}_d(F)]$ by (\ref{mapferm}) is 
$\sum_i[a_i]$, as $(n+2)(d-1)^{n+1}-a(n,d)$ is even. 
Consequently, the map (\ref{mapferm}) is injective in this case. Hence, 
it suffices to show the equalities in $(a)$ and $(b)$ as characters 
in $H^1(\Tilde{U}_{a_\bullet,2d},\mathbb{Z}/2)$. 
Then the assertion follows from Lemma \ref{discferm}.2. 

$2$. When $d$ is odd,  the map (\ref{mapferm}) sends $[{\rm disc}_d(F)]$ to 
$[d]$. Hence its kernel 
consists of the two elements $0$ and $[d]+[{\rm disc}_d(F)]$. 
As ${\rm disc}B_{F,dt}$ is contained in the kernel, it is equal to 
$0$ or $[d]+[{\rm disc}_d(F)]$ in $H^1(\Tilde{U}_{2d},\mathbb{Z}/2)$ in this case. 
\qed}

To determine ${\rm disc}B_{F,dt}$ when $d$ is odd, we use the following lemma. 
\begin{lm}\label{epF}
Assume that $d$ is odd. 
Let $k$ be a finite field whose characteristic is prime to $2d$. We assume that $k$ contains a primitive $d$-th power root of unity. 
Then we have the following equalities. 
\begin{enumerate}
\item When $n$ is odd, we have 
\begin{equation*}
\bigl(\frac{{\rm disc}B_{F,dt}}{k}\bigr)=1. 
\end{equation*}
\item When $n$ is even, we have 
\begin {equation*}
\bigl(\frac{{\rm disc}B_{F,dt}}{k}\bigr)=\det H^n(Y_{\bar{k}},\Ql(\frac{n}{2}))({\rm Frob}_k). 
\end{equation*}
\end{enumerate}
\end{lm}

Let us postpone the proof of this lemma and complete the proof of Proposition \ref{dischomog}. 

\vspace{2mm}

(Proof of Proposition \ref{dischomog}, admitting Lemma \ref{epF})

\vspace{2mm}
{The assertion $1$ in the proposition is proved in Lemma \ref{lmhomog}. By the same lemma, we already know that 
${\rm disc}B_{F,dt}=0$ or $[d]+[{\rm disc}_d(F)]$ when $d$ is odd. For a finite field $k$ with the characteristic prime to $2d$, 
set $\Tilde{U}_k:=\Tilde{U}\times_{\mathbb{Z}}k$. As the divided discriminant ${\rm disc}_d(F)\in k[\{C_I\}]$ with $k$-coefficient is 
geometrically irreducible (\cite[Proposition 2.12]{hypdet}), the map $H^1(\Tilde{U}_{2d},\mathbb{Z}/2)\to H^1(\Tilde{U}_{k},\mathbb{Z}/2)$ is injective on the subspace 
$\{0,[d]+[{\rm disc}_d(F)]\}$. Hence it is enough to determine the image of ${\rm disc}B_{F,dt}$. 

In the sequel of the proof, We choose $k$ which contains a primitive $d$-th power root of unity and also 
 the square roots of $-1,d$. For such a $k$, 
we have equalities in $H^1(\Tilde{U}_k,\mathbb{Z}/2)$
\begin{equation*}
[d]+[{\rm disc}_d(F)]=[{\rm disc}_d(F)]=\det Rf_\ast\Ql(\frac{n}{2}), 
\end{equation*}
where $f\colon Y\to \Tilde{U}_k$ is the structure map of the universal hypersurface. 
The second equality is the main result of \cite{hypdet}, which we recall in Theorem \ref{hypdeet}. 

By the Chebotarev density, the assertion $2$ in Proposition \ref{dischomog} is reduced to the following: for 
any finite extension $k'/k$ and any $k'$-valued point $z\in \Tilde{U}_k(k')$, 
\begin{enumerate}
\item[$(a)'$] when $n$ is odd, we have $\bigl(\frac{{\rm disc}B_{F_z,dt}}{k'}\bigr)=1$. 
\item[$(b)'$] when $n$ is even, we have $\bigl(\frac{{\rm disc}B_{F_z,dt}}{k'}\bigr)=\det H^n(Y_{z,\bar{k}},\Ql(\frac{n}{2}))({\rm Frob}_{k'})$. 
\end{enumerate}
Here $F_z$ is the homogeneous polynomial defined by the point $z$ and 
$Y_z$ is the hypersurface defined by $F_z$. Hence Lemma \ref{epF} completes the proof. 
\qed}

The rest of this subsection is devoted to the proof of Lemma \ref{epF}. 
For this purpose, we compute the local epsilon factors of the vanishing cycles complexes. The 
following method by blowing up is suggested to the author by T. Saito. 

Let $n\geq0,d\geq2$ be integers and $k$ be a perfect field whose characteristic is prime to $2d$. For a homogeneous polynomial 
$F=\sum_{|I|=d}C_IT^I\in k[T_0,\dots,T_{n+1}]$ satisfying the equivalent conditions in Lemma \ref{condhom}, we describe the vanishing 
cycles complex $R\Phi_F(\mathbb{Q}_\ell)_0$ of $F\colon\mathbb{A}^{n+2}_k\to\A^1_k$ as follows. 
Let $S={\rm Spec}(k\{t\})$ be the henselization of $\A^1_k$ at the origin, whose closed and generic points are denoted by $s$ and $\eta$. 
Let $X$ be the base change $\A^{n+2}_{k}\times_{F,\A^1_k}S$. We 
use the same symbol $F\colon X\to S$ for the base change by abuse of notation. 
Let $\Tilde{X}\to X$ be the 
blow-up at the origin and let $\Tilde{F}$ be the structure map $\Tilde{X}\to S$. 
We write $E$ for the exceptional divisor of the blow-up and $D_F$ for the strict transform of the special fiber $X_s$ in $\Tilde{X}$. 
Using the coordinate functions $T_0,\dots,T_{n+1}$, we identify $E$ with $\mathbb{P}^{n+1}_k$ and $E\cap D_F$ with the 
hypersurface $Y=\{F=0\}\subset \mathbb{P}^{n+1}_k=E$. 

Let $S':={\rm Spec}(\OO_S[u]/(u^d-t))$ be a tamely ramified extension of $S$. 
We define $\Tilde{X}'$ to be the normalization of $\Tilde{X}\times_SS'$. 
Let $E'$ be the pullback of $E\subset \Tilde{X}$ to $\Tilde{X}'$. For each $0\leq i\leq n+1$, $\Tilde{X}$ has the open 
 locus where we have $T_i\cdot\OO_{\Tilde{X}}=\OO_{\Tilde{X}}(-E)$.  Over this open subset, $\Tilde{X}'$ is defined by the $d$-th power root 
 $\sqrt[d]{F(T_0/T_i,\dots,T_{n+1}/T_i)}$. Therefore, the finite coverings $\Tilde{X}'\to\Tilde{X}$ and $E'\to E$ are tamely ramified along 
 $D_F$ and $Y=E\cap D_F$ respectively and finite \'etale elsewhere. Moreover, the inclusion $D_F\hookrightarrow \Tilde{X}$ lifts 
 to an inclusion $D_F\hookrightarrow\Tilde{X}'$, which defines a smooth divisor $D_{F}'$ of $\Tilde{X}'$. 
 The divisor $E'\cup D_F'$ is simple normal crossings in $\Tilde{X}'$ with $E'\cap D_F'$ isomorphic to $Y$. 
 
 Write $s'$ for the closed point of $S'$. Since the special fiber $\Tilde{X}'_{s'}$ is equal to $E'\cup D_F'$ as divisors, 
 $\Tilde{X}'$ is strictly semi-stable over $S'$, by which we mean that, Zariski locally 
 on $\Tilde{X}'$, it is \'etale over ${\rm Spec}(\OO_{S'}[t_1,\dots,t_m]/(t_1\cdots t_r-\pi))$ for a uniformizer $\pi$ of $\OO_{S'}$ and integers $1\leq r\leq m$. 
 For a strictly semi-stable scheme over a trait, its nearby cycles complex is computed in \cite{WSS}, which we recall for $\Tilde{X}'$ in Lemma \ref{nctilde}. 
 
 To proceed, we give necessary definitions and notations. 
 Let $\bar{k}$ be an algebraic closure of $k$ with its affine spectrum denoted by $\bar{s}$. Write $S^{\rm ur}$ for the 
 strict henselization of $S$ at $\bar{s}\to S$. 
 Write $S'^{{\rm ur}}:=S'\times_SS^{\rm ur}$ and let $\bar{s}'$ be its closed point defined by the diagonal map $\bar{s}\to
 S'\times_SS^{\rm ur}$. We write $\Tilde{X}'_{\bar{s}'}=\Tilde{X}'\times_{S'}\bar{s}'$ for the geometric special fiber of $\Tilde{X}'$ 
 over $S'$. 
 Let $\overline{\eta}$ be the spectrum of a separable closure of the function field of $S'^{\rm ur}$, where 
 $k(\overline{\eta})$ is also regarded as a separable closure of the function field $k(\eta)$ of $S$.  
 We have canonical morphisms
 \begin{equation}\label{morX'}
\xymatrix{
 \Tilde{X}'\times_{S'}\overline{\eta}\ar[r]^-{\bar{\tilde{j}}'}& \Tilde{X}'\times_{S'}S'^{\rm ur}&
  \Tilde{X}'_{\bar{s}'}\ar[l]_-{\tilde{i}'}. 
  }
  \end{equation}
  The Galois group ${\rm Gal}(\overline{\eta}/\eta)$ acts on 
 $\Tilde{X}'\times_{S'}S'^{\rm ur}$ as follows. It acts on $\Tilde{X}\times_SS'^{\rm ur}$ via the right component of the product, 
 hence acts on its normalization $\Tilde{X}'\times_{S'}S'^{\rm ur}$. 
 Under the isomorphisms $\Tilde{X}'\times_{S'}\overline{\eta}\cong (\Tilde{X}'\times_{S'}S'^{\rm ur})\times_{S'^{\rm ur}}\overline{\eta}, \Tilde{X}'_{\bar{s}'}\cong (\Tilde{X}'\times_{S'}S'^{\rm ur})\times_{S'^{\rm ur}}\bar{s}'$, these schemes also admit Galois actions, by which 
 (\ref{morX'}) are ${\rm Gal}(\overline{\eta}/\eta)$-equivariant. Let us write $\Tilde{j}'\colon (\Tilde{X}'\times_{S'}S'^{\rm ur})\setminus\Tilde{X}'_{\bar{s}'} \to\Tilde{X}'
 \times_{S'}S'^{\rm ur}$ for the open complement of $\tilde{i}'$.

 For an immersion $\iota\colon Z\to T$ of schemes, $\mathbb{Q}_{\ell, Z}$ denotes the $0$-extension $\iota_!\mathbb{Q}_\ell$. 
 Set $\mQ_{\ell,E'_{\bar{k}},D'_{F,\bar{k}}}:=\mQ_{\ell,E'_{\bar{k}}}\oplus\mQ_{\ell,D'_{F,\bar{k}}}$, which 
 is an $\ell$-adic sheaf on $\Tilde{X}'_{\bar{s}'}$.  
 Let $\theta'_1\colon
 \mQ_{\ell,E'_{\bar{k}},D'_{F,\bar{k}}}\to \tilde{i}'^\ast R^1\tilde{j}'_{\ast}\mQ_\ell(1)$ be the map defined by 
 the divisors $E'_{\bar{k}},D'_{F,\bar{k}}$ (\cite[1.1]{WSS}). 
  By cup product, it induces a map $\theta'_m\colon\wedge^{m}\mQ_{\ell,E'_{\bar{k}},D'_{F,\bar{k}}}\to \tilde{i}'^\ast R^{m}\tilde{j}'_{\ast}\mQ_\ell(m)$ for $m\geq1$. 
 
 \begin{lm}\label{nctilde}
 Let $R\psi\mathbb{Q}_\ell=\tilde{i}'^\ast R\bar{\tilde{j}}'_{\ast}\mathbb{Q}_\ell$ be the nearby cycles complex of $\Tilde{X}'$ 
 as an $S'$-scheme. 
 \begin{enumerate}
 \item (\cite[Proposition 1.1.2.1]{WSS}) The map $\theta'_m$ is an isomorphism for $m\geq1$. 
 \item (\cite[Proposition 1.1.2.2]{WSS}) Consider the composition of 
 \begin{equation*}
 \theta\colon\mQ_{\ell,\Tilde{X}'_{\bar{s}'}}\xrightarrow{\delta}\mQ_{\ell,E'_{\bar{k}},D'_{F,\bar{k}}}
\xrightarrow{\theta'_1} \tilde{i}'^{\ast}R^1\tilde{j}'_{\ast}\Ql(1) 
\end{equation*}
where $\delta$ is the canonical map. 
Then, for $m\geq0$, the cup product 
$\theta\cup\colon\tilde{i}'^{\ast}R^m\tilde{j}'_\ast\Ql(m)\to\tilde{i}'^{\ast}R^{m+1}\tilde{j}'_\ast\Ql(m+1)$ factors as 
$\tilde{i}'^{\ast}R^m\tilde{j}'_\ast\Ql(m)\to R^m\psi\Ql(m) \to\tilde{i}'^{\ast}R^{m+1}\tilde{j}'_\ast\Ql(m+1)$ 
where the former map is the canonical one, 
and this factorization identifies $R\psi^m\Ql(m)$ with the image of the map $\theta\cup$. 
  \item We have 
  $\mQ_{\ell,\Tilde{X}'_{\bar{s}'}}\cong R^0\psi\Ql$, $\mQ_{\ell,E'_{\bar{k}}\cap D'_{F,\bar{k}}}(-1)\cong 
  R^1\psi\Ql$, and $R^m\psi\Ql=0$ for $m\neq0,1$ as sheaves with equivariant ${\rm Gal}(\overline{\eta}/\eta)$-actions. 
 \end{enumerate}
 \end{lm}
\proof{
The assertions $1,2$ are proved in \cite[Proposition 1.1.2]{WSS}. Then the assertion $3$ readily follows from the 
commutativity of the diagram 
\begin{equation*}
\xymatrix{
\wedge^m\mQ_{\ell,E'_{\bar{k}},D'_{F,\bar{k}}}\ar[r]^-{\delta\wedge}\ar[d]^-{\theta'_{m}}&\wedge^{m+1}\mQ_{\ell,E'_{\bar{k}},D'_{F,\bar{k}}}
\ar[d]^-{\theta'_{m+1}}\\
\tilde{i}'^\ast R^m\tilde{j}'_{\ast}\mQ_\ell(m)\ar[r]^-{\theta\cup}&\tilde{i}'^\ast R^{m+1}\tilde{j}'_{\ast}\mQ_\ell(m+1), 
}
\end{equation*}
where $\theta'_{m}$ for $m=0$ means the canonical isomorphism $\mQ_{\ell,\Tilde{X}'_{\bar{s}'}}\to 
\tilde{i}'^\ast \tilde{j}'_{\ast}\mQ_\ell$. 
\qed}

 Let $R\Phi_F(\Ql)_0$ be the vanishing cycles complex of $X$ with respect to $F$, supported on the origin. 
 Using Lemma \ref{nctilde}, we describe $R\Phi_F(\Ql)_0$ as follows. The morphisms (\ref{morX'}) fit into the 
 commutative diagram 
 \begin{equation}\label{cartnc}
 \xymatrix{
 \Tilde{X}'\times_{S'}\overline{\eta}\ar[r]^-{\bar{\tilde{j}}'}\ar[d]^-\cong& \Tilde{X}'\times_{S'}S'^{\rm ur}\ar[d]&
  \Tilde{X}'_{\bar{s}'}\ar[d]^-{\alpha}\ar[l]_-{\tilde{i}'}\\
  \Tilde{X}\times_S\overline{\eta}\ar[d]^-\cong\ar[r]^-{\bar{\tilde{j}}}& \Tilde{X}\times_SS^{\rm ur}\ar[d]& \Tilde{X}_{\bar{s}}\ar[d]^-{\pi}\ar[l]_-{\tilde{i}}\\
  X\times_S\overline{\eta}\ar[r]^-{\bar{j}}&X\times_SS^{\rm ur}&X_{\bar{s}},\ar[l]_-i
 }
 \end{equation}
 on which ${\rm Gal}(\overline{\eta}/\eta)$ acts equivariantly. Here the left vertical arrows are isomorphisms 
 since the blow-up has its center in the special fiber and $\Tilde{X}\times_SS'$ is regular on its generic fiber.

\begin{lm}\label{Rphi_F}
In the Grothendieck group of 
 $\ell$-adic representations of ${\rm Gal}(\overline{\eta}/\eta)$, we have an equality 
 \begin{equation}\label{eqphi}
 [R\Phi_F(\Ql)_0]=[R\Gamma(E'_{\bar{k}},\Ql)]-[R\Gamma(Y_{\bar{k}},\Ql)(-1)]-[\Ql]. 
 \end{equation}
 Here the Galois action on the cohomology groups in the right-hand side is given as follows. 
 \begin{enumerate}
 
 \item Let $\mu_d$ be the group of $d$-th power roots of unity in $\bar{k}$. 
The inertia subgroup acts on $R\Gamma(E'_{\bar{k}},\Ql)$ through the quotient $\mu_d$, 
 whose action is induced from the geometric action on $E'_{\bar{k}}$ 
 given by $d$-th power roots of $F(T_0/T_i,\dots,T_{n+1}/T_i)$. 
 Let $\eta'$ be the generic point of $S'$. 
 The action of the subgroup 
 ${\rm Gal}(\overline{\eta}/\eta')$ is given by the composition of the 
 unramified quotient ${\rm Gal}(\overline{\eta}/\eta')\to 
 {\rm Gal}(\bar{k}/k)$ and the canonical action of ${\rm Gal}(\bar{k}/k)$ on $R\Gamma(E'_{\bar{k}},\Ql)$.  
  \item The Galois group ${\rm Gal}(\overline{\eta}/\eta)$ acts on $R\Gamma(Y_{\bar{k}},\Ql)(-1)$ through the unramified quotient ${\rm Gal}(\bar{k}/k)$. 
 \end{enumerate}
\end{lm}
\proof{
First, we decompose $\alpha_{\ast}R\psi\Ql\cong \tilde{i}^\ast R\bar{\tilde{j}}_\ast\Ql=:R\Psi_{\Tilde{F}}(\Ql)$ in the Grothendieck group of $\ell$-adic sheaves on $\Tilde{X}_{\bar{s}}$ with equivariant ${\rm Gal}(\overline{\eta}/\eta)$-actions. 
By Lemma \ref{nctilde}, we have $[\alpha_{\ast}R^0\psi\Ql]=[\alpha_{\ast}\mQ_{\ell,E'_{\bar{k}}}]+[\mQ_{\ell,D_{F,\bar{k}}}]-[\mQ_{\ell,
Y_{\bar{k}}}]$, $[\alpha_{\ast}R^1\psi\Ql]=[\mQ_{\ell,Y_{\bar{k}}}(-1)]$, and $\alpha_{\ast}R^i\psi\Ql=0$ for the other degrees. 
Hence we have 
\begin{equation}\label{rphi}
[R\Psi_{\Tilde{F}}(\Ql)]=[\alpha_{\ast}R\psi(\Ql)]=
[\alpha_{\ast}\mQ_{\ell,E'_{\bar{k}}}]+[\mQ_{\ell,D_{F,\bar{k}}\setminus Y_{\bar{k}}}]
-[\mQ_{\ell,Y_{\bar{k}}}(-1)]. 
\end{equation}

  The right squares in (\ref{cartnc}) are cartesian up to nilpotent thickening. 
Hence the proper base change theorem gives us $R\pi_{\ast}R\Psi_{\Tilde{F}}(\Ql)\cong R\Psi_F(\Ql)$. 
 Noting the isomorphism $R\pi_{\ast}\mQ_{\ell,D_{F,\bar{k}}\setminus Y_{\bar{k}}}\cong \mQ_{\ell,F^{-1}(0)\setminus\{0\}}
$, we have 
\begin{equation*}
[R\Phi_F(\Ql)_0]=[R\Psi_{F}(\Ql)]-[\mQ_{\ell, F^{-1}(0)}]=
[R(\pi\alpha)_{\ast}\mQ_{\ell,E'_{\bar{k}}}]-[R\pi_\ast\mQ_{\ell,Y_{\bar{k}}}(-1)]-[\mQ_{\ell,0}]
\end{equation*}
in the Grothendieck group of $\ell$-adic sheaves on $X_{\bar{s}}$. Taking the stalk at the origin, we obtain the desired equality. 
 The Galois actions on $E'_{\bar{k}},Y_{\bar{k}}$ are induced 
from the Galois action on (\ref{cartnc}). Hence the assertions $1,2$ are verified. 
\qed}

To compute the local epsilon factor of $R\Gamma(E'_{\bar{k}},\Ql)$, we need to 
recall well-known results on the cohomology groups of a variety with a finite group action. 
\begin{lm}\label{detvar}
Let $k$ be a finite field with $q$ elements. For elements $a,b\in\overline{\Ql}^\times$, let us write $a\equiv b$ if 
$ab^{-1}\in q^{\mathbb{Z}}$. 

Let $Z$ be a proper smooth $k$-variety of dimension $m$ with an admissible action of a finite group $G$. Write 
$Z_{\bar{k}}/G$ for the quotient scheme of $Z_{\bar{k}}$ by $G$. 
For an irreducible $\overline{\Ql}$-representation $\rho$ of $G$, let $H^i_\rho$ denote the $\rho$-isotypic component of $H^i(Z_{\bar{k}},\overline{\Ql})$. 
\begin{enumerate}
\item The cup product induces an isomorphism  
$H^{2m-i}_{\rho^{\vee}}\cong(H^i_\rho)^\vee(-m)$ of Galois representations. Consequently, we have 
\begin{equation*}
\det(\pm{\rm Frob}_k,H^i_\rho)\cdot\det(\pm{\rm Frob}_k,H^{2m-i}_{\rho^{\vee}})=q^{mb_{\rho,i}}, 
\end{equation*}
where the signs are taken to be equal to each other and we set $b_{\rho,i}={\rm rk}H^i_\rho={\rm rk}H^{2m-i}_{\rho^{\vee}}$. 
\item Suppose that an irreducible representation $\rho$ is self-dual. Then we have 
\begin{equation*}
\det({\rm Frob}_k,H^m_\rho)=\pm q^{\frac{mb_{\rho,m}}{2}}. 
\end{equation*}
Further if $m$ is odd, then $b_{\rho,m}$ is even and the sign is $1$. 
\item For the trivial representation $\rho=1$, we have an isomorphism $H^i_1\cong H^i(Z_{\bar{k}}/G,\overline{\Ql})$ of 
Galois representations. 
\item
If $G$ is abelian of odd order, we have 
\begin{equation*}
\det({\rm Frob}_k,R\Gamma(Z_{\bar{k}},\Ql))\equiv\det({\rm Frob}_k,H^m(Z_{\bar{k}}/G,\Ql))^{(-1)^m}. 
\end{equation*}
\end{enumerate}
\end{lm}
\proof{
The pairing 
\begin{equation}\label{perh}
H^i_{\rho}\times H^{2m-i}_{\rho^{\vee}}\xrightarrow{\cup}H^{2m}(Z_{\bar{k}},\Ql)\xrightarrow{{\rm Tr}}\Ql(-m)
\end{equation}
is perfect by the Poincar\'e duality and the fact that the cup product is $G$-equivariant, which shows the assertion $1$. 
Applying $1$ to a self-dual $\rho$, the equality in $2$ follows. 
When $m$ is odd and $\rho$ is self-dual, (\ref{perh}) equips $H^m_\rho$ with an alternating perfect pairing, which shows 
the remaining part of $2$. 

 Let $f\colon Z\to Z/G$ be the quotient map. 
The isomorphism $R\Gamma(G,R\Gamma(Z_{\bar{k}}/G,-))\cong R\Gamma(Z_{\bar{k}}/G,R\Gamma(G,-))$ of derived functors gives 
\begin{equation*}
R\Gamma(Z_{\bar{k}},\Ql)^G\cong (R\Gamma(Z_{\bar{k}}/G,f_\ast\Ql))^G\cong R\Gamma(Z_{\bar{k}}/G,(f_\ast\Ql)^G). 
\end{equation*}
As the canonical map $\Ql\to(f_\ast\Ql)^G$ is an isomorphism, which can be verified stalk-wise, the assertion $3$ follows. 

We show $4$. The assumption on $G$ implies that any irreducible representation of $G$ (which is one dimensional as $G$ is abelian) 
is not self-dual unless it is trivial. 
Hence by $1$, we have 
\begin{equation*}
\det({\rm Frob}_k,R\Gamma(Z_{\bar{k}},\Ql))\equiv\det({\rm Frob}_k,H^m_1)^{(-1)^m}
\end{equation*}
in $\overline{\Ql}^\times/q^{\mathbb{Z}}$. 
Then the assertion $4$ follows from $3$. 
\qed}

Let us go back to our situation. Let $k$ be a finite field in which $2d$ is invertible. Let $F\colon \A^{n+2}_k\to\A^1_k$ be a 
homogeneous function of degree $d$ which is smooth outside the origin. 
We collect necessary results on $Y,E'$ as follows. 
\begin{lm}\label{E'}
Assume that $d$ is odd. 
\begin{enumerate}
\item The integer $\chi(\Proj^{n+1}_{\bar{k}})+\chi(Y_{\bar{k}})$ is odd. 
\item Let $R\Gamma(E'_{\bar{k}},\Ql)$ be the complex of tamely ramified 
${\rm Gal}(\overline{\eta}/\eta)$-representations considered in Lemma \ref{Rphi_F}. Further assume that $k$ contains $\mu_d$. We have 
\begin{equation*}
\varepsilon_0(\A^1_{k,(0)},R\Gamma(E'_{\bar{k}},\Ql),dt)\equiv (-1)^{\chi(\Proj^{n+1}_{\bar{k}})}. 
\end{equation*}
Here we write $a\equiv b$ for $a,b\in\overline{\Ql}^\times$ if $ab^{-1}\in q^{\mathbb{Z}}$. 
\end{enumerate}
\end{lm}
\proof{
$1$. We write $T_X$ for the tangent bundle of a smooth variety $X$. 
Let $h:=c_1(\OO_{\Proj^{n+1}_k}(1))\in {\rm CH}^1(\Proj^{n+1}_k)$ be the first Chern class of the very ample sheaf. 
Let $\iota\colon Y\to\Proj^{n+1}_k$ be the closed immersion. 
The canonical exact sequence
\begin{equation*}
0\to T_Y\to\iota^\ast T_{\Proj^{n+1}_k}\to \iota^\ast\OO_{\Proj^{n+1}_k}(d)\to0
\end{equation*}
and the projection formula give a formula for the total Chern class $\iota_\ast c(T_Y)$ (\cite[3.2.12]{Ful}): 
\begin{equation*}
\iota_\ast c(T_Y)=((1+h)^{n+2}/(1+dh))\cap[Y]=(1+h)^{n+2}dh/(1+dh)
\end{equation*}
in the Chow ring of $\Proj^{n+1}_k$. As $d$ is odd, 
taking modulo $2$ gives $\iota_\ast c(T_Y)\equiv (1+h)^{n+1}h$, hence $\chi(Y_{\bar{k}})\equiv n+1$. On the other hand, 
we have $\chi(\Proj^{n+1}_{\bar{k}})=n+2$. The assertion $1$ follows. 

$2$. 
Note that 
the action of $\mu_d$ on $E'_{\bar{k}}$ is defined over $k$ by the assumption $\mu_d\subset k$. Hence we can apply Lemma \ref{detvar} to $(Z,G)=(E',\mu_d)$. 

By \cite[(5.4)]{Del}, we have 
\begin{equation*}
\varepsilon_0(\A^1_{k,(0)},R\Gamma(E'_{\bar{k}},\Ql),dt)\equiv\det(
R\Gamma(E'_{\bar{k}},\Ql))(t)\cdot\varepsilon_0(\A^1_{k,(0)},R\Gamma(E'_{\bar{k}},\Ql),\frac{dt}{t}). 
\end{equation*}
Here the determinant in the product is considered as a character of the function field of $S$ by the local class field theory. 
Note that the norm map $\OO_{S'}\to\OO_S$ sends 
$u=\sqrt[d]{t}$ to $t$. 
Hence the functoriality of the local class field theory and Lemma 
\ref{Rphi_F}.1 show that $\det(
R\Gamma(E'_{\bar{k}},\Ql))(t)$ is equal to $
\det({\rm Frob}_k,R\Gamma(E'_{\bar{k}},\Ql))$, which 
is equal to $\det({\rm Frob}_k,H^{n+1}(\Proj^{n+1}_{\bar{k}},\Ql))^{(-1)^{n+1}}
\equiv 1$ by Lemma \ref{detvar}.3. 

It remains to show $\varepsilon_0(\A^1_{k,(0)},R\Gamma(E'_{\bar{k}},\Ql),\frac{dt}{t})\equiv(-1)^{\chi(\Proj^{n+1}_{\bar{k}})}$. 
For a character $\chi\colon\mu_d\to\overline{\Ql}^\times$, 
set $V_\chi=(R\Gamma(E'_{\bar{k}},\overline{\Ql})\otimes\chi^{-1})^{\mu_d}$. By \cite[5.10]{Del}, we have 
\begin{equation*}
\varepsilon_0(\A^1_{k,(0)},R\Gamma(E'_{\bar{k}},\Ql),\frac{dt}{t})\equiv 
\prod_\chi(-\tau(\chi_k,\psi))^{{\rm rk}V_\chi}. 
\end{equation*}
Here we set $\chi_k(a)=\chi(a^{\frac{q-1}{d}})$ for $a\in k^\times$ and set $\tau(\chi_k,\psi)=-\sum_{a\in k^\times}
\chi_k(a)^{-1}\psi(a)$. For $\chi\neq1$, 
we have \cite[5.7]{Del} 
\begin{equation*}
\tau(\chi_k,\psi)\cdot\tau(\chi_k^{-1},\psi)=\chi_k(-1)q=q. 
\end{equation*}
The equality $\chi_k(-1)=1$ holds as $\chi$ is of odd order. Since we have ${\rm rk}V_\chi={\rm rk}V_{\chi^{-1}}$ and 
${\rm rk}V_1=\chi(\Proj^{n+1}_{\bar{k}})$ by Lemma \ref{detvar}, we have 
\begin{equation*}
\varepsilon_0(\A^1_{k,(0)},R\Gamma(E'_{\bar{k}},\Ql),\frac{dt}{t})\equiv (-\tau(1_k,\psi))^{{\rm rk}V_1}=(-1)^{\chi(\Proj^{n+1}_{\bar{k}})}. 
\end{equation*}
The assertion follows. 
\qed}

\vspace{2mm}

(Proof of Lemma \ref{epF}) 

\vspace{2mm}
Since there only appear signs in the lemma, we can work in $\overline{\Ql}^\times/q^{\mathbb{Z}}$, rather than $\overline{\Ql}^\times$. As in Lemmas \ref{detvar}, \ref{E'}, we write $a\equiv b$ for $a,b\in\overline{\Ql}^\times$ if $ab^{-1}\in q^{\mathbb{Z}}$.

{From the continuity of Milnor numbers (Proposition \ref{Milconti}) and Lemma \ref{discferm}.3, the Milnor number $\mu(F,0)$ is equal to $(d-1)^{n+2}$, which 
is also equal to ${\rm dimtot}R\Phi_F(\Ql)_0$ up to sign, by the Milnor formula (\ref{intromil}). Hence they are multiples of $4$ as $d$ is odd. 
Then 
Corollary \ref{Milcor}.1 takes the following form in this case: 
\begin{equation*}
\varepsilon_0(\A^1_{k,(0)},R\Phi_F(\Ql)_0,dt)=q^{(-1)^{n+3}\frac{(n+2)\mu(F,0)}{2}}\cdot\bigl(\frac{{\rm disc}B_{F,dt}}{k}\bigr), 
\end{equation*}
where we use $\tau_\psi^2=(\frac{-1}{k})q$ and ${\rm disc}B_{F,-2dt}=(-2)^{(n+2)\mu(F,0)}{\rm disc}B_{F,dt}={\rm disc}B_{F,dt}$ as $\mu(F,0)$ is even. 
We compute the local epsilon factor modulo $q^{\mathbb{Z}}$. 
By Lemma \ref{Rphi_F} and the multiplicativity of local epsilon factor, $\varepsilon_0(\A^1_{k,(0)},R\Phi_F(\Ql)_0,dt)$ is equal to 
\begin{equation*}
\varepsilon_0(\A^1_{k,(0)},R\Gamma(E'_{\bar{k}},\Ql),dt)\cdot\det(-{\rm Frob}_k,R\Gamma(Y_{\bar{k}},\Ql)(-1))^{-1}\cdot
(-1),
\end{equation*}
where we use, for an unramified representation $V$, $\varepsilon_0(\A^1_{k,(0)},V,dt)=\det(-{\rm Frob}_k,V)^{-1}$. By Lemma \ref{E'}.2, and 
the fact that the cohomology groups of projective spaces are either $0$ or Tate twists, 
we proceed 
\begin{align*}
\varepsilon_0(\A^1_{k,(0)},R\Phi_F(\Ql)_0,dt)&\equiv(-1)^{\chi(\Proj^{n+1}_{\bar{k}})+\chi(Y_{\bar{k}})+1}
\cdot\det({\rm Frob}_k,H^n(Y_{\bar{k}},\Ql))^{-1}\\
&=\det({\rm Frob}_k,H^n(Y_{\bar{k}},\Ql))^{-1}, 
\end{align*}
as $\chi(\Proj^{n+1}_{\bar{k}})+\chi(Y_{\bar{k}})+1$ is even (Lemma \ref{E'}.1). 
Hence the assertion follows by applying Lemmas \ref{detvar}.1,2 to $(Z,G)=(Y,\{1\})$. 
\qed}

The proof of Proposition \ref{dischomog} is completed. 

  \section*{Acknowledgement}
The author would like to thank his advisor Professor Takeshi Saito for many useful advices and his careful reading 
of the manuscript. 
The author also thanks Professor Tomoyuki Abe, who asked him a question which led him to this work. 
This work was supported by the Program for Leading Graduate Schools, MEXT, Japan and by 
JSPS KAKENHI Grant Number 19J11213.

Daichi Takeuchi\\
Department of Mathematics,
Institute of Science Tokyo,
2-12-1, Ookayama, Meguro-ku, Tokyo, 152-8551, Japan. \\
daichi.takeuchi4@gmail.com \\

\end{document}